\documentclass[11pt,a4,paper]{article}
\usepackage{amssymb}
\usepackage{bbm}
\newcommand{\nc}[2]{\newcommand{#1}{#2}}

\newcommand{\Boxneu}{\square}
\newcommand{\C}{\mathbbm{C}}
\newcommand{\N}{\mathbbm{N}}
\newcommand{\R}{\mathbbm{R}}
\newcommand{\Z}{\mathbbm{Z}}
\nc{\bsa}{\begin{satz}}
\nc{\bpr}{\begin{pr}}
\nc{\bth}{\begin{them}}
\nc{\ble}{\begin{lem}}
\nc{\bco}{\begin{corollary}}
\nc{\bre}{\begin{remark}}
\nc{\bex}{\begin{example}}
\nc{\bde}{\begin{de}}
\nc{\ede}{\end{de}}
\nc{\esa}{\end{satz}}
\nc{\epr}{\end{pr}}
\nc{\ethe}{\end{them}}
\nc{\ele}{\end{lem}}
\nc{\eco}{\end{corollary}}
\nc{\ere}{\hfill\mbox{$\Diamond$}\end{remark}}
\nc{\eex}{\end{example}}
\nc{\beq}{\begin{equation}}
\nc{\eeq}{\end{equation}}
\nc{\ot}{\otimes}
\nc{\lra}{\longrightarrow}
\nc{\ci}{\circ}
\nc{\hG}{\Gamma}
\begin{document}
\title{Connections on locally trivial quantum principal fibre bundles }
\author{Dirk Calow$^1$\thanks{supported by Deutsche Forschungsgemeinschaft,
e-mail
Dirk.Calow@itp.uni-leipzig.de}
 \ and Rainer Matthes$^{1,2}$ \thanks{supported by S\"achsisches
Staatsministerium   f\"ur Wissenschaft und Kunst,\newline\hspace*{.5cm} e-mail
   Rainer.Matthes@itp.uni-leipzig.de or rmatthes@mis.mpg.de}\\[.5cm]
\normalsize $^1$Institut f\"ur Theoretische Physik
\normalsize der Universit\"at Leipzig\\
\normalsize Augustusplatz 10/11,
\normalsize D-04109 Leipzig,
\normalsize Germany
 \\[.5cm]
\normalsize $^2$Max-Planck-Institut f\"ur Mathematik
\normalsize in den Naturwissenschaften\\
\normalsize Inselstra{\ss}e 22-26,
\normalsize D-04103 Leipzig,
\normalsize Germany}
\date{}
\maketitle
\newtheorem{pr}{Proposition}
\newtheorem{de}{Definition}
\newtheorem{lem}{Lemma}
\newtheorem{them}{Theorem}
\begin{abstract}
Budzy\'nski and Kondracki (\cite{Kon}) have introduced a notion of locally  
trivial quantum principal fibre bundle making use of an algebraic notion of
covering, which allows a reconstruction of the bundle from local pieces.  
Following this approach, 
we construct covariant differential algebras and connections on locally 
trivial quantum principal fibre bundles by gluing together such locally 
given geometric objects.
We also consider covariant derivatives, connection
forms, curvatures and curvature forms and explore the relations between these
notions. As an example, a $U(1)$ quantum principal bundle over a glued quantum
sphere as well as a connection in this bundle is constructed. The connection
may be considered as a $q$-deformed Dirac monopole.
\end{abstract}
{\small 1991 MSC: 81R50, 46L87\\
Keywords: quantum principal bundle, differential structure, covariant
derivative, connection, $q$-monopole}
\section{Introduction}
Since the appearance of quantum groups there has been a hope that it should
be
possible to use them instead of the classical symmetry groups of physical
theories, in particular for quantum field theories. It was expected that
the 
greater variety of group-like structures should lead, perhaps, to greater 
flexibility in the formulation of physical theories, thereby paving the way
to a better understanding of fundamental problems of quantum theory and 
gravitation.\\
In (Lagrangian) quantum field theory, symmetry groups can be considered to 
appear in a very natural geometrical scheme: They are structure groups of 
principal fibre bundles. Moreover, on the classical level, all fields are
geometrical
objects living on the principal bundle or on associated fibre bundles.
Thus, 
it is natural to
ask for a generalization of the notion of principal bundle to a
noncommutative
situation. Thereby, in order to avoid unnecessary restrictions, one should 
replace not only the structure group by a quantum group, but also the base 
manifold (space-time) by a  noncommutative space, which may even be
necessary
for physical reasons (see \cite{dop}, \cite{dfr}, \cite{fgr} and
\cite{kem4}).

In recent years, there have been several attempts to define such quantum  
principal bundles and the usual geometric objects that are needed to
formulate 
gauge field theories on them, see \cite{brma}, \cite{Kon}, \cite{dur}, 
\cite{haj}, \cite{mue}, \cite{pf1}, \cite{scha} and \cite{schn}. Roughly  
following the same idea (``reversing the arrows''), the approaches  differ
in 
the details of the definitions. Closest to the classical idea that
a locally trivial bundle should be imagined as being glued together from  
trivial pieces is the definition given in \cite{Kon}. There, one starts with
the notion of a covering of a quantum space. Being in the context of 
C*-algebras, a covering is defined to be a (finite) family of closed ideals
with zero intersection, which is easily seen to correspond to finite
coverings by closed sets in the commutative case. 
C*-algebras which have such a covering can be reconstructed from their 
``restriction'' to the elements of the covering by a gluing procedure.
Such a reconstruction is not always possible for general (not C*-)algebras,
as was noticed in \cite{cama}. The aim of \cite{cama} was to introduce
differential calculi over algebras with covering. Leaving the C*-category,
one is confronted with the above difficulty, called ``noncompleteness of a
covering''. Nevertheless, making use of ``covering completions'', if
necessary, a general scheme for differential calculi on quantum spaces with
covering was developed, and the example of the gluing of two quantum discs,
being homeomorphic to the quantum sphere $S^2_{\mu c}$, $c>0$, including
the gluing of suitable differential calculi on the discs, was described in
detail.

In \cite{Kon}, a locally trivial quantum principal fibre bundle having as
base $B$ such a quantum space
with covering, and as fibre a compact quantum group $H$, is defined as a
right
$H$-comodule algebra with a covering adapted to the covering of the base.
``Adapted'' means that the ideals defining the covering appear as
kernels 
of ``locally trivializing'' homomorphisms such that
the intersections of these kernels with the embedded base are just the 
embeddings of the ideals defining the covering of $B$. Given such a locally
trivial principal fibre bundle, one can define analogues of the classical 
transition functions which have the usual cocycle properties. Reversely,
given
such a cocycle one can reconstruct the bundle. The transition functions are
algebra homomorphisms $H\rightarrow B_{ij}$, where $B_{ij}$ is the algebra
corresponding to the ``overlap'' of two elements of the covering of $B$. It
turns out that they must have values in the center of $B_{ij}$, which is
related to the fact that principal bundles with structure group $H$ are 
determined by bundles which have as structure group the classical subgroup
of $H$, see \cite{Kon}.

The aim of the present paper is to introduce notions of differential
geometry
on locally trivial bundles in the sense
of \cite{Kon} in such a way that all objects can be glued together from
local pieces. 

Let us describe the contents of the paper:
In Section \ref{local}, locally trivial principal bundles are defined 
slightly different from \cite{Kon}. Not assuming C*-algebras, we add to
the definition of \cite{Kon} the assumption that the
``base'' algebra is embedded
as the algebra of right invariants into the ``total space'' algebra. This 
assumption has to be made in order to come back to the usual notion in the
classical case, as is shown by an example. We prove a technical proposition
about the kernels of the local trivializing homomorphisms which in turn
makes
it possible to prove a reconstruction theorem for locally trivial principal
bundles in terms of transition functions in the context of general
algebras. 

The aim of Section \ref{adco} is to introduce differential calculi on 
locally trivial
quantum principal bundles. They are defined in such a way that they are 
uniquely determined by giving differential calculi on the ``local pieces''
of the base and a right covariant differential calculus on the Hopf
algebra (assuming that the calculi on the trivializations are graded tensor
products).
Uniqueness follows from the assumption that the local trivializing 
homomorphisms should be differentiable and that the kernels of their 
differential extensions should form a covering of the differential calculus
on the total space, i. e. the differential calculus is ``adapted'' in the
sense of \cite{cama}. This covering 
need
not be complete. Thus, in order to have reconstructability, one has to use
the covering completion, which in general is only a differential algebra.

Section \ref{con} is the central part of the paper. Whereas in the
classical situation 
there is a canonically given vertical part in the tangent space of a
bundle,
in the dual algebraic situation there is a canonically given horizontal
subbimodule in the bimodule of forms of first degree on the bundle space.
We start with the definition of left (right) covariant derivatives, which
involves a Leibniz rule, a covariance condition, invariance of the submodule
of horizontal forms, and a locality condition. Covariant derivatives can be
characterized by families of linear maps $A_i:H\lra \hG(B_i)$ satisfying
$A_i(1)=0$ and a compatibility condition being analogous to the classical
relation between local connection forms. At this point a bigger differential
algebra on the basis $B$ appears, which is maximal among all the (LC)
differential algebras being embeddable into the differential structure
of the total space.  
Next we define left (right) connections as a choice of a projection of the left
(right) ${\cal P}$-module of one-forms onto the submodule of horizontal forms
being covariant under the right coaction and satisfying a locality condition.
This is equivalent to the choice of a vertical complement to the submodule
of horizontal forms.
Left and right connections are equivalent. With this definition it is possible
to reconstruct a connection from connections on the local pieces of the
bundle. The corresponding linear 
maps $A_i:H\rightarrow \Gamma(B_i)$  satisfy the conditions for the $A_i$ of 
covariant derivatives, and in addition $R\subset ker A_i$ ($S^{-1}(R)\subset
ker A_i$), where $R$ is the right ideal in $H$ defining the right covariant
differential calculus there. Thus, connections are special cases of covariant
derivatives.
There is a  corresponding notion of connection form as well as a
corresponding notion of an exterior covariant derivative. The curvature
can be defined as the
square of the exterior covariant derivative, and is nicely related to a 
curvature form being defined by analogues of the structure equation.
The local
components of the curvature are related to the local connection forms in a
nice way, and they are related among themselves by a homogeneous formula
analogous to the classical one.

Finally, in Section \ref{ex}, we give an example of a locally trivial
principal bundle
with a connection. The basis of the bundle, constructed in \cite{cama}, is a
C*-algebra glued together from
two copies of a quantum disc. The structure group is the classical
group $U(1)$, and the bundle is defined by giving one transition
function, which is sufficient because the covering of the basis has
only two elements. Since all other coverings appearing in the example
then have also two elements, there 
are no problems with noncomplete coverings. 
The differential calculus on the total space is determined by 
differential ideals in the universal differential calculi over the
two copies of the quantum disc and the structure group. 
For the group, the ideal is chosen in nonclassical way. Then, a connection
is
defined by 
giving explicitely the two local connection forms on the generator of 
$P(U(1))$ and extending them using the properties a local 
connection form should have. The curvature of this 
connection is nonzero. 

In the appendix, the relevant facts about coverings and gluings of algebras
and differential algebras are collected, for the convenience of the reader.
Details can be found in \cite{cama}. Moreover, we recall there some
well-known
facts about covariant differential calculi on quantum groups.

In the following, algebras are always assumed to be over $\C$,
associative
and unital. Ideals are assumed to be two-sided, up to some occasions, where 
their properties are
explicitely specified.

\section{\label{local} Locally trivial quantum principal fibre bundles}
Following the ideas of \cite{Kon} we introduce in this section the
definition
of a locally trivial quantum principal fibre bundle and prove propositions 
about the existence of trivial subbundles and about the reconstruction
of the bundle. Essentially, this is contained in \cite{Kon},
up to some modifications: We do not assume C*-algebras, and we add to the
axioms the condition that the embedded base algebra coincides with the
subalgebra of coinvariants. As structure group we take a general Hopf
algebra.
\\In the sequel we use the results of \cite{cama}, see in the appendix. We 
recall here that for an algebra $B$ with a covering $(J_i)_{i\in I}$, there are 
canonical mappings $\pi_i : B \longrightarrow B_i:=B/J_i$, $\pi^i_j: B_i 
\longrightarrow B_{ij}:=B/(J_i + J_j)$, $\pi_{ij}: B \longrightarrow B_{ij}$.
\begin{de} \label{F} A locally trivial quantum principal fibre bundle
(QPFB)
is a tupel \begin{equation} ({\cal P},\Delta_{\cal P},H,B,\iota,
(\chi_i,J_i)_{i \in I}) \end{equation}
where $B$ is an algebra, $H$ is a Hopf algebra, $\cal P$ is a right $H$
comodule algebra with coaction $\Delta_{\cal P}$, $(J_i)_{i \in I}$ is a  
complete covering of $B$, and $\chi_i$ and $\iota$ 
are homomorphisms with the following properties:
\begin{eqnarray*} \chi_i : {\cal P} &\longrightarrow & B_i \otimes H \;\;\;
\mbox{surjective}, \\
\iota : B & \longrightarrow & {\cal P} \;\;\; \mbox{injective}, \\
(id \otimes \Delta) \circ \chi_i &=& (\chi_i \otimes id) \circ \Delta_{\cal
P}, \\
\chi_i \circ \iota(a) &=& \pi_i(a) \otimes 1\;\;a \in B, \\
(ker \chi_i )_{i \in I} & & \mbox{complete covering of ${\cal P}$}, \\
\iota(B)&=&\{f \in {\cal P}|\Delta_{\cal P}(f)=f \otimes I\}.
\end{eqnarray*} 
\end{de}
Such a tupel we often denote simply by $\cal P$. Occasionally, ${\cal P}$,
$B$ and $H$ are called total space, base space and structure group of the
bundle. 
\\The last assumption in Definition \ref{F} does not appear in the
definition of QPFB given in \cite{Kon}. It is however used by other 
authors ( \cite{br}, \cite{dur}, \cite{pf1}). 
Already in the classical case this condition is needed to 
guarantee the transitive action of the structure group on the fibres, as
shows the following example.
\\\\Example: Let $M$ be a compact topological space covered by two closed 
subsets $U_1$ and $U_2$ being the closure of two open subsets covering	
$M$. Define $M_0=U_1 \dot{\bigcup} U_2$ (disjoint union). $M$ is obtained 
from $M_0$ identifying all corresponding points of $U_1$ and $U_2$.
There is a natural projection $M_0 \longrightarrow M$.
Let us consider the algebras of continuous functions 
$C(M)$ and $C(M_0)$ over $M$ and $M_0$ respectively.
There exists an injective homomorphism $\kappa:C(M) \longrightarrow C(M_0)$
being the pull back of the natural projection $M_0 \longrightarrow M$.
Suppose 
we have constructed a principal fibre bundle $P$ over $M_0$ with structure 
group $G$, which is trivial on each of the disjoint components. Then we
have an injective homomorphism $\iota_o: C(M_0) \longrightarrow C(P)$ 
and two trivialisations $\chi_{1,2}:C(P) \longrightarrow C(U_{1,2}) \otimes
C(G)$
with the properties assumed in Definition \ref{F}. The	injective 
homomorphismus $\iota :C(M) \longrightarrow 
C(P)$, $\iota:=\iota_0 \circ \kappa$, fullfills all the assumptions in
Definition \ref{F} up to the last one, and one obtains a fibration $P$ over
the base manifold $M$ which is not a principal fibre bundle. 
\begin{pr} Let ${\cal P}_c$ be the covering completion of $\cal P$ with
respect 
to the complete covering $(ker \chi_i)_{i \in I}$. Let $K: {\cal P} 
\longrightarrow {\cal P}_c$ be the corresponding isomorphism. The tupel
\[ ( {\cal P}_c, \Delta_{{\cal P}_c}, H,B, 
\iota_c, (\chi_{i_c},J_i)_{i \in I}),  \] where 
\begin{eqnarray*}  \Delta_{{\cal P}_c} &=& (K \otimes id) \circ
\Delta_{\cal P} \circ K^{-1}, \\ \chi_{i_c} &=& \chi_i \circ K^{-1}, \\ 
\iota_c &=& K \circ \iota, \end{eqnarray*}
is a locally trivial QPFB. 
\end{pr} 
The proof ist trivial (transport of the structure using $K$).
\begin{de} \label{trivP} A locally trivial QPFB $\cal P$ is called trivial 
if there exists an
isomorphism $\chi : {\cal P} \longrightarrow B \otimes H$ such that
\begin{eqnarray*} \chi \circ \iota &=& id \otimes 1, \\
(\chi \otimes id) \circ \Delta_{\cal P} &=& (id \otimes \Delta) \circ \chi.
\end{eqnarray*} \end{de}
Remark: A locally trivial QPFB with $card I=1$, i.e. with trivial covering
of $B$, is trivial.
\\Triviality of the covering means that it consists of only one ideal 
$J=0$. Moreover, there is only one trivializing epimorphism $\chi: {\cal P} 
\longrightarrow B \otimes H$ which necessarily fulfills $ker \chi=0$.
\\\\There are several trivial QPFB related to a locally trivial QPFB.
Define 
${\cal P}_i:= {\cal P}/ker \chi_i$. Then $\tilde{\chi_i}: {\cal P}_i 
\longrightarrow B_i \otimes H$ defined by 
\begin{equation} \label{monst} \tilde{\chi}_i(f + ker \chi_i):= 
\chi_i(f) \end{equation} is a well defined isomorphism. $\iota_i: B_i \longrightarrow
{\cal 
P}_i$ defined by \[\iota_i(b):= \tilde{\chi}^{-1}_i(b \otimes 1) \] is 
injective and 
fulfills $\tilde{\chi}_i \circ \iota_i = id \otimes 1$. Moreover
$\Delta_{{\cal 
P}_i}: {\cal P}_i \longrightarrow {\cal P}_i \otimes H$ is well defined by 
\[ \Delta_{{\cal P}_i}(f + ker \chi_i):= \Delta_{\cal P}(f) + ker \chi_i
\otimes H, \] because from $ (id \otimes \Delta) \circ \chi_i = (\chi_i 
\otimes id) \circ 
\Delta_{\cal P}$ follows $\Delta_{\cal P}(ker \chi_i) \subset ker \chi_i
\otimes H$. Obviously, $\Delta_{{\cal P}_i}$ is a right coaction. Moreover, 
$(\tilde{\chi}_i \otimes id) \circ \Delta_{{\cal P}_i} = ( id \otimes
\Delta) \circ \tilde{\chi}_i$, and $\iota_i(B_i)=\{f \in {\cal 
P}_i|~\Delta_{{\cal P}_i}(f) = f \otimes 1\}$. Thus $({\cal P}_i, 
\Delta_{{\cal P}_i}, H,B_i,\iota_i,(\tilde{\chi}_i,0))$ is a trivial QPFB.
\\Let ${\cal P}_{ij}:= {\cal P}/(ker \chi_i + ker \chi_j)$. Then there is
an isomorphism $\tilde{\chi}^i_{ij}: {\cal P}_{ij} \longrightarrow (B_i
\otimes H)/\chi_i(ker \chi_j)$ given by \begin{equation} \label{shiijt}
\tilde{\chi}^i_{ij} (f + ker \chi_i + ker \chi_j):=\chi_i(f) + 
\chi_i(ker \chi_j).
\end{equation}
It is natural to expect that $P_{ij}$ should be a trivial bundle isomorphic
to $B_{ij} \otimes H$. In fact, we will show that there is a natural 
isomorphism 
$(B_i \otimes H)/\chi_i(ker \chi_j) \simeq B_{ij} \otimes H$, leading to
trivialization maps  $\chi^i_{ij}:{\cal P}_{ij} \longrightarrow B_{ij}
\otimes H$.
Let us introduce the natural projections $\pi_{i_{\cal P}}: {\cal P} 
\longrightarrow {\cal P}_i$, $\pi_{{ij}_{\cal P}}: {\cal P} \longrightarrow
{\cal P}_{ij}$ and $\pi^i_{j_{\cal P}}: {\cal P}_i \longrightarrow {\cal  
P}_{ij}$. Obviously, $\tilde{\chi}_i \circ \pi_{i_{\cal P}}= \chi_i$, $ 
\pi_{i_{\cal P}}= \tilde{\chi}^{-1}_i \circ \chi_i$ and
$\pi_{{ij}_{\cal P}}= \pi^i_{j_{\cal P}} \circ \pi_{i_{\cal P}}$.
We will need the following lemma.
\begin{lem} \label{lemisubb} Let $B$ be an algebra and $H$ be a Hopf
algebra. Let $J \subset B \otimes H$ be an ideal with the property 
\[ (id \otimes \Delta)J \subset J \otimes H. \] Then there exists an 
ideal $I \subset B$ such that $ J =I \otimes H$. This ideal is uniquely 
determined and equals $(id \otimes \varepsilon)(J)$.\end{lem}
Proof: Let $m_H : H \otimes H \longrightarrow H$ be the
algebra product in $H$. It is not difficult to verify that $I:=(id \otimes 
\varepsilon)(J)$ is an ideal in $B$. We will show $J=I \otimes H$ 
First we prove 
$J \subset I \otimes H$. Because of $(id \otimes \varepsilon \otimes id)  
\circ (id \otimes \Delta)= id$ and $(id \otimes \Delta)J \subset J \otimes
H$ we have $(id \otimes \varepsilon \otimes id) \circ (id \otimes \Delta)J=J
\subset I \otimes H$. $I \otimes H \subset J$ is a consequence of $I \otimes 1 
\subset J$, which is proved as follows:
A general element of $I$ has the form $\sum_k a_k \varepsilon(h_k)$ where 
$\sum_k a_k \otimes h_k \in J$. Because of \[ \sum_k a_k \varepsilon(h_k) 
\otimes 1 = \sum_k \sum (a_k \otimes h_{k_1})(1 \otimes S(h_{k_2})) \]
and 
\[ (id \otimes \Delta)(\sum_k a_k \otimes h_k)=\sum_k a_k \otimes 
h_{k_1} \otimes h_{k_2} \in J \otimes H, \] $\sum_k a_k \varepsilon (h_k)
\otimes 1$ is an element of $J$.\hfill$\Boxneu$
\begin{pr} ${\cal P}_{ij}$ is a trivial QPFB, i.e. there exist	
\begin{eqnarray*} \chi^i_{ij} : {\cal P}_{ij} &\longrightarrow& B_{ij}
\otimes H, 
\\ \Delta_{{\cal P}_{ij}} : {\cal P}_{ij} &\longrightarrow& {\cal P}_{ij}
\otimes 
H, \\  \iota_{ij} : B_{ij} &\longrightarrow & {\cal P}_{ij}, \end{eqnarray*}
such 
that the conditions of Definitions \ref{F} and \ref{trivP} are satisfied. 
\end{pr}
Remark: ${\cal P}_{ij}$ is a trivial QPFB in two ways by choosing $\chi^i_{ij}$ 
or $\chi^j_{ij}$. The composition of these maps just gives the transition 
functions.
\\\\Proof: Applying $\chi_i \otimes id$ to $\Delta_{\cal P}(ker
\chi_i) \subset ker 
\chi_i \otimes H$ and using $(id \otimes \Delta) \circ \chi_i = (\chi_i
\otimes id) 
\circ \Delta_{\cal P}$ it follows that $(id \otimes \Delta) \circ
\chi_i(ker 
\chi_j) \subset \chi_i( ker \chi_j) \otimes H$. By Lemma \ref{lemisubb},
there exist ideals $\tilde{K}^i_j \subset B_i$ such that 
$\chi_i(ker \chi_j)=\tilde{K}^i_j \otimes H$. $K^i_j :=\pi^i_j(\tilde{K}^i_j)$ 
is an ideal in $B_{ij}$.			     
\\We show now $\pi_i(J_j) \subset \tilde{K}^i_j$: According to Lemma 
\ref{lemisubb}, we have $\tilde{K}^i_j= (id \otimes \varepsilon)(\chi_i(ker 
\chi_j))$. We have to show that for a $b \in J_j$ there exists 
$\tilde{b} \in ker 
\chi_j$ with $(id \otimes \varepsilon) \circ \chi_i(\tilde{b})= \pi_i(b)$. 
It is obvious that we can take $\tilde{b}=\iota(b)$.
\\Using this inclusion, one finds that there is
a canonical isomorphism $(B_i \otimes H)/\chi_i(ker \chi_j)
\simeq (B_{ij}/K^i_j) \otimes H$ given by $b \otimes h+ \chi_i(ker \chi_j)
\longrightarrow (\pi_i(b) +K^i_j) \otimes h$. Composing with 
$\tilde{\chi}^i_{ij}$ (see (\ref{shiijt})), 
there results an isomorphism $\chi^i_{ij}: {\cal P}_{ij} \longrightarrow  
B_{ij}/K^i_j \otimes H$ given by \[ \chi^i_{ij}(f + ker \chi_i + ker 
\chi_j):=(\pi^i_j \otimes id) \circ \chi_i(f)+ K^i_j \otimes H.\]
Our goal is now to show $K^i_j= \pi^i_j(\tilde{K}^i_j)=0$, so that $\chi^i_{ij}$ 
will become the isomorphism ${\cal P}_{ij} \longrightarrow B_{ij} \otimes H$ 
wanted in the proposition. 
As a first step we will prove $K^i_j=K^j_i$.
To this end, we note that \[ \tilde{\phi}_{ji}:= \chi^j_{ij} 
\circ \pi^i_{j_{\cal P}} \circ \tilde{\chi_i}^{-1} \] is a homomorphism  
 $\tilde{\phi}_{ji}: B_i \otimes H 
\longrightarrow B_{ij}/K^j_i \otimes H$ 
with $ker \tilde{\phi}_{ji} =\tilde{K}^i_j \otimes H$. In terms of this 
homomorphism we define a homomorphism $\psi_{ji}: B_{ij} \longrightarrow  
B_{ij}/K^j_i$ by \[\psi_{ji}(a + J_i+ J_j):=(id \otimes \varepsilon) \circ 
\tilde{\phi}_{ji}((a+J_i) \otimes 1).\]
It is easy to show that $ker\psi_{ji}=K^i_j$. On the other 
hand one shows that  $\psi_{ji} : B_{ij} 
\longrightarrow B_{ij}/K^j_i$ is the natural projection, and therefore 
$K^i_j =K^j_i$:
\\We calculate 
\begin{eqnarray*} \psi_{ji}(a+J_i + J_j)&=&(id \otimes \varepsilon) \circ 
\tilde{\phi}_{ji} ((a+J_i) \otimes 1)\\ 
&=& (id \otimes \varepsilon) \circ \chi^j_{ij} \circ 
\pi^i_{j_{\cal P}} \circ \tilde{\chi}^{-1}_i((a+J_i) \otimes 1)\\
&=& (id \otimes \varepsilon) \circ \chi^j_{ij} \circ \pi_{{ij}_{\cal P}}
 (\iota(a)) \\
&=& (id \otimes \varepsilon) \circ \chi^j_{ij}(\iota(a) + ker \chi_i +  
ker \chi_j)\\ &=& (id \otimes \varepsilon) ( (\pi^j_i \otimes id)\circ  
\chi_j(\iota(a)+ K^j_i \otimes H \\ &=& (id \otimes \varepsilon)
((\pi^j_i \otimes id)(\pi_j(a) \otimes 1)+ K^j_i \otimes H) \\	
&=& \pi_{ij}(a) + K^j_i. \end{eqnarray*}
\\For showing $K^i_j=0$ we use the completeness of the covering 
$(ker \chi_i)_{i \in I}$. The covering 
completion of ${\cal P}$ is by definition \[{\cal P}_c=\{ ( f_i)_{i \in
I} \in \bigoplus_{i \in I} {\cal P}/ker\chi_i| 
\pi^i_{j_{\cal P}}(f_i)=\pi^j_{i_{\cal P}}(f_j) \}. \]
We introduce a locally trivial QPFB $\breve{\cal P} \simeq {\cal P}_c$ such
that a comparison of $\breve{\cal P}^{coH}=\{f \in \breve{\cal P}| 
\Delta_{\breve{\cal P}}(f)=f \otimes 1\}$ with $B \simeq B_c$
allows to read off $ker \psi_{ij}=K^i_j=0$.
Let $\phi_{ij} : B_{ij}/K^i_j \otimes H \longrightarrow B_{ij}/K^i_j
\otimes H$ be the isomorphisms defined by \[ \phi_{ij}:=\chi^i_{ij} \circ 
{\chi^j_{ij}}^{-1}. \] 
Using the identities \[ \chi^i_{ij} \circ \pi^i_{j_{\cal P}} \circ 
\tilde{\chi}^{-1}_i=(\psi_{ij} \otimes id) \circ (\pi^i_j \otimes id) \]
it is easy to verify that the algebra ${\cal P}_c$ is isomorphic 
to the algebra \[ \breve{\cal P}=\{ (g_i)_{i \in I} \in \bigoplus_{i \in I} 
(B_i \otimes H)|
(\psi_{ji} \otimes id)\circ (\pi^i_j \otimes id)(g_i)=\phi_{ij} \circ
(\psi_{ji} \otimes id) \circ (\pi^j_i \otimes id)(g_j) \} \]
(cf. Lemma 1 in \cite{cama}), and the corresponding
isomorphism 
$\chi :{\cal P}_c \longrightarrow \breve{\cal P}$ is defined by 
$\chi (( f_i)_{i \in I}):= ( \tilde{\chi}_i(f_i))_{i \in I}$. Transporting the 
homomorphisms $\Delta_{{\cal P}_c}$, $\chi_{i_c}$ and $\iota_c$ to 
$\Delta_{\breve{\cal P}}:=(\chi \otimes id) \circ \Delta_{{\cal P}_c} 
\circ \chi^{-1}$, $\breve{\chi}_i:=\chi_{i_c} \circ 
\chi^{-1}$ and $\breve{\iota}:=\chi \circ \iota_c$ respectively, one obtains 
a locally trivial 
QPFB again. Explicitly, \begin{eqnarray*}
\Delta_{\breve{\cal P}}((g_i)_{i \in I})&=&((id \otimes \Delta)
(g_i))_{i \in I},\\
\breve{\chi}_i((g_k)_{k \in I})&=& g_i, \\ \breve{\iota}(a) &=& (\pi_i(a) 
\otimes 1)_{i \in I}. \end{eqnarray*}
\\We note that the isomorphisms $\phi_{ij}$ fulfill 
\begin{eqnarray} \label{deltaij} (id \otimes \Delta) \circ \phi_{ij} &=&
(\phi_{ij} \otimes id) \circ ( id \otimes \Delta),\\ \label{iotaij}
\phi_{ij}(a \otimes 1) &=& a \otimes 1,\;\; a \in B_{ij}/K^i_j. 
\end{eqnarray}
Using (\ref{deltaij}) and (\ref{iotaij}) it follows that the subalgebra 
$\breve{\cal P}^{coH}=\{f \in \breve{\cal P}|\Delta_{\breve{\cal P}}(f)=
f \otimes 1 \}$ is isomorphic to \[ \breve{\cal P}^{coH} =\{\{ (a_i \otimes
1)\}_{i \in I} \in \bigoplus_{i \in I} B_i \otimes 1| \psi_{ji}\circ \pi^i_j(a_i) 
\otimes 1 =\psi_{ji} \circ \pi^j_i(a_j) \otimes 1\}.\]
This algebra is by Definition \ref{F} isomorphic to 
\[ B \simeq B_c=\{ (a_i)_{i \in I} \in \bigoplus_{i \in I} B_i| 
\pi^i_j(a_i)=\pi^j_i(a_j)\} \] (see \cite{cama} and the appendix). It follows
that the $\psi_{ij}$ have to be isomorphisms, i.e. $ker \psi_{ij}=K^i_j=0$,
which means in fact $\psi_{ij}=id$. Thus, the $\chi^i_{ij}:{\cal P}_{ij} 
\longrightarrow B_{ij} \otimes H$ are isomorphisms. 
\\Further define $\Delta_{{\cal P}_{ij}}: {\cal P}_{ij} \longrightarrow
{\cal P}_{ij} \otimes H$ by \[\Delta_{{\cal P}_{ij}}(f + ker \chi_i + ker 
\chi_j):=\Delta_{\cal P}(f) + (ker \chi_i + ker \chi_j) \otimes H \] and 
$\iota_{ij}: B_{ij} \longrightarrow {\cal P}_{ij}$ by \[ 
\iota_{ij}(\pi_{ij}(a)):= \iota(a)+ ker \chi_i + ker \chi_j.\]
It is easy to verify that all the conditions of Definition \ref{trivP} are 
satisfied.\hfill$\Boxneu$
\\\\Notice that due to $K^i_j=0$, we have $\tilde{K}^i_j= \pi_i(J_j)$. This 
means \begin{equation} \chi_i(ker \chi_j)=\pi_i(J_j) \otimes H. \end{equation} 
The isomorphisms $\chi^i_{ij}$ satisfy 
\begin{equation} 
\label{blubblub} \chi^i_{ij} \circ \pi_{{ij}_{\cal P}}= (\pi^i_j 
\otimes id) \circ \chi_i, \end{equation} and the $\phi_{ij}$ defined above 
are isomorphisms $B_{ij} \otimes H \longrightarrow B_{ij} \otimes H$ 
fulfilling (\ref{deltaij}), (\ref{iotaij}) and $\phi_{ij} \circ \phi_{ji}=id$.
\begin{pr}\label{transf} (cf. \cite{Kon}) Locally trivial QPFB's over a
basis $B$ with complete covering 
$(J_i)_{i \in I}$ and with structure group H are in one-to-one 
correspondence with families of 
homomorphisms  \[\tau_{ij} :H \longrightarrow B_{ij}, \] called transition 
functions, satisfying the conditions
\begin{eqnarray*} \tau_{ii}(h)&=& 1\varepsilon(h) \;\; \forall h \in
H,\\ \tau_{ji}(S(h)) &=& \tau_{ij} (h) \;\; \forall h \in H,\\
\tau_{ij}(h)a &=& a\tau_{ij}(h) \;\; \forall a \in B_{ij} \;h \in H, \\
\pi^{ij}_k \circ \tau_{ij}(h)&=& m_{B_{ijk}} \circ (\pi^{ik}_j \circ
\tau_{ik} \otimes \pi^{jk}_i \circ \tau_{kj}) \circ \Delta(h) \;\; \forall h 
\in H. \end{eqnarray*} \end{pr}
Proof: Let a bundle ${\cal P}$ be given and let the $\phi_{ij}: B_{ij}
\otimes H \longrightarrow B_{ij} \otimes H$ be defined as above. 
Define homomorphisms $\tau_{ij} : H \longrightarrow B_{ij}$ by 
\begin{equation} \label{monst1} \tau_{ji}(h):=(id \otimes \varepsilon) 
\circ \phi_{ij}(1 \otimes h).\end{equation} (There is another possible choice, 
$\tau_{ij}(h):=(id \otimes \varepsilon) \phi_{ij}(1 \otimes h)$, which 
correspond to another form of the cocycle condition.) 
One shows that this 
is equivalent to \begin{equation} \label{phi} 
\phi_{ij}(a \otimes h)=\sum a \tau_{ji}(h_1) \otimes h_2: \end{equation}
Using (\ref{deltaij}), (\ref{iotaij}) and $(\varepsilon \otimes id) \circ
\Delta =id$ it follows from (\ref{monst1}) that
\begin{eqnarray*} \sum a \tau_{ji}(h_1) \otimes h_2&=& \sum (a 
\otimes 1) ((id \otimes \varepsilon) \circ \phi_{ij}(1 \otimes
h_1) \otimes h_2) ,\\
&=&(a \otimes 1)(id \otimes \varepsilon \otimes id) \circ (\phi_{ij}
\otimes id)\circ (id \otimes \Delta)(1 \otimes h),\\
&=&(a \otimes 1)(id \otimes \varepsilon \otimes id) \circ (id \otimes
\Delta) \circ \phi_{ij} (1 \otimes h), \\
&=& (a \otimes 1)\phi_{ij}(1 \otimes h), \\ 
&=& \phi_{ij}(a \otimes h). \end{eqnarray*}
Conversely, if (\ref{phi}) is satisfyed, the choice $a=1$ gives 
(\ref{monst1}).
$\tau_{ii}(h)=\varepsilon(h)1$ follows from $\phi_{ii}= id$.
Every homomorphism $\tau_{ij}: H \longrightarrow B$ is convolution 
invertible with convolution inverse $\tau^{-1}_{ij}= \tau_{ij} \circ S$.
On the other hand from $\phi_{ij} \circ \phi_{ji}=id$ easily follows 
$\tau^{-1}_{ij}=\tau_{ji}$: \[\phi_{ij} \circ \phi_{ji}(1 \otimes 
h)=\phi_{ij}(\sum \tau_{ij}(h_1) \otimes h_2)=\sum
\tau_{ij}(h_1)\tau_{ji}(h_2) 
\otimes h_3=1 \otimes h. \] Therefore $\sum \tau_{ij}(h_1) 
\tau_{ji}(h_2)=\varepsilon(h)1$, i.e. $\tau_{ji}=\tau_{ij} \circ
S$. $\tau_{ij}$ has values in the center of $B_{ij}$:
\begin{eqnarray*} a\tau_{ij}(h)-\tau_{ij}(h)a&=&a(id \otimes 
\varepsilon)\phi_{ji}(1 \otimes h)-(id \otimes 
\varepsilon)\phi_{ji}(1 \otimes h)a \\
&=& (id \otimes \varepsilon)((a \otimes 1)\phi_{ji}(1 \otimes 
h)-\phi_{ji}(1 \otimes h)(a \otimes 1)) \\
&=& (id \otimes \varepsilon)(\phi_{ji}((a \otimes h)-(a \otimes h))=0. 
\end{eqnarray*}
To prove the last relation of the proposition, define isomorphisms 
$\phi^k_{ij} :B_{ijk} \otimes H 
\longrightarrow B_{ijk} \otimes H$ by \[ \phi^k_{ij}((a \otimes 
h)+\pi_{ij}(J_k) \otimes H):=\phi_{ij}(a \otimes h) +\pi_{ij}(J_k) \otimes 
H \] (using $B_{ijk} \simeq B_{ij}/\pi_{ij}(J_k)$).
$\phi^k_{ij}$ are well defined because of $\phi_{ij}(a \otimes 1)=a
\otimes 1$. Now, a lengthy but simple 
computation leads to \[ \phi^k_{ij}=\phi^j_{ik} \circ \phi^i_{kj}. \]
The idea of this computation is to consider the isomorphism $\chi^{ijk}_i: 
{\cal P}/(ker \chi_i +ker \chi_j +ker \chi_k) \longrightarrow B_{ijk} \otimes H$ 
induced by $\chi_i$ and to prove $\phi^k_{ij}= \chi^{ijk}_i \circ 
{\chi^{ijk}_j}^{-1}$.
\\Combining the definition of $\phi^k_{ij}$ with (\ref{phi}), one obtains \[ 
\phi^k_{ij}(a \otimes h)=\sum a\pi^{ij}_k \circ \tau_{ji}(h_1) \otimes
h_2.\]
Therefore, \[\pi^{ij}_k \circ \tau_{ji}(h)= (id \otimes \varepsilon) \circ 
\phi^k_{ij}(1 \otimes h).\]
Inserting here $\phi^k_{ij}(1 \otimes h)=\phi^j_{ik} \circ
\phi^i_{kj} 
(1 \otimes h)$ one obtains \[ \pi^{ij}_k \circ \tau_{ji}(h)= \sum
\pi^{jk}_i \circ \tau_{jk}(h_1) \pi^{ik}_j \circ \tau_{ki}(h_2). \]
This ends the proof of one direction of the proposition.
\\We will not give the details of reconstruction of the bundle from the
transition 
functions. We only remark, that, for a given family of transition functions 
$\tau_{ij}$, we define the isomorphisms $\phi_{ij}$ by formula (\ref{phi}),
which 
gives rise to the gluing 
\begin{equation}
\label{P}
\breve{\cal P}=\{(f_i)_{i\in I} \in \bigoplus_{i \in I} (B_i \otimes H)|
(\pi^i_j \otimes id)(f_i)=\phi_{ij} \circ (\pi^j_i \otimes id)(f_j) \}.
\end{equation}
One verifies that the formulas	\begin{eqnarray}
\Delta_{\breve{\cal P}} ((f_i)_{i \in I})&=&(id \otimes \Delta(f_i))_{i \in I},~
\forall (f_i)_{i \in I} \in \tilde{\cal P},\\
\breve{\chi}_k((f_i)_{i \in I})&=&f_k,~~ \forall (f_i)_{i \in I} \in 
\tilde{\cal P},\\
\breve{\iota}(a) &=& (\pi_i(a) \otimes 1)_{i \in I},~~\forall a \in B \end{eqnarray}
define a locally trivial QPFB $(\breve{\cal P}, \Delta_{\breve{\cal P}}, H, B, 
\breve{\iota}, (\breve{\chi}_i, J_i)_{i\in I})$. 
If the $\tau_{ij}$ stem from a given locally trivial QPFB $\cal P$, 
applying the isomorphism $\chi^{-1}$ defined as above (proof of proposition 2) 
leads to ${\cal P}_c \simeq \breve{\cal P}$.\hfill$\Boxneu$
\section{\label{adco}Adapted covariant differential structures on locally
trivial QPFB}
In the sequel we will use the skew tensor product of differential calculi.
Let $\Gamma(A)$ and $\Gamma(B)$ be two differential calculi. We define the
differential calculus $\Gamma(A) \hat{\otimes} \Gamma(B)$ as the vector space 
$\Gamma(A) \otimes \Gamma(B)$ equipped with the product 
\begin{equation} (\gamma \hat{\otimes} \rho)(\omega \hat{\otimes} \tau)=
(-1)^{mn}(\gamma \omega \hat{\otimes} \rho \tau), \;\;\omega \in 
\Gamma^n(A),\rho 
\in \Gamma^m(B),\gamma \in \Gamma(A) , \tau \in \Gamma(B) \end{equation}
and the differential \begin{equation} d(\gamma \hat{\otimes} \rho)=(d\gamma
\hat{\otimes} \rho)+(-1)^n(\gamma \hat{\otimes} d\rho),\;\; \gamma \in 
\Gamma^n(A),\rho \in \Gamma(B). \end{equation}
\begin{pr} \label{hat} Let $\Gamma(A)$ and $\Gamma(B)$ be two differential
calculi and let $J(A) \subset \Omega(A)$ and $J(B) \subset \Omega(B)$ be the 
corresponding differential ideals respectively. Let $id \otimes 1 :A  
\longrightarrow
A \otimes B$ and $1 \otimes id: B \longrightarrow A \otimes B$ be the 
embedding homomorphisms. The differential ideal $J(A \otimes B)
\subset \Omega(A \otimes B)$ corresponding to 
$\Gamma(A) \hat{\otimes} \Gamma(B)$ if it is generated by the sets 
\[ (id \otimes 1)_{\Omega}(J(A));
(1 \otimes id)_{\Omega}(J(B)) \] \begin{equation}\label{set} \{(a \otimes
1)d(1 \otimes b)-(d(1 \otimes b))(a \otimes 1)|a \in A, \, b \in B \}. 
\end{equation} \end{pr}
Proof: First we define a homomorphism $\tilde{\psi}: \Omega(A \otimes B)  
\longrightarrow \Gamma(A) \hat{\otimes} \Gamma(B)$ by 
\begin{eqnarray*} \tilde{\psi}(\sum_k (a^0_k \otimes 1)d(a^1_k \otimes
1)) &=& \sum_k a^0_kda^1_k \hat{\otimes} 1, \\
\tilde{\psi}(\sum_k (1 \otimes b^0_k)d(1 \otimes b^1_k)) &=& \sum_k 1
\hat{\otimes} b^0_kdb^1_k. \end{eqnarray*}
It is easy to verify that the differential ideal $\tilde{J}(A \otimes B)$ 
generated
by the sets (\ref{set}) satisfies $\tilde{J}(A \otimes B) \subset ker 
\tilde{\psi}$.
Let $\tilde{\Gamma}(A \otimes B) =\Omega(A \otimes B)/\tilde{J}(A \otimes B)$.	
Note that there are the following relations in $\tilde{\Gamma}(A \otimes
B)$. \begin{eqnarray*} (d(a \otimes 1))(1 \otimes b)&=&(1 \otimes 
b)d(a \otimes 1),\\ d(a \otimes 1)d(1 \otimes b) &=& -d(1 \otimes 
b)d(a \otimes 1) \end{eqnarray*}
Therefore, an element $\gamma \in \tilde{\Gamma}(A \otimes B)$ has the
general form \[ \gamma=\sum_k \sum_l (a^0_k \otimes 1)d(a^1_k \otimes
1)...d(a^n_k 
\otimes 1)(1 \otimes b^0_l)d(1 \otimes b^1_l)...d(1 \otimes b^m_l).
\]
There exist homomorphisms $\Upsilon_A :\Gamma(A)
\longrightarrow 
\tilde{\Gamma}(A \otimes B)$ and $\Upsilon_B :\Gamma(B) \longrightarrow 
\tilde{\Gamma}(A \otimes B)$ defined by \begin{eqnarray*} 
\Upsilon_A(a_0da_1...da_n)&:=&(a_0 \otimes 1)d(a_1 \otimes 1)...d(a_n
\otimes 1),\\
\Upsilon(b_0db_1...db_n)&:=&(1 \otimes b_0)d(1 \otimes b_1)...d(1
\otimes b_n). \end{eqnarray*}
Because of $\tilde{J}(A \otimes B) \subset ker\tilde{\psi}$ the homomorphism $\psi
:\tilde{\Gamma}(A \otimes B) \longrightarrow \Gamma(A) \hat{\otimes}
\Gamma(B)$ defined by \begin{eqnarray*} \psi(d(a \otimes 1)) &=& da
\hat{\otimes} 1, \\ \psi(d(1 \otimes b)) &=& 1 \hat{\otimes} db, 
\end{eqnarray*} exists.
Since $\Gamma(A) \hat{\otimes} \Gamma(B)$ is isomorphic to $\Gamma(A)
\otimes \Gamma(B)$ as vector space, we can define a linear map 
$\psi^{-1}: \Gamma(A) \hat{\otimes} \Gamma(B) \longrightarrow
\tilde{\Gamma}(A \otimes B)$, \[ \psi^{-1}(\alpha \hat{\otimes} \beta):= 
\Upsilon_A(\alpha)\Upsilon_B(\beta). \] 
Since this linear map is a homomorphism and fulfills \begin{eqnarray*} 
\psi^{-1} \circ d &=& d \circ \psi^{-1}, \\ \psi^{-1} \circ \psi &=& id, \\
\psi \circ \psi^{-1} &=& id, \end{eqnarray*}
$\psi$ is an isomorphism, and $\tilde{J}(A \otimes B) = ker \tilde{\psi}$. 
Therefore, $\Gamma(A) \hat{\otimes} \Gamma(B) \simeq \Omega(A \otimes 
B)/\tilde{J}(A \otimes B) \simeq \Omega(A \otimes B)/J(A \otimes B)$ and the 
proposition follows from uniqueness of the differential ideal corresponding to 
a differential calculus.
\hfill$\Boxneu$
\\\\Remark: If we are in the converse situation, i.e. if a 
differential calculus $\Gamma(A \otimes B)$ with corresponding differential
ideal $J(A \otimes B)$ is given, there exist differential ideals $J(A):= J(A
\otimes B) \cap \Omega(A \otimes 1)$ and $J(B):=J(A \otimes B) \cap 
\Omega(1 \otimes B)$. By Proposition \ref{hat}, the differential calculus 
is isomorphic to an algebra of the form $\Gamma(A) \hat{\otimes} \Gamma(B)$ 
if and only if $J(A \otimes B)$ is generated by the sets (\ref{set}).
\\\\In the sequel we always identify ${\cal P}/ker \chi_i$ with $B_i \otimes H$, 
by means of the isomorphisms $\tilde{\chi}_i$ (see (\ref{monst})). 
\\Our goal is now to define differential structures on $\cal P$.
By Proposition (\ref{diff}), a family of differential calculi $\Gamma(B_i)$ and
right covariant differential calculus $\Gamma(H)$ determine unique differential
calculi $\Gamma(B)$ and $\Gamma({\cal P})$ such that
$(\Gamma(B),(\Gamma(B_i))_{i \in I})$ and $(\Gamma({\cal P}),(\Gamma(B_i)
\hat{\otimes} \Gamma(H))_{i \in I}$ are adapted to $(B,(J_i)_{i \in I})$ and
$({\cal P}, (ker \chi_i)_{i \in I})$ respectively. $\Gamma({\cal P})$ and
$\Gamma(B)$ are given in the following way: One has the extensions
$\chi_{i_{\Omega \rightarrow \Gamma}} : \Omega({\cal P}) \longrightarrow
\Gamma(B_i) \hat{\otimes} \Gamma(H)$ and $\pi_{i_{\Omega \rightarrow \Gamma}}:
\Omega(B) \longrightarrow \Gamma(B_i)$ of the $\chi_i$ and $\pi_i$ respectively.
These extensions form differential ideals $ker \chi_{i_{\Omega \rightarrow
\Gamma}} \subset \Omega({\cal P})$ and $ker \pi_{i_{\Omega \rightarrow \Gamma}}
\subset \Omega(B)$, thus $J({\cal P}):=\cap_{i \in I} ker \chi_{i_{\Omega
\rightarrow \Gamma}}$ and $J(B):=\cap_{i \in I} ker \pi_{i_{\Omega
\rightarrow \Gamma}}$ are differential ideals. By construction, $\Gamma({\cal
P}):=\Omega({\cal P})/J({\cal P})$ and $\Gamma(B):=\Omega(B)/J(B)$ are adapted,
i.e. the extensions $\chi_{i_\Gamma} : \Gamma({\cal P}) \longrightarrow
\Gamma(B_i) \hat{\otimes} \Gamma(H)$ and $\pi_{i_\Gamma}:
\Gamma(B) \longrightarrow \Gamma(B_i)$ of the $\chi_i$ and $\pi_i$
exist and fulfill $\cap_{i \in I} ker \chi_{i_\Gamma}=0$ and $\cap_{i \in
I} ker \pi_{i_\Gamma}=0$ respectively.
\begin{de} \label{strukstrik} A differential structure on a locally trivial 
QPFB is a differential calculus $\Gamma({\cal P})$ defined by a 
family of differential calculi $\Gamma(B_i)$ and a right covariant differential 
calculus $\Gamma(H)$, as described above. \end{de}

\begin{pr} \label{blau} Let $\Gamma({\cal P})$ be a differential structure on 
$\cal P$, and let $\Gamma(B)$ be determined by the corresponding $\Gamma(B_i)$
as above. $\Gamma({\cal P})$ is covariant. The $\chi_{i_\Gamma}$ satisfy \begin{equation}
\label{stoned} \Delta^{\Gamma}_{\cal P}(ker \chi_{i_{\Gamma}}) \subset ker
\chi_{i_{\Gamma}} \otimes H~~\forall i \in I. \end{equation}
The extension $\iota_{\Gamma}: \Gamma(B) \longrightarrow \Gamma({\cal P})$ of 
$\iota$ exists, fulfills \[ \chi_{i_{\Gamma}} \circ \iota_{\Gamma}
(\gamma)=\pi_{i_{\Gamma}}(\gamma) \hat{\otimes} 1,\;\; \forall 
\gamma \in \Gamma(B), \] and is injective. \end{pr}
Proof: As explained before definition \ref{strukstrik}, the differential ideal 
corresponding to $\Gamma({\cal P})$ is 
$J({\cal P})=\cap_i ker \chi_{i_{\Omega \rightarrow \Gamma}} \subset 
\Omega({\cal P})$. Using the right covariance of $\Gamma(H)$ and Definition 
\ref{mocov} one finds that the extensions $\chi_{i_{\Omega \rightarrow \Gamma}}$ 
fulfill \[ (\chi_{i_{\Omega \rightarrow \Gamma}} \otimes id) \circ 
\Delta^{\Omega}_{\cal P} = (id \otimes \Delta^{\Gamma}) \circ 
\chi_{i_{\Omega \rightarrow \Gamma}}, \] where $\Delta^{\Gamma}$ is the 
right coaction of 
$\Gamma(H)$. Due to this formula the differential ideals $ker \chi_{i_{\Omega
\rightarrow \Gamma}}$ are covariant under the coaction of $H$, i.e. 
$\Delta^{\Omega}_{\cal P}(ker \chi_{i_{\Omega \rightarrow \Gamma}}) 
\subset ker \chi_{i_{\Omega 
\rightarrow \Gamma}} \otimes H$, thus, the differential ideal $J({\cal 
P}):=\cap_{i \in I} ker \chi_{i_{\Omega
\rightarrow \Gamma}}$ corresponding to $\Gamma({\cal P})$ is covariant 
and it follows that $\Gamma({\cal P})$ is covariant. This also gives
(\ref{stoned}). 
\\The differential ideal corresponding to $\Gamma(B)$ is 
$J(B)=\cap_i ker \pi_{i_{\Omega \rightarrow \Gamma}} 
\subset  \Omega(B)$. It is easy 
to see that $\iota_{\Omega}(J(B)) \subset J({\cal P})$, thus the extension
$\iota_{\Gamma}$ of $\iota$ with respect to $\Gamma(B)$ and 
$\Gamma({\cal P})$ exists. Clearly $\iota_{\Gamma}$ satisfies 
\[ \chi_{i_{\Gamma}} \circ \iota_{\Gamma}(\gamma)=\pi_{i_{\Gamma}}(\gamma) 
\hat{\otimes} 1,\;\; \forall \gamma \in \Gamma(B). \]
Because of this formula and $\cap_i ker \pi_{i_{\Gamma}}=0$, 
$\iota_{\Gamma}$ is injective. \hfill$\Boxneu$
\\\\The differential structure on a locally trivial QPFB determines the 
covering completion $\Gamma_c({\cal P})$ of $\Gamma({\cal P})$ 
with respect to the covering $(ker \chi_{i_{\Gamma}})_{i \in I}$ (see appendix). 
$\Gamma_c({\cal P})$ is an LC differential algebra (see appendix) with 
local differential calculi $\Gamma(B_i) \hat{\otimes} \Gamma(H)$. It will be 
shown that $\Gamma_c({\cal P})$ is a right $H$-comodule algebra and that the 
covering completion $\Gamma_c(B)$ of $\Gamma(B)$ is embedded in $\Gamma_c({\cal 
P})$. But first we need some facts about differential calculi over $B_{ij} 
\otimes H$ appearing in our context. For the moment we can even assume that 
we have a general differential calculus $\Gamma(B_i \otimes H)$.  
Over the algebras $B_{ij} \otimes H$ there exist two isomorphic
differential calculi $\Gamma^i(B_{ij} \otimes H)= \Gamma(B_i \otimes H)/
\chi_{i_{\Gamma}}(ker \chi_{j_{\Gamma}})$ and $\Gamma^j(B_{ij} \otimes 
H)=\Gamma(B_j \otimes H)/ \chi_{j_{\Gamma}}(ker\chi_{i_{\Gamma}})$, and two
corresponding differential ideals $J^i(B_{ij} \otimes H) \subset
\Omega(B_{ij} \otimes H)$ and $J^j(B_{ij} \otimes H) \subset \Omega(B_{ij} 
\otimes H)$.
\begin{pr} \label{kraut1} The differential ideals $J^i(B_{ij} \otimes H)$ and 
$J^j(B_{ij} \otimes H)$ have the following form: 
\begin{equation} \label{si}
J^i(B_{ij} \otimes H) = (\pi^i_j \otimes id)_{\Omega}(J(B_i \otimes 
H))+\phi_{{ij}_{\Omega}} \circ (\pi^j_i \otimes id)_{\Omega}(J(B_j \otimes
H)) \end{equation}
\begin{equation} \label{si1} J^j(B_{ij} \otimes H) = (\pi^j_i \otimes
id)_{\Omega}(J(B_j \otimes 
H))+\phi_{{ji}_{\Omega}} \circ (\pi^i_j \otimes id)_{\Omega}(J(B_i \otimes
H)), \end{equation}  
where $\phi_{{ij}_{\Omega}}$ are the the extensions of the isomorphisms 
$\phi_{ij}$ corresponding 
to the transition functions $\tau_{ji}$. \end{pr}
For the proof we need 
\begin{lem}\label{kraut} \[ (\pi^i_j \otimes id)\circ \chi_{i} = 
\phi_{ij} \circ (\pi^j_i \otimes id) \circ \chi_{j} \] 
\end{lem}
Proof of the lemma: Using the identities $\phi_{ij}=\chi^i_{ij} \circ 
{\chi^j_{ij}}^{-1}$ and $\chi^i_{ij} \circ \pi_{{ij}_{\cal P}} = 
(\pi^i_j \otimes id) \circ \chi_i$, one has \begin{eqnarray*}
(\pi^i_j \otimes id) \circ \chi_i &=& \chi^i_{ij} \circ 
\pi_{{ij}_{\cal P}}\\
&=& \phi_{ij} \circ \chi^j_{ij} \circ \pi_{{ij}_{\cal P}} \\
&=& \phi_{ij} \circ (\pi^j_i \otimes id) \circ \chi_j. \end{eqnarray*}
\  \hfill$\Boxneu$
\\\\Proof of the proposition: The differential calculus 
$\Gamma^i(B_{ij} \otimes H)=\Gamma(B_i \otimes H)/\chi_{i_{\Gamma}}
(ker\chi_{j_{\Gamma}})$ is isomorphic to $\Gamma({\cal P})/(ker 
\chi_{i_{\Gamma}} + ker \chi_{j_{\Gamma}})$, which in turn is isomorphic to
$\Omega({\cal P})/(ker \chi_{i_{\Omega \rightarrow \Gamma}} + 
ker \chi_{j_{\Omega \rightarrow 
\Gamma}})$. Thus the differential calculi $\Gamma^i(B_{ij} \otimes H)$ and
$\Gamma^j(B_{ij} \otimes H)$ can be identified with 
$\Omega(B_i \otimes H)/\chi_{i_{\Omega}}(ker \chi_{i_{\Omega \rightarrow
\Gamma}} + ker \chi_{j_{\Omega \rightarrow \Gamma}})$ and 
$\Omega(B_j \otimes H)/\chi_{j_{\Omega}}(ker \chi_{i_{\Omega \rightarrow
\Gamma}} + ker 
\chi_{j_{\Omega \rightarrow \Gamma}})$ respectively, and one obtains the  
differential ideals \begin{eqnarray*} 
J^i(B_{ij} \otimes H) &=& (\pi^i_j \otimes id)_{\Omega} \circ
\chi_{i_{\Omega}}(ker 
\chi_{i_{\Omega \rightarrow \Gamma}} + ker \chi_{j_{\Omega \rightarrow 
\Gamma}}), \\
J^j(B_{ij} \otimes H)&=& (\pi^j_i \otimes id)_{\Omega} \circ
\chi_{j_{\Omega}}(ker 
\chi_{i_{\Omega \rightarrow \Gamma}} + ker \chi_{j_{\Omega \rightarrow 
\Gamma}}). \end{eqnarray*}
Now, $\chi_{i_{\Omega}}(ker \chi_{i_{\Omega \rightarrow
\Gamma}})=J(B_i 
\otimes H)$ and $\chi_{j_{\Omega}}(ker \chi_{j_{\Omega \rightarrow 
\Gamma}})=J(B_j \otimes H)$ yields \begin{equation} \label{kraukru}
 J^i(B_{ij} \otimes H) = (\pi^i_j
\otimes id)_{\Omega}(J(B_i \otimes H) + (\pi^i_j \otimes id)_{\Omega} \circ 
\chi_{i_{\Omega}}(\chi^{-1}_{j_{\Omega}}(J(B_j \otimes H))). \end{equation} 
Due to Lemma \ref{kraut} the two homomorphisms $(\pi^i_j \otimes id)_{\Omega} 
\circ \chi_{i_{\Omega}} : \Omega({\cal P}) \longrightarrow \Omega(B_{ij}
\otimes H)$ and $ (\pi^j_i \otimes id)_{\Omega} \circ \chi_{j_{\Omega}} :
\Omega({\cal 
P}) \longrightarrow \Omega(B_{ij} \otimes H)$ are connected by 
\[ (\pi^i_j \otimes id)_{\Omega} \circ \chi_{i_{\Omega}} = \phi_{{ij}_{\Omega}} \circ 
(\pi^j_i \otimes id)_{\Omega} \circ \chi_{j_{\Omega}}, \] thus 
\[ (\pi^i_j \otimes id)_{\Omega} \circ
\chi_{i_{\Omega}}(\chi^{-1}_{j_{\Omega}}(J(B_j \otimes H)))= 
\phi_{{ij}_{\Omega}} \circ (\pi^j_i \otimes id)_{\Omega}(J(B_j \otimes H)).\]
Inserting this formula in (\ref{kraukru}) gives (\ref{si}). (\ref{si1}) results 
by exchanging $i,j$. \hfill$\Boxneu$
\\\\Due to $J^i(B_{ij} \otimes H)=\phi_{{ij}_{\Omega}}(J^j(B_{ij}
\otimes H))$ (immediate from Proposition \ref{kraut1}) the isomorphism 
$\phi_{ij}$ is differentiable with respect to
$\Gamma^j(B_{ij} \otimes H)$ and $\Gamma^i(B_{ij} \otimes H)$.
\\From now on we consider the case $\Gamma(B_i \otimes H)=\Gamma(B_i) \hat{\otimes} 
\Gamma(H)$.
\\Denoting by $(\pi^i_j \otimes id)_{\Gamma^i} : \Gamma(B_i)
\hat{\otimes}\Gamma(H) \longrightarrow \Gamma^i(B_{ij}
\otimes H)$ the natural projection, $\Gamma_c({\cal P})$ has the following
explicit form: \begin{equation} \label{glue} \Gamma_c({\cal P})=
\{ ( \gamma_i)_{i \in I} \in \bigoplus_{i \in I} \Gamma(B_i) \hat{\otimes} 
\Gamma(H) | (\pi^i_j \otimes id)_{\Gamma^i} (\gamma_i)=\phi_{{ij}_{\Gamma}}
\circ (\pi^j_i \otimes id)_{\Gamma^j}(\gamma_j) \}. \end{equation}
Remark: Later we will need \begin{equation} \label{krukru}
ker (\pi^i_j \otimes id)_{\Gamma}=\chi_{i_{\Gamma}}(ker \chi_{j_{\Gamma}})=
\chi_{i_{\Gamma_c}}(ker \chi_{j_{\Gamma_c}}) \end{equation}
\begin{pr} Let $\Gamma({\cal P})$ be a differential 
structure on $\cal P$, let $\Gamma_c({\cal P})$ be the covering completion
of $\Gamma(\cal P)$ and let $\Gamma_c(B)$ be the covering completion of 
$\Gamma(B)$. Let $\chi_{i_{\Gamma_c}}$ and $\pi_{i_{\Gamma_c}}$ be
the restrictions of the respective i-th projections. 
\\Then there exist a unique right coaction $\Delta^{\Gamma_c}_{\cal P} :
\Gamma_c({\cal P}) \longrightarrow \Gamma_c({\cal P}) \otimes H$ and a 
unique injective homomorphism 
$\iota_{\Gamma_c} :\Gamma_c(B) \longrightarrow \Gamma_c({\cal P})$ 
such that
\begin{eqnarray} \label{krakra} (\chi_{i_{\Gamma_c}} \otimes id) \circ 
\Delta^{\Gamma_c}_{\cal P} &=& (id_i \otimes \Delta^{\Gamma}) \circ 
\chi_{i_{\Gamma_c}} \\ \chi_{i_{\Gamma_c}} \circ \iota_{\Gamma_c}(\gamma)
&=& \pi_{i_{\Gamma_c}}(\gamma) \otimes 1,\;\;\forall \gamma \in
\Gamma_c(B). 
\end{eqnarray} \end{pr}
Remark: Indeed, the $\chi_{i_{\Gamma_c}}: \Gamma_c({\cal P}) \longrightarrow
\Gamma(B_i) \hat{\otimes} \Gamma(H)$ and $\pi_{i_{\Gamma_c}}: \Gamma_c(B) 
\longrightarrow \Gamma(B_i)$ coincide with the differential extensions of $\chi_i$
and $\pi_i$. 
\\\\Proof: The covariance of the
ideals $\chi_{i_{\Gamma}}(ker \chi_{j_{\Gamma}})$ under the $H$-coaction
$(id_i \otimes \Delta^{\Gamma})$ follows 
from the covariance of the ideals $ker \chi_{i_{\Gamma}}$ under 
the $H$-coaction $\Delta^{\Gamma}_{\cal P}$. 
Therefore there exist $H$-coactions $(id \otimes  
\Delta)^{\Gamma^i}$ on $\Gamma^i(B_{ij} \otimes H)$ satisfying	
\begin{eqnarray} \label{coi}
(id \otimes \Delta)^{\Gamma^i} \circ (\pi^i_j \otimes id)_{\Gamma^i} &=&
((\pi^i_j \otimes id)_{\Gamma^i})\otimes id) \circ (id \otimes
\Delta^{\Gamma}), 
\\ \label{coi1} (id \otimes \Delta)^{\Gamma^i} &=& (\phi_{{ij}_{\Gamma}} 
\otimes id)
\circ (id \otimes \Delta)^{\Gamma^j}. \end{eqnarray} 
Thus there exists a $H$-coaction $\Delta^{\Gamma_c}_{\cal P}$ on 
$\Gamma_c({\cal P})$ defined by 
\begin{equation} \Delta^{\Gamma_c}_{\cal P}(( \gamma_i)_{i \in I})=(  
(id_i \otimes \Delta^{\Gamma})(\gamma_i))_{i \in I},\;\; \forall 
(\gamma_i)_{i \in I} \in \Gamma_c({\cal P}). \end{equation}
Further one defines an injective homomorphism $\iota_{\Gamma_c}:
\Gamma_c(B) \longrightarrow \Gamma_c({\cal P})$ by 
\begin{equation}
\iota_{\Gamma_c} ((\rho_i)_{i \in I})=( \rho_i \hat{\otimes} 1)_{i
\in I},\;\;\forall ( \rho_i)_{i \in I} \in \Gamma_c(B). \end{equation}
Both homomorphisms are uniquely determined by the assumptions of
the proposition. \hfill$\Boxneu$
\\\\In general the differential calculi $\Gamma^i(B_{ij} \otimes H)$
and $\Gamma^j(B_{ij} \otimes H)$ seem not to be isomorphic to 
differential calculi of the form $\Gamma(B_{ij}) \hat{\otimes }
\Gamma(H)$.
This is suggested by a look at the generators of the differential ideal
$J^i(B_{ij} \otimes H)$: 
\\Let $\iota_{i_{\Omega}} : \Omega(B_i) \longrightarrow \Omega(B_i
\otimes H)$ be the extension of $\iota_i:=id \otimes 1$ and let
$\phi_{i_ {\Omega}}: \Omega(H) \longrightarrow \Omega(B_i \otimes H)$ be
the extension
of $\phi_i:=1 \otimes id$. By Proposition \ref{hat} the
differential  ideal $J(B_i \otimes H)$ corresponding to $\Gamma(B_i)
\hat{\otimes} \Gamma(H)$ is generated by the sets
\[\iota_{i_{\Omega}}(J(B_i)), ~~
 \phi_{i_{\Omega}}(J(H)) \] \begin{equation} \label{genera} \{(a \otimes
 1)d(1
 \otimes h)-(d(1 \otimes h))(a \otimes 1),\;\; a\in B_i,\;
 h\in H\}, \end{equation} where the differential ideals $J(B_i)$ and $J(H)$
 correspond to the differential calculi $\Gamma(B_i)$ and $\Gamma(H)$.
Assume that the differential ideal $J(H)$ is determined 
by a right ideal $R \subset ker \varepsilon \subset H$ in the sense that
$J(H)$ is generated by the set $\{\sum S^{-1}(r_2)dr_1|~r \in R\}$ 
(see also the 
appendix). Using (\ref{phi}),(\ref{genera}) and (\ref{si}) one obtains the
following generators of $J^i(B_{ij} \otimes H)$:
\begin{eqnarray} \label{gener1} && (\pi^i_j \otimes id)_{\Omega} \circ 
\iota_{i_{\Omega}}(J(B_i)),
(\pi^j_i \otimes id)_{\Omega} \circ \iota_{j_{\Omega}}(J(B_j)), 
\\ \label{gener2} && \{ \sum (1 \otimes S^{-1}(r_2)d(1
\otimes r_1)| r\in R \}, 
\\ && \label{gener3} \{ \sum (\tau_{ij}(r_4) \otimes
S^{-1}(r_3)d(\tau_{ji}(r_1) 
\otimes r_2)| r \in R\} , \\ 
&&\label{gener4} \{ (a \otimes 1)d(1 \otimes h)
-(d(1 \otimes 
h))(a \otimes 1)| a\in B_{ij}, h \in H \} , \end{eqnarray} \begin{equation}
\label{gener5} \{(a \otimes 1)d(\sum 
\tau_{ji}(h_1) \otimes h_2)-(d(\sum \tau_{ji}(h_1) \otimes h_2))(a \otimes
1) | a \in B_{ij}, h \in H \}. \end{equation}
Observe that 
\begin{eqnarray*} && \sum (\tau_{ij}(r_4) \otimes S^{-1}(r_3)d(\tau_{ji}(r_1) 
\otimes r_2) - 
\sum (\tau_{ij}(r_2) \otimes 1)d(\tau_{ji}(r_1) \otimes 1) \\
&&- \sum (\tau_{ji}(S(r_4)r_1) \otimes S^{-1}(r_3))d(1 \otimes
r_2) \in J^i(B_{ij} \otimes H);~~ r \in R, \end{eqnarray*} thus one can replace 
the generators (\ref{gener3}) by \begin{equation}  \label{kuh}
\sum (\tau_{ij}(r_2) \otimes 1)d(\tau_{ji}(r_1) \otimes 1) + 
\sum (\tau_{ji}(S(r_4)r_1) \otimes S^{-1}(r_3))d(1
r_2) \in J^i(B_{ij} \otimes H);~~ r \in R. \end{equation} 
\\Using the Leibniz rule, the fact that the image of $\tau_{ji}$ lies in
the center of $B_{ij}$, and the generators (\ref{gener4}), one can replace
(\ref{gener5}) 
by the set of generators \[ \{(a \otimes 1)d(\tau_{ji}(h) \otimes 
1)-d(\tau_{ji}(h) \otimes 1)(a \otimes 1)|~~a \in B_{ij},~h \in
H\}.\]
\begin{pr} \label{hatot} Let the differential calculus $\Gamma(H)$ be
determined by a right ideal $R \subset ker \varepsilon \subset H$ and let
$\tau_{ji}$ be the transition function corresponding to the isomorphism 
$\phi_{ij}$. Assume that the right ideal has the property \begin{equation} 
\label{jbij} \sum \tau_{ij}(S(r_1)r_3) \otimes r_2 \in B_{ij} \otimes R~ 
\forall r \in R;~\forall i,j \in I. \end{equation}
Then there exist differential ideals $J_m(B_{ij}) \subset \Omega(B_{ij})$
such that 
\[ \Gamma^i(B_{ij} \otimes H)=\Gamma^j(B_{ij} \otimes H) \cong 
(\Omega(B_{ij})/J_m(B_{ij})) \hat{\otimes} \Gamma(H). \] \end{pr}
Proof: Because of (\ref{jbij}) the second term of (\ref{kuh})
lies already in the part of $J^i(B_{ij} \otimes H)$ generated by the set  
(\ref{gener2}), thus $J^i(B_{ij} \otimes H)$ is generated by the sets	
\begin{eqnarray*} &&(\pi^i_j \otimes id)_{\Omega} \circ \iota_{i_{\Omega}}
(J(B_i)),(\pi^j_i \otimes id)_{\Omega} \circ \iota_{j_{\Omega}}(J(B_j)),  
\\ && \{ \sum (1 \otimes S^{-1}(r_2)d(1 \otimes r_1)|
r\in R \}, 
\\ &&\{ \sum (\tau_{ij}(r_2) \otimes 1)d(\tau_{ji}(r_1) \otimes 1)| r  
\in R\} , \\ && \{ (a \otimes 1)d(1 \otimes h) -(d(1
\otimes 
h))(a \otimes 1)| a\in B_{ij},h \in H \} , \\ &&\{(a \otimes 1)d( 
\tau_{ji}(h) \otimes 1)-(d(\tau_{ji}(h) \otimes 1))(a \otimes 1) 
| a \in B_{ij},h \in H \}. \end{eqnarray*}
One can see that the differential ideal $J^i(B_{ij} \otimes H)$ is of the
form
(\ref{set}), where the differential ideal $J_m(B_{ij})$ corresponding to  
$\Omega(B_{ij})/J_m(B_{ij})$ is 
generated by the following sets: \begin{eqnarray} \label{fish1} &&
\pi^i_{j_{\Omega}}(J(B_i)), 
\pi^j_{i_{\Omega}}(J(B_j)), \\ \label{fish2} && \{ \sum \tau_{ji}(r_1)d
\tau_{ij}(r_2)|r \in 
R\},\\ \label{fish3} &&\{ (d \tau_{ji}(h))a -a d \tau_{ji}(h)| a \in
B_{ij};h \in H \}. \end{eqnarray}
Replacing $\tau_{ji}$ with $\tau_{ij}$ we get the same differential
ideal $J_m(B_{ij})$. This is clear because of the relation
$\tau_{ji}(S(h))=\tau_{ij}(h)$ 
and the following calculation. From the identity
\[ \sum \tau_{ij}(r_1)\tau_{ji}(r_2)d(\tau_{ij}(r_3) 
\tau_{ji}(r_4))=0;~ r \in R \] one obtains  \[ \sum 
\tau_{ji}(S(r_1)r_4)\tau_{ji}(r_2)d\tau_{ij}(r_3)+ \sum 
\tau_{ij}(r_1)d(\tau_{ji}(r_2) \in J_m(B_{ij}). \]
Due to (\ref{jbij}) the first term lies already in $J_m(B_{ij})$,
thus  $\{ \sum \tau_{ij}(r_1)d(\tau_{ji}(r_2)| r \in R \} \subset
J_m(B_{ij})$. \hfill$\Boxneu$
\\\\Remark: All right ideals $R$ determining a bicovariant differential
calculus 
$\Gamma(H)$ have the property (\ref{jbij}), because such right ideals are 
Ad-invariant, i.e. $\sum S(r_1)r_3 \otimes r_2 \in H \otimes R;~ \forall r
\in R$.
\\\\Observe that in the case described in the previous proposition the 
differential ideal $J_m(B_{ij})$ is in general larger than the differental 
ideal $J(B_{ij})$ (see (\ref{set})), thus the differential calculi 
$\Gamma_m(B_{ij}):=\Omega(B_{ij})/J_m(B_{ij})$ and $\Gamma(B)/
(ker \pi_{i_{\Gamma}} + ker 
\pi_{j_{\Gamma}})$ are in general not isomorphic. 
\\This gives rise to the differential algebra  
\[ \Gamma_m(B):= \{ (\gamma_i)_{i \in I} \in \bigoplus_i \Gamma(B_i)| 
\pi^i_{j_{\Gamma_m}}(\gamma_i \hat{\otimes} 1)=\pi^j_{i_{\Gamma_m}} 
(\gamma_j \hat{\otimes} 1) \}, \] where the homomorphism $\pi^j_{i_{\Gamma_m}}: 
\Gamma(B_i) \longrightarrow \Gamma_m(B_{ij})$ are the composition of the map
$\Gamma(B_{ij}) \longrightarrow \Gamma_m(B_{ij})$ induced by the embedding
$J(B_{ij}) \subset J_m(B_{ij})$ and  $\pi^ i_{j_{\Gamma}}$. Because of
$J(B_{ij}) \subset J_m(B_{ij})$ the LC differential algebra $\Gamma_c(B)$ is a
subalgebra of $\Gamma_m(B)$. Further, $\Gamma_m(B)$ is an LC-differential
algebra naturally embedded in  $\Gamma_c({\cal P})$ by $(\gamma_i)_{i \in I}
\longrightarrow (\gamma_i \otimes  1)_{i \in I}$. If (\ref{jbij}) is fulfilled
one has the identity \[ (\pi^i_j  \otimes id)_{\Gamma^i}=\pi^i_{j_{\Gamma_m}}
\otimes id.\] \\If the right ideal $R$ determining $\Gamma(H)$ does not fulfill
(\ref{jbij}), one can nevertheless construct such  a LC-differential algebra
$\Gamma_m(B)$ with $\Gamma_c(B)$ as subalgebra, and this LC-differential 
algebra on $B$ will play the role of a differential structure on $B$
uniquely induced from the differential structure on $\cal P$.
For an equivalent definition of this LC-differential algebra, we need the
following remark about 
the differential calculus induced on a subalgebra:\\
Let $C$ be an algebra and let $A \subset C$ be a subalgebra. From a 
differential calculus $\Gamma(C)$ one obtains a differential calculus
$\Gamma(A)$ 
by \[ \Gamma^n(A):= \{ \sum_k a^k_0da^k_1...da^k_n \in \Gamma(C)| a^k_i \in
A \}.\] Let $J(C) \subset \Omega(C)$ be the differential ideal corresponding
to the differential calculus $\Gamma(C)$. It is easy to verify that the
differential ideal $J(A) \subset \Omega(A)$ corresponding to $\Gamma(A)$ is
$J(C) \cap \Omega(A)$. \\Now recall that
there are differential calculi $\Gamma^i(B_{ij} \otimes H)$ and
$\Gamma^j(B_{ij}
\otimes H)$. Since $B_{ij} \otimes 1$ is a subalgebra of $B_{ij} \otimes
H$ we obtain differential calculi $\Gamma^i(B_{ij})$ and $\Gamma^j(B_{ij})$,
with  corresponding differential ideals $J^i(B_{ij})$ and $J^j(B_{ij})$
defined by \begin{eqnarray*} J^i(B_{ij}) &=& J^i(B_{ij} \otimes H) \cap
\Omega(B_{ij} \otimes 
1), \\ J^j(B_{ij}) &=& J^j(B_{ij} \otimes H) \cap \Omega(B_{ij} \otimes
1). \end{eqnarray*}
Since $\phi_{{ij}_{\Omega}}(J^j(B_{ij} \otimes H)) =J^i(B_{ij} \otimes H)$
one concludes the identity $\phi_{{ij}_{\Omega}}(J^j(B_{ij}))=J^i(B_{ij})$,
and because of $\phi_{ij}(a \otimes 1)=a \otimes 1$ it follows that
$J^i(B_{ij})=J^j(B_{ij})$, i.e. $\Gamma^i(B_{ij})=\Gamma^j(B_{ij})=
\Gamma_m(B_{ij})$. There are injective homomorphisms $\iota^i_{{ij}_{\Gamma_m}}: 
\Gamma_m(B_{ij}) \longrightarrow \Gamma^i(B_{ij} \otimes H)$ given by
\begin{equation} \label{dadumm} \iota^i_{{ij}_{\Gamma_m}} 
(a_0da_1 d a_2 ... da_n)=
(a_0 \otimes 1) d(a_1 \otimes 1) d (a_2 \otimes 1) ... d(a_n \otimes 1). 
\end{equation} One has the idenitity 
\begin{equation} \label{dadadu} \iota^i_{{ij}_{\Gamma_m}} = 
\phi_{{ij}_{\Gamma}} \circ \iota^j_{{ij}_{\Gamma_m}}.\end{equation}
Let us define the projections $\pi^i_{j_{\Gamma_m}} : \Gamma(B_i) 
\longrightarrow \Gamma_m(B_{ij})$ and $\pi^j_{i_{\Gamma_m}} : \Gamma(B_j) 
\longrightarrow \Gamma_m(B_{ij})$ by \begin{eqnarray}
\label{Bij} \iota^i_{{ij}_{\Gamma_m}} \circ
\pi^i_{j_{\Gamma_m}}(\gamma_i) &=& (\pi^i_j \otimes
id)_{\Gamma^i} (\gamma_i \hat{\otimes} 1), \;\gamma_i \in \Gamma(B_i),\\
\label{Bji} \iota^j_{{ij}_{\Gamma_m}} \circ
\pi^j_{i_{\Gamma_m}}(\gamma_j) &=& (\pi^j_i \otimes
id)_{\Gamma^j} (\gamma_j \hat{\otimes} 1),\; \gamma_j \in\Gamma(B_j). 
\end{eqnarray} 
Obviously, these projections are extensions of $\pi^i_j$ and $\pi^j_i$ 
respectively. In terms of these projections the LC-differential
algebra $\Gamma_m(B)$ is defined as \begin{equation}  \label{glb}
\Gamma_m(B):=\{ (\gamma_i)_{i \in I} \in \bigoplus_{i \in I} \Gamma(B_i)| 
\pi^i_{j_{\Gamma_m}}(\gamma_i)=\pi^j_{i_{\Gamma_m}}(\gamma_j)\}.
\end{equation}
$\Gamma_c(B)$ is a subalgebra of $\Gamma_m(B)$, and there 
exists an injective homomorphism $\iota_{\Gamma_m} : \Gamma_m(B)
\longrightarrow 
\Gamma_c({\cal P})$ defined by \[ \iota_{\Gamma_m}(( \gamma_i)_{i \in 
I})=(\gamma_i \hat{\otimes} 1 )_{i \in I}. \]
\\\\Example: \\In this example we consider a $U(1)$ bundle over the sphere 
$S^2$. Assume that the algebra of differentiable functions $C^{\infty}
(U(1))$ 
over $U(1)$ is the closure in some Fr\'echet topology of the algebra 
generated by the elements $\alpha$ and $\alpha^*$ satisfying
\[\alpha \alpha^*= \alpha^* \alpha =1.\] With
$\Delta(\alpha)= \alpha \otimes \alpha$, $\varepsilon(\alpha)=1$ and 
$S(\alpha)=\alpha^*$,  
this is a Hopf algebra. Let $U_N$ and $U_S$ be the (closed) northern 
and the southern hemisphere respectively, $\{U_N,U_S\}$ is a
covering of $S^2$. We have a complete covering $(I_N,I_S)$ of
$C^{\infty}(S^2)$, $I_N \subset C^{\infty} (S^2)$ and $I_S \subset
C^{\infty}(S^2)$ being the functions vanishing on the subsets $U_N$ and
$U_S$ respectively. Elements of 
$C^{\infty}(U_N)=C^{\infty}(S^2)/I_N$ and
$C^{\infty}(U_S)=C^{\infty}(S^2)/I_S$ can be identified with restrictions 
of elements of $C^{\infty}(S^2)$ to the subsets $U_N$ and $U_S$
respectively. Since $U_N
\cap U_S=S^1$, a transition function $\tau_{NS} :
C^{\infty}(U(1)) 
\longrightarrow C^{\infty}(S^1)$ defines a locally trivial QPFB $\cal P$. We 
choose	
\begin{eqnarray*} \tau_{NS}(\alpha)(e^{i \phi})&=&e^{i \phi}, \\
		  \tau_{NS}(\alpha^*)(e^{i \phi})&=& e^{-i \phi} 
\end{eqnarray*} (Hopf bundle).
\\Now we construct a differential structure on this bundle by giving the 
differential calculi\\ $\Gamma(C^{\infty}(U_N))$, $\Gamma(C^{\infty}(U_S))$
and $\Gamma(C^{\infty}(U(1)))$.
$\Gamma(C^{\infty}(U_N))$ and $\Gamma(C^{\infty}(U_S))$ are taken to be the
usual exterior differential 
calculi where the corresponding differential ideals are generated by all  
elements of the form $adb-dba$. For the right covariant differential 
calculus $\Gamma(C^{\infty}(U(1)))$ we assume a noncommutative form. We
choose as
the right ideal $R$ determining $\Gamma(C^{\infty}(U(1)))$ the right ideal
generated by the element \[\alpha +\nu \alpha^*-(1+\nu)1 \] where $0<\nu
\leq 1$. (One obtains the usual exterior differential calculus for $\nu=1$.)
\\Now we are interested in the LC-differential algebra
$\Gamma_m(C^{\infty}(S^2))$ 
coming from this differential structure on $\cal P$ for $q<1$.
\\It is easy to verify that the right ideal $R$ has the property
(\ref{jbij}), 
thus the differential ideal $J_m(C^{\infty}(S^1))$ is generated by the sets
(\ref{fish1})-(\ref{fish3}). The sets of generators (\ref{fish1}) and 
(\ref{fish3}) give the usual exterior differential calculus on $S^1$, but
the
set of generators (\ref{fish2}) leads to $d \phi =q d\phi$, i.e. $d \phi=0$
for $q<1$. One obtains for the LC-differential algebra
$\Gamma_m(C^{\infty}(S^2))$ 
\begin{eqnarray*} \Gamma^0_m(C^{\infty}(S^2)) &=& C^{\infty}(S^2), \\
\Gamma^n_m(C^{\infty}(S^2)) &=& \Gamma^n(C^{\infty}(U_N)) \bigoplus
\Gamma^n(C^{\infty}(U_S));~~n>0. \end{eqnarray*}
\\\\The foregoing considerations suggest the following definition.
\begin{de} Let $\Gamma({\cal P})$ be a differential structure on the the locally 
trivial QPFB $\cal P$. An LC-differential algebra $\Gamma_g(B)$ over $B$ is 
called embeddable into $\Gamma_c({\cal P})$ if the local differential 
calculi of $\Gamma_g(B)$ are  
$\Gamma(B_i)$ and if there exists the extension $\iota_{\Gamma_g}:  
\Gamma_g(B) \longrightarrow \Gamma_c({\cal P})$ of $\iota$ such 
that \begin{equation} \label{rotesauce} \chi_{i_{\Gamma_c}} \circ \iota_{\Gamma_g}(\gamma) = 
\pi_{i_{\Gamma_g}}(\gamma) \hat{\otimes} 1;~~\forall \gamma \in \Gamma_g(B) 
\end{equation} ($\pi_{i_{\Gamma_g}}: \Gamma_g(B) 
\longrightarrow \Gamma(B_i)$ is the extension of $\pi_i$). \end{de}
Remark: From $\bigcap_{i \in I} ker \pi_{i_{\Gamma_g}}=\{0\}$ follows 
immediately that $\iota_{\Gamma_g}$ is injective.
\begin{pr} The LC-differential algebra $\Gamma_m(B)$ defined above is the 
maximal embeddable LC-differential algebra, i.e every 
embeddable LC-differential 
algebra $\Gamma_g(B)$ is embedded in $\Gamma_m(B)$ as a subalgebra of the direct 
sum of the $\Gamma(B_i)$ by $\gamma \longrightarrow 
(\pi_{i_{\Gamma_g}}(\gamma))_{i \in I}$. \end{pr}
Proof: Let $\Gamma_g(B)$ be an embeddable LC-differential algebra. It is clear 
from $\bigcap_{i \in I} ker \pi_{i_{\Gamma_g}}=\{0\}$ that the mapping is 
injective. To show that its image is in $\Gamma_m(B)$ one has to prove that 
for $\gamma \in \Gamma_g(B)$ 
\begin{equation} \label{chili} \pi^i_{j_{\Gamma_m}} \circ \pi_{i_{\Gamma_g}}
(\gamma)=\pi^j_{i_{\Gamma_m}} \circ \pi_{j_{\Gamma_g}}(\gamma) \end{equation} 
(see (\ref{glb})). By (\ref{rotesauce}) $\iota_{\Gamma_g}$ has the form 
\[\iota_{\Gamma_g}(\gamma)=(\pi_{i_{\Gamma_g}}(\gamma) \hat{\otimes} 1)_{i \in 
I} ;~~\forall \gamma \in \Gamma_g(B).\] By definition, the image of 
$\iota_{\Gamma_g}$ lies in $\Gamma_c({\cal P})$, i.e. 
\[ (\pi^i_j \otimes id)_{\Gamma}(\pi_{i_{\Gamma_g}}(\gamma) \hat{\otimes} 1)
=\phi_{{ij}_{\Gamma}} \circ (\pi^j_i \otimes id)_{\Gamma}(\pi_{j_{\Gamma_g}}
(\gamma) \hat{\otimes} 1).\]
Using (\ref{dadadu}), (\ref{Bij}) and (\ref{Bji}) one obtains (\ref{chili}). 
\hfill$\Boxneu$
\section{\label{con}Covariant derivatives and connections on locally 
trivial QPFB}
First we define covariant derivatives, which are more general objects then
connections. This is done on the covering completion  of the differential
structure on $\cal P$, which is necessary to obtain a one to one
correspondence between covariant derivatives on $\cal P$ and certain  families
of covariant derivatives on the trivializations of $\cal P$. \begin{de}
\label{kovohnz} Let $\Gamma({\cal P})$ be the differential structure on $\cal
P$ and let  $\Gamma_c({\cal P})$ be the covering completion of $\Gamma({\cal
P})$. Let  $hor \Gamma_c({\cal P}) \subset \Gamma_c({\cal P})$ be the
subalgebra defined  by \begin{equation} hor\Gamma_c({\cal P}):=\{ \gamma \in
\Gamma_c({\cal P})| \chi_{i_{\Gamma_c}}(\gamma) \in \Gamma(B_i)
\hat{\otimes} H\; \forall i \in I\}.\end{equation}
A linear map $D_{l,r}: hor\Gamma_c({\cal P}) 
\longrightarrow hor \Gamma_c({\cal P})$ ist called left (right) covariant 
derivative if it satisfies \begin{eqnarray}
\label{nasowas} && D_{l,r}(hor \Gamma^n_c({\cal P})) \subset  
hor \Gamma^{n+1}_c ({\cal P}), \\
\label{kov1} && D_{l,r}(1)=0, \\
\label{kov2} && D_l(\iota_{\Gamma_c}(\gamma) \alpha) = (d(\iota_{\Gamma_c}
\gamma)) \alpha + (-1)^n \gamma D_l(\alpha);~~\gamma \in
\Gamma^n_c(B);~\alpha \in hor \Gamma_c({\cal P}),\\
\label{kov5} && D_r(\alpha \iota_{\Gamma_c} (\gamma)) = D_r(\alpha) 
\iota_{\Gamma_c}(\gamma) + (-1)^n 
\alpha (d\iota_{\Gamma_c} (\gamma));~~\gamma \in
\Gamma_c(B);~\alpha \in hor \Gamma^n_c({\cal P}),\\
\label{kov3} &&(D_{l,r} \otimes id) \circ \Delta_{{\cal P}^{\Gamma_c}} =
\Delta_{{\cal P}^{\Gamma_c}} \circ D_{l,r}, \\
\label{kov4} &&D_{l,r}(ker \chi_{i_{\Gamma_c}}|_{hor \Gamma_c({\cal P})}) 
\subset ker \chi_{i_{\Gamma_c}}|_{hor\Gamma_c({\cal P})};~~\forall i \in I . 
\end{eqnarray}  \end{de}
In this definition the lower indices $l$ or $r$ indicate the left or the right 
case. The appearence $l,r$ means that the corresponding condition is fulfilled 
for both the left and the right case. This convention will be used in the 
sequel permanently.
\\\\Remark: In the case of trivial bundles $B \otimes H$ with differential
structure  $\Gamma(B) \hat{\otimes} \Gamma(H)$, where $hor (\Gamma(B)
\hat{\otimes}  \Gamma(H)) =\Gamma(B) \hat{\otimes} H$, condition (\ref{kov4})
is trivial.  Condition (\ref{kov2}) (respectively (\ref{kov5}) has the form :
\begin{eqnarray*} D_l(\gamma \hat{\otimes} h)&=& d \gamma \hat{\otimes} h + 
(-1)^n (\gamma \hat{\otimes} 1) D_l (1 \otimes h);~~\gamma \in \Gamma^n(B), \\
D_r(\gamma \hat{\otimes} h)&=& D_r(1 \otimes h) (\gamma \hat{\otimes} 1) +
d \gamma \hat{\otimes} h. \end{eqnarray*}
\begin{pr} \label{kovbij} Left (right) covariant derivatives are in 
bijective correspondence to families of linear maps 
$A_{{l,r}_i} : H \longrightarrow \Gamma^1(B_i)$ with the properties 
\begin{eqnarray} \label{abl1} A_{{l,r}_i} 
(1) &=& 0, \\ \label{abl2} \pi^i_{j_{\Gamma_m}}(A_{{l,r}_i}(h))&=& \sum
\tau_{ij}(h_1)\pi^j_{i_{\Gamma_m}}(A_{{l,r}_j}(h_2)) \tau_{ji}(h_3) +\sum
\tau_{ij}(h_1)d\tau_{ji}(h_2).\end{eqnarray} \end{pr}
Remark: Note that (\ref{abl2}) is a condition in $\Gamma_m(B_{ij})$ (See the 
considerations at the end of the forgoing section.).
\\\\Proof: Because of (\ref{kov4}) a given left covariant derivative on $hor 
\Gamma_c({\cal P})$ determines a family of left covariant derivatives $D_{l_i} : 
\Gamma(B_i) \hat{\otimes} H \longrightarrow \Gamma(B_i) \hat{\otimes} H $ by 
\begin{equation} \label{grob} D_{l_i} \circ \chi_{i_{\Gamma_c}} = 
\chi_{i_{\Gamma_c}} \circ D_l.\end{equation} It follows the identity 
$D_l((\gamma_i)_{i \in I})=(D_{l_i}(\gamma_i))_{i \in I}$. Since 
$(D_{l_i}(\gamma_i))_{i \in I} \in \Gamma_c({\cal P})$, the $D_{l_i}$ satisfy 
\begin{equation} \label{smar} (\pi^i_j \otimes id)_{\Gamma}\circ D_{l_i}
(\gamma_i)= \phi_{{ij}_{\Gamma}} \circ (\pi^j_i \otimes id)_{\Gamma}\circ 
D_{l_j}(\gamma_j),~~ (\gamma_i)_{i \in I} \in hor \Gamma_c({\cal 
P}).\end{equation}
One obtains a family of linear maps $A_{l_i}: 
H \longrightarrow \Gamma^1(B_i)$ by  
\[ A_{l_i}(h):=-(id \otimes \varepsilon) \circ D_{l_i}(1 \otimes h).\]
Now we need: 
\begin{lem} \label{id} \[ ((id \otimes \varepsilon)_{\Gamma} \otimes id)
\circ 
\Delta^{\Gamma}_{\cal P}{\mid_{\Gamma(B) \hat{\otimes} H}}=id \]
\end{lem}
Proof of the lemma: An element $\gamma \in \Gamma(B) \hat{\otimes} H$ has
the general form \[ \gamma = \sum_k a^k_0da^k_1 \hat{\otimes} h^k.\]
We obtain \begin{eqnarray*} ((id \otimes \varepsilon)_{\Gamma} \otimes 
id)\circ \Delta^{\Gamma}_{\cal P}(\gamma) &=& ((id \otimes \varepsilon) 
\otimes id)(\sum_k \sum a^k_0da^k_1 \hat{\otimes} h^k_1 \otimes h^k_2) \\
&=&\sum_k a^k_0da^k_1 \hat{\otimes} h^k.\end{eqnarray*}
{}\hfill$\Boxneu$\\
By the foregoing lemma, (\ref{kov2}) and (\ref{kov3}) one computes the identity 
\begin{equation} \label{jajahaha} D_{l_i}(\gamma \hat{\otimes} h)=d\gamma 
\hat{\otimes} h + (-1)^{n+1} \sum \gamma
A_{l_i}(h_1) \hat{\otimes} h_2;~~\gamma \in \Gamma^n(B_i);~~h \in H. 
\end{equation} 
Because of (\ref{kov1}) the $A_{l_i}$ fulfill (\ref{abl1}).
To prove the property (\ref{abl2}) we need:
\begin{lem} \label{frosch} Let $B$ be an algebra, $H$ be a Hopf algebra, 
$\Gamma(B)$ be a 
differential calculus over $B$ and $\Gamma(H)$ be a right covariant differential 
calculus over $H$. Let $D_l : \Gamma(B) \hat{\otimes} H \longrightarrow 
\Gamma(B) \hat{\otimes} H$ be a left covariant derivative on the trivial bundle 
$B \otimes H$. Let $J \subset \Gamma(B) 
\hat{\otimes} \Gamma(H)$ be a differential ideal with the property $(id \otimes 
\Delta^{\Gamma})(J) \subset J \otimes H.$ Then one has 
\begin{equation} \label{knoch} D_l(J \cap (\Gamma(B) \hat{\otimes} H)) \subset 
J \cap (\Gamma(B) \hat{\otimes} H). \end{equation} \end{lem}
Proof of the Lemma: By Lemma \ref{lemisubb} there is an ideal  
$\tilde{J} \subset \Gamma(B)$ such that 
\[J \cap (\Gamma(B) \hat{\otimes} H) = \tilde{J} \hat{\otimes} H. \]
$\tilde{J}$ is an differential ideal: Let
$\sum_k \gamma_k \hat{\otimes} h_k \in
\tilde{J} \hat{\otimes} H \subset J$. Since $J$ is a differential ideal one obtains
\[ \sum_k d \gamma_k \hat{\otimes} h_k + (-1)^n\sum_k \gamma_k \hat{\otimes}
dh_k \in J;~~\gamma_k \in \Gamma^n(B).\]
The second summand lies in $ \tilde{J} \hat{\otimes} \Gamma^1(H) \subset
J$. It follows that 
$\sum_k d \gamma_k \hat{\otimes} h_k \in d\tilde{J} \hat{\otimes} H
\subset  J \cap (\Gamma(B) \hat{\otimes} H)$ and one obtains  
$d\tilde{J} \subset \tilde{J}$, thus $\tilde{J}$ is a differential ideal.
\\Applying $D_l$ to $\sum_k \gamma_k \hat{\otimes} h_k \in
\tilde{J} \hat{\otimes} H \subset J$ leads to \[ D_l (\sum_k \gamma_k 
\hat{\otimes} h_k)=
\sum_k d\gamma_k \hat{\otimes}h_k +(-1)^{n+1}\sum_k \gamma_k
D_l(1 \otimes h_k),~~ \gamma_k \in \Gamma^n(B).\]
Since the image of $D_l$ lies in $\Gamma(B) \hat{\otimes} H$, the right hand 
side of 
this formula is an element of $\tilde{J} \hat{\otimes} H$. \hfill$\Boxneu$
\\\\Since the $ker (\pi^i_j \otimes id)_{\Gamma} \subset \Gamma(B_i) 
\hat{\otimes} \Gamma(H)$ are coinvariant differential ideals 
(see (\ref{coi})), by the foregoing lemma follows 
$D_{l_i}(ker (\pi^i_j \otimes id)_{\Gamma} \cap (\Gamma(B_i) \hat{\otimes} H)) 
\subset ker (\pi^i_j \otimes id)_{\Gamma} \cap (\Gamma(B_i) \hat{\otimes} H)$.
This allows to define linear maps $D^{ij}_{l_i}$ by \[D^{ij}_{l_i} \circ (\pi^i_j 
\otimes id)_{\Gamma}=(\pi^i_j \otimes id)_{\Gamma} \circ D_{l_i}.\]
Applying $(\pi^i_j \otimes id)_{\Gamma}$ to (\ref{jajahaha}) 
one obtains \begin{equation} \label{zaub} 
D^{ij}_{l_i}(a \otimes h)= (d(a \otimes 
1))(1 \otimes h) -(a \otimes 1) \sum (\pi^i_j \otimes 
id)_{\Gamma}(A_{l_i}(h_1) \otimes 1) (1 \otimes h_2). \end{equation} 
Let $(\gamma_i)_{i \in I} \in hor \Gamma_c({\cal P})$, in particular 
\begin{equation} \label{hex}(\pi^i_j \otimes id)_{\Gamma} (\gamma_i) = 
\phi_{{ij}_{\Gamma}}\circ (\pi^j_i \otimes id)_{\Gamma}(\gamma_j). 
\end{equation} 
Since $D_l(\gamma_i)_{i \in I} = (D_{l_i}(\gamma_i))_{i \in I} \in hor 
\Gamma_c({\cal P})$ it follows that
\begin{equation} \label{hexe} 
D^{ij}_{l_i} \circ (\pi^i_j \otimes id)_{\Gamma}(\gamma_i)=
\phi_{{ij}_{\Gamma}} 
\circ  D^{ij}_{l_j} \circ (\pi^j_i \otimes id)_{\Gamma}(\gamma_j).\end{equation}
Combining (\ref{hex}) and (\ref{hexe}), one obtains 
\begin{equation} \label{bluber} D^{ij}_{l_i} = 
\phi_{{ij}_{\Gamma}} \circ D^{ij}_{l_j} \circ 
\phi_{{ji}_{\Gamma}}.\end{equation} 
Taking advantage of (\ref{zaub}), (\ref{bluber}), (\ref{iotaij}), (\ref{dadadu})
and (\ref{phi}) one computes (see also (\ref{dadumm})) \begin{eqnarray} 
\nonumber D^{ij}_{l_i}(1 \otimes h)&=& -\sum \iota^i_{{ij}_{\Gamma_m}}
(\pi^i_{j_{\Gamma_m}}(A_{l_i}(h_1))) (1 \otimes h_2) \\ \nonumber
&=& \phi_{{ij}_{\Gamma}} \circ D^{ij}_{l_j} \circ
\phi_{{ji}_{\Gamma}}(1 \otimes h)  \\ \label{tja}
&=& \phi_{{ij}_{\Gamma}} \circ D^{ij}_{l_j}(\sum \tau_{ij}(h_1) \otimes h_2) \\
\nonumber
&=& \phi_{{ij}_\Gamma}( \sum \iota^j_{{ij}_{\Gamma_m}}(d\tau_{ij}(h_1))
(1 \otimes h_2) - \sum \iota^j_{{ij}_{\Gamma_m}}
(\tau_{ij}(h_1) \pi^j_{i_{\Gamma_m}}(A_{l_j}(h_2)) (1 \otimes h_3)) \\ \nonumber
&=& \sum \iota^i_{{ij}_{\Gamma_m}}
((d\tau_{ij}(h_1))\tau_{ji}(h_2))(1 \otimes h_3) \\ \nonumber
&&-\sum \iota^i_{{ij}_{\Gamma_m}}
(\tau_{ij}(h_1) (\pi^j_{i_{\Gamma_m}}(A_{l_j}(h_2)) \tau_{ji}(h_3))
(1 \otimes h_4). \end{eqnarray}
Applying the Leibniz rule to the first term of the last row and using 
$\sum \tau_{ij}(h_1)\tau_{ji}(h_2)= \varepsilon(h)1$ one obtains the identity 
\begin{eqnarray}
\label{tran2} \sum \iota^i_{{ij}_{\Gamma_m}}
(\pi^i_{j_{\Gamma_m}}(A_{l_i}(h_1)))(1 \otimes h_2)  
 &=& \sum \iota^i_{{ij}_{\Gamma_m}}
 (\tau_{ij}(h_1)\pi^j_{i_{\Gamma_m}} (A_{l_j}(h_2))\tau_{ji}(h_3))(1 \otimes h_4) 
 \\ \nonumber &&+ \sum \iota^j_{{ij}_{\Gamma_m}}
 (\tau_{ij}(h_1)d(\tau_{ji}(h_2)))(1 \otimes h_3).\end{eqnarray}
In order to arrive at (\ref{abl2}) we need to kill the $1 \otimes h$-factor.
This is achieved by using a projection 
$P_{inv} : \Gamma^i(B_{ij} \otimes H)
\longrightarrow \{ \gamma \in \Gamma^i(B_{ij} \otimes H)|(id \otimes
\Delta)^{\Gamma^i}(\gamma)=\gamma \otimes 1\}$ on the elements of 
$\Gamma^i(B_{ij} \otimes H)$ being coinvariant under the right H coaction 
$(id \otimes \Delta)^{\Gamma^i}: \Gamma^i(B_{ij} \otimes H) \longrightarrow
\Gamma^i(B_{ij} \otimes H) \otimes H$ (see also (\ref{coi}) and (\ref{coi1})).	
$P_{inv}$ is defined by \begin{equation}
\label{pinv} P_{inv}(\rho)=\sum \rho_0S(\rho_1), \; \rho \in 
\Gamma^i(B_{ij} \otimes H). \end{equation}
\\Applying $P_{inv}$ to the identity (\ref{tran2}) leads to 
\begin{eqnarray*} \iota^i_{{ij}_{\Gamma_m}}
(\pi^i_{j_{\Gamma_m}}(A_{l_i}(h))) &=&
\sum \iota^i_{{ij}_{\Gamma_m}} (\tau_{ij}(h_1)\pi^j_{i_{\Gamma_m}}(A_{l_j}(h_2)) 
\tau_{ji}(h_3))\\ &&+ 
\sum \iota^i_{{ij}_{\Gamma_m}} (\tau_{ij}(h_1)d\tau_{ji}(h_2)). 
\end{eqnarray*}
Due to the injectivity of $\iota^i_{{ij}_{\Gamma_m}}$, this is identical to 
\begin{equation} \pi^i_{j_{\Gamma_m}}(A_{l_i}(h))=\sum 
\tau_{ij}(h_1)\pi^j_{i_{\Gamma_m}}(A_{l_j})(h_2))\tau_{ij}(h_3)+ \sum 
\tau_{ij}(h_1)d\tau_{ji}(h_2) \end{equation} in $\Gamma_m(B_{ij})$.
\\\\Now we prove the converse assertion. Assume there is given a family of linear 
maps $A_{l_i} : H \longrightarrow \Gamma(B_i)$ which fulfill (\ref{abl1}) and 
(\ref{abl2}). Every $A_{l_i}$ 
defines by \[ D_{l_i}(\gamma \hat{\otimes} h)=d\gamma
\hat{\otimes} h + (-1)^{n+1} \sum \gamma
A_{l_i}(h_1) \hat{\otimes} h_2;~~\gamma \in \Gamma^n(B_i);~~h \in H \]
a left covariant derivative  
$D_{l_i}$ on $\Gamma(B_i) \hat{\otimes} H$. The properties 
(\ref{nasowas})-(\ref{kov2}) and (\ref{kov3}) of 
$D_{l_i}$, are easily derived from the above formula. One has to show 
that $D_{l}((\gamma_i)_{i \in I}):=(D_{l_i}(\gamma_i))_{i\in I},~
(\gamma_i)_{i \in I} \in hor\Gamma_c({\cal P})$ 
is a covariant derivative on $hor\Gamma_c({\cal P})$.
Because of (\ref{abl1}), $D_l$ fulfills (\ref{kov1}). 
The conditions (\ref{kov2}) and (\ref{kov3}) follows from the corresponding 
properties of $D_{l_i}$. It remains to prove, that the image of $D_l$ lies in 
$\Gamma_c({\cal P})$, because then it also lies in $hor \Gamma_c({\cal P})$. 
(This is due to the 
fact that all the images of the $D_{l_i}$ obviously are in $\Gamma(B_i) 
\hat{\otimes} H$.) Then it is also obvious from the fact that the 
$\chi_{i_{\Gamma_c}}$ are the projections to the $i$-th components that 
condition (\ref{kov4}) is fulfilled. The image of $D_l$ lies in 
$hor\Gamma_c({\cal P})$ if the family of the $D_{l_i}$ fulfills 
\begin{equation} \label{drin} (\pi^i_j \otimes
id)_{\Gamma}\circ D_{l_i}(\gamma_i)=\phi_{{ij}_{\Gamma}}
\circ  (\pi^j_i \otimes id)_{\Gamma} \circ D_{l_j}(\gamma_j), ~\forall 
(\gamma_i)_{i \in I} \in  hor\Gamma_c({\cal P}). \end{equation}
\\By Lemma \ref{frosch}, the covariant derivatives $D_{l_i}$ give rise to  
maps $D^{ij}_{l_i}$ defined by	
\[ D^{ij}_{l_i} \circ (\pi^i_j \otimes id)_{\Gamma} = (\pi^i_j \otimes
id)_{\Gamma} \circ D_{l_i}.\] One has \begin{equation} \label{krak}
D^{ij}_{l_i}(1 \otimes h)=-\sum (\pi^i_j \otimes
id)_{\Gamma}(A_{l_i}(h_1) \otimes 1) (1 \otimes h_2), \end{equation} 
and we will show that (\ref{abl2}) yields the identity 
\[ D^{ij}_{l_i} = \phi_{{ij}_{\Gamma}} \circ D^{ij}_{l_j} 
\circ \phi_{{ji}_{\Gamma}}: \]
One computes for $\gamma \in \Gamma^n_m(B_{ij})$ 
\begin{eqnarray*} D^{ij}_{l_i}(\iota^i_{{ij}_{\Gamma_m}}(\gamma)(1 \otimes h))
&=&(d\gamma) (1 \otimes h)+(-1)^{n+1}\iota^i_{{ij}_{\Gamma_m}}
(\gamma) D^{ij}_{l_i}(1 \otimes h)\\
&=&\iota^i_{{ij}_{\Gamma_m}} 
(d\gamma) (1 \otimes h)+(-1)^{n+1}\sum \iota^i_{{ij}_{\Gamma_m}}
(\gamma \pi^i_{j_{\Gamma_m}} (A_{l_i}(h_1)))(1 
\otimes h_2)\\ &=& \iota^i_{{ij}_{\Gamma_m}} (d\gamma)
(1 \otimes h)\\&&+(-1)^{n+1}\sum \iota^i_{{ij}_{\Gamma_m}}
(\gamma \tau_{ij}(h_1) 
\pi^j_{i_{\Gamma_m}}(A_{l_j}(h_2))\tau_{ji}(h_3))(1 \otimes h_4) \\
&&+ (-1)^{n+1} \sum \iota^i_{{ij}_{\Gamma_m}}
(\gamma (d \tau_{ij}(h_1))\tau_{ji}(h_2))(1 
\otimes h_3) \\
&=& \phi_{{ij}_{\Gamma}} \circ D^{ij}_{l_j} \circ \phi_{{ji}_{\Gamma}}(\gamma(1 
\otimes h)). \end{eqnarray*}
Thus, one obtains for $(\gamma_i)_{i \in I} \in hor \Gamma_c({\cal P})$ 
\begin{eqnarray*} D^{ij}_{l_i} \circ (\pi^i_j \otimes id)_{\Gamma}(\gamma_i) 
&=& D^{ij}_{l_i} \circ \phi_{{ij}_{\Gamma}}\circ (\pi^j_i \otimes 
id)_{\Gamma}(\gamma_j) \\
&=& \phi_{{ij}_{\Gamma}} \circ D^{ij}_{l_j} \circ (\pi^j_i \otimes 
id)_{\Gamma}(\gamma_j), \end{eqnarray*} and (\ref{drin}) follows.
\\It is immediate from the construction (using Lemma \ref{id}) that the 
correspondence is bijective. 
\\The proof for right covariant derivatives is analogous. In this case one uses 
\begin{equation} D_{r_i}(\gamma \hat{\otimes} h)= d \gamma \hat{\otimes}h + 
(-1)^{n+1} \sum A_{r_i}(h_1) \gamma \hat{\otimes} h_2 \end{equation} for 
$\gamma \in \Gamma^n(B_i)$ 
\hfill$\Boxneu$
\\\\Remark: Obviously, a family of linear maps $A_i: H \longrightarrow 
\Gamma^1(B_i)$ fulfilling  (\ref{abl1}) and (\ref{abl2}) determines at the 
same time a left and a right covariant derivative. Consequently, there is also a 
bijective correspondence between left and right covariant derivatives.
\begin{pr} Let $D_{l,r}: hor \Gamma_c({\cal P}) \longrightarrow  hor 
\Gamma_c({\cal P})$ be a left (right) covariant derivative and let $\Gamma_g(B)$ 
be embeddable into $\Gamma_c({\cal P})$. $D_{l,r}$ fulfills 
\begin{eqnarray}\label{kov2g} && 
D_l(\iota_{\Gamma_g}(\gamma) \alpha) = 
(d(\iota_{\Gamma_g}(\gamma)) \alpha + (-1)^n \iota_{\Gamma_g}
(\gamma) D_l(\alpha);~~\gamma \in
\Gamma^n_g(B);~\alpha \in hor \Gamma_c({\cal P}),\\
\label{kov5g} && D_r(\alpha \iota_{\Gamma_g} (\gamma)) = D_r(\alpha)
\iota_{\Gamma_g}(\gamma) + (-1)^n
\alpha (d\iota_{\Gamma_g} (\gamma));~~\gamma \in
\Gamma_g(B);~\alpha \in hor \Gamma^n_c({\cal P}). \end{eqnarray} \end{pr}
Proof: Let $(\gamma_i)_{i \in I} \in hor \Gamma_c({\cal P})$ and $\rho \in 
\Gamma^n_g(B)$. One has $\iota_{\Gamma_g}(\rho)=
(\pi_{i_{\Gamma_g}}(\rho) \hat{\otimes} 1)_{i \in I}$ and 
$D_l((\gamma_i)_{i \in I})=(D_{l_i}(\gamma_i))_{i \in I}$. One calculates 
\begin{eqnarray*} D_l(\iota_{\Gamma_g}(\rho) (\gamma_i)_{i \in I})&=& 
D_l(((\pi_{i_{\Gamma_g}}(\rho) \hat{\otimes} 1)\gamma_i)_{i \in I}) \\
&=&(D_{l_i}(\pi_{i_{\Gamma_g}}(\rho) \hat{\otimes} 1)\gamma_i))_{i \in I} \\
&=& ((d(\pi_{i_{\Gamma_g}}(\rho) \hat{\otimes} 1))\gamma_i))_{i \in I}+ (-1)^n
((\pi_{i_{\Gamma_g}}(\rho) \hat{\otimes} 1)D_{l_i}(\gamma_i))_{i \in I} \\
&=& (d(\iota_{\Gamma_g}(\rho))) (\gamma_i)_{i \in I} + (-1)^n 
\iota_{\Gamma_g}(\rho)	D_l((\gamma_i)_{i \in I}). \end{eqnarray*}
The proof for right covariant derivatives is anlog. \hfill$\Boxneu$
\\\\Now we are going to define connections on locally trivial QPFB. 
It turns out that connections are special cases of covariant derivatives. 
We start with a definition dualizing the classical case in a certain sense.  
\begin{de} \label{zusi} Let $\Gamma({\cal P})$ be a differential structure on 
$\cal P$ and let $\Gamma_c({\cal P})$ be the covering completion of 
$\Gamma(\cal P)$. 
A left (right) connection is a surjective left (right) $\cal P$-module homomorphism 
$hor_{l,r}: \Gamma^1_c{\cal P} \longrightarrow hor \Gamma^1_c({\cal P})$ such 
that:
\begin{equation} \label{quad} {hor_{l,r}}^2=hor_{l,r}, \end{equation}
\begin{equation} \label{krid}
(hor_{l,r} \otimes id) \circ \Delta^{\Gamma_c}_{\cal P} = \Delta^{\Gamma_c}_{\cal P} 
\circ hor_{l,r} \end{equation} and 
\begin{equation} \label{triv1} hor_{l,r}(ker \chi_{i_{\Gamma_c}}) \subset
ker \chi_{i_{\Gamma_c}},~\forall i \in I. \end{equation} \end{de}
Remark: Conditions (\ref{triv1}) in this definition are needed to have 
the one-to-one 
correspondence between connections on $\cal P$ and certain families 
of connections on the trivial bundles $B_i \otimes H$. On a trivial bundle $B 
\otimes H$ condition (\ref{triv1}) is obsolete.
\\\\Remark: For a given left (right) connection there is a vertical left (right) 
$\cal P$-submodule $ver_{l,r} \Gamma^1_c({\cal P})$ such that 
\[\Gamma^1_c({\cal P})=  ver_{l,r} \Gamma^1_c({\cal P}) \oplus 
hor\Gamma^1_c({\cal P}),\] where the projection $ver_{l,r}: 
\Gamma^1_c({\cal P}) \longrightarrow ver_{l,r}\Gamma^1_c({\cal P})$ is 
defined by $ver_{l,r}:=id-hor_{l,r}$.
\\\\On a trivial bundle $B \otimes H$ with differential structure $\Gamma(B) 
\hat{\otimes} \Gamma(H)$ exists always the canonical connection $hor_c$, which 
is at the same time left and right. The existence of $hor_c$ comes from the 
decomposition 
\[ (\Gamma(B) \hat{\otimes} \Gamma(H))^1 =(\Gamma^1(B) \hat{\otimes} H) 
\oplus (B \hat{\otimes} \Gamma^1(H)) \] (direct sum of $(B \otimes H)$-bimodules), 
which allows to define \begin{eqnarray*} hor_c(\gamma \hat{\otimes} h) &=& \gamma 
\hat{\otimes} h;~~\gamma \in \Gamma^1(B),~h \in H, \\ 
hor_c(a \hat{\otimes} \theta)&=& 0;~~a \in B,~\theta \in \Gamma^1(H). 
\end{eqnarray*}
\begin{lem} For a given connection $hor_{l,r}$ on $\cal P$ there exists a family 
of connections $hor_{{l,r}_i}$ on the trivilizations 
$B_i \otimes H$ such that 
\begin{equation}
\label{hori1} 
\chi_{i_{\Gamma_c}} \circ hor_{l,r} = hor_{{l,r}_i} 
\circ \chi_{i_{\Gamma_c}} . \end{equation} \end{lem}
Proof: The existence of linear map $hor_{l_i}$ satisfying 
(\ref{hori1}) follows from (\ref{triv1}). The $hor_{{l,r}_i}$ are connections 
on $B_i \otimes H$:
Because of the surjectivity of the $\chi_{i_{\Gamma_c}}$ the $hor_{l_i}$ map 
onto $\Gamma^1(B_i) \hat{\otimes} H$. To prove condition (\ref{quad}) 
one computes
\begin{eqnarray*}  
{hor_{{l,r}_i}}^2 \circ \chi_{i_{\Gamma_c}} 
&=&hor_{{l,r}_i} \circ \chi_{i_{\Gamma_c}} \circ hor_{l,r} \\ 
&=&\chi_{i_{\Gamma_c}} \circ {hor_{l,r}}^2 \\
&=&\chi_{i_{\Gamma_c}} \circ hor_{l,r}\\
&=&hor_{{l,r}_i} \circ \chi_{i_{\Gamma_c}} \end{eqnarray*} 
The condition (\ref{krid}) is fulfilled because of (\ref{krakra}). \hfill$\Boxneu$
\\By  Definition \ref{zusi} and the foregoing lemma a connection $hor_{l,r}$ has 
the following form \begin{equation} \label{hoho} hor_{l,r}((\gamma_i)_{i\in I})
=(hor_{{l,r}_i}(\gamma_i))_{i \in I}, ~ (\gamma_i)_{i\in I} \in hor 
\Gamma_c({\cal P}), \end{equation}
which also means that the family of linear maps $hor_{i_{l_r}}$ satisfies 
\begin{equation} \label{joe}  (\pi^i_j \otimes id)_{\Gamma} \circ 
hor_{{l,r}_i}(\gamma_i) = \phi_{{ij}_{\Gamma}} \circ 
(\pi^j_i \otimes id)_{\Gamma} \circ hor_{{l,r}_j} (\gamma_j) \end{equation}
for $(\gamma_i)_{i \in I} \in \Gamma^1_c({\cal P})$.
\begin{pr} \label{vi} Let $R \subset H$ be the right ideal corresponding
to the right covariant differential calculus $\Gamma(H)$. 
Left (right) connections on a locally trivial QPFB $\cal P$ are in 
one-to-one correspondence to left (right) covariant derivatives with the following 
property: The corresponding linear maps $A_{{l,r}_i}$ fulfill 
\begin{eqnarray} \label{kerrl}	R \subset ker A_{l_i},~~\forall i \in I \\
\label{kerrr} S^{-1}(R) \subset ker A_{r_i},~~\forall i \in I. \end{eqnarray}
\end{pr}
Remark: Thus left (right)connections are in one to one correspondence to linear 
maps $A_{{l,r}_i}$ fulfilling (\ref{abl1}), (\ref{abl2}) and (\ref{kerrl}) 
(respectively (\ref{kerrr})).
\\\\Proof: We prove the assertion only for left connections. The proof is fully 
analogous for right connections. 
\\\\A left connection $hor_l$ determines a family of linear maps $A_{l_i} : 
H \longrightarrow \Gamma^1(B_i)$ by 
\[ A_{l_i}(h):= -(id \otimes \varepsilon ) hor_{l_i}(1 \hat{\otimes} dh).\] 
From (\ref{krid}) and Lemma \ref{id} one has te identity 
\begin{equation} \label{soeinhor} hor_{l_i} (1 \otimes dh)=-\sum A_{l_i}(h_1) 
\hat{\otimes} h_2. \end{equation} Therefore $A_{l_i}$ have the property 
$R \subset ker A_{l_i}$: \[ 0 = \sum hor_{l_i}(1 \hat{\otimes} 
S^{-1}(r_2)dr_1)=-A_{l_i}(r) \hat{\otimes}1,~~r \in R.\]  
It remains to show that this family of linear maps fulfills (\ref{abl1}) and 
(\ref{abl2}).
\\(\ref{abl1}) is fulfilled by definition ( $A_{l_i}(1):=(id \otimes 
\varepsilon)\circ hor_{l_i} (1 \otimes d1)=0$).
\\Because of (\ref{triv1}), (\ref{hori1}) and (\ref{krukru}) one has $hor_{l_i} 
(ker (\pi^i_j \otimes id)_{\Gamma}) \subset  ker (\pi^i_j \otimes id)_{\Gamma}$, 
and the linear maps $hor^{ij}_{l_i}$ defined by \[  hor^{ij}_{l_i} \circ
(\pi^i_j  \otimes id)_{\Gamma}=(\pi^i_j \otimes id)_{\Gamma} \circ hor_{l_i}
\] exist.  It follows that \[hor^{ij}_{l_i}(d(1 \otimes h))=
-\sum \pi^i_{j_{\Gamma_m}}(A_{l_i}(h_1)) (1 \otimes h_2). \] 
On the other hand, by an analogoue of the computation leading to (\ref{bluber}) 
(using (\ref{joe})), one obtains 
\[ hor^{ij}_{l_i}= \phi_{{ij}_{\Gamma}} \circ hor^{ij}_{l_j} 
\circ \phi_{{ji}_{\Gamma}}.\] Now using the last two formulas, on can repeat the 
arguments written after formula (\ref{bluber}) to obtain formula (\ref{abl2}).
\\\\Now assume that there is given a left covariant derivative $D_l$, whose 
corresponding linear maps $A_{l_i}$ satisfy $R \subset ker A_{l_i}$. There 
exist left connections $hor_{l_i}: (\Gamma^1(B_i) \hat{\otimes} H) \oplus (B_i 
\hat{\otimes} \Gamma^1(H)) \longrightarrow \Gamma^1(B_i) \hat{\otimes} H $ 
defined by \begin{eqnarray} \nonumber
hor_{l_i}(\gamma \hat{\otimes} h)&:=& \gamma \hat{\otimes}h,\\ 
\label{nasoeinhor1} hor_{l_i}( a \hat{\otimes} hdk) &:=& -\sum a A_{l_i}(k_1) 
\hat{\otimes}hk_2. \end{eqnarray} To verify this assertion we define linear 
maps $hor^{\Omega}_{l_i}: (\Gamma(B_i) \hat{\otimes} 
\Omega(H))^1 \longrightarrow \Gamma^1(B_i) \hat{\otimes} H$ by 
\begin{eqnarray*} hor^{\Omega}_{l_i}(a_0da_1 \hat{\otimes} h) &=& a_0da_1
\hat{\otimes} h, \\ hor^{\Omega}_{l_i}(a \hat{\otimes} h^0dk)&=&-\sum
aA_l(k_1) \hat{\otimes} hk_2. \end{eqnarray*}
The $(B_i \otimes H)$ subbimodules $B_i \hat{\otimes} J^1(H)$ are generated by
the sets $\{1 \hat{\otimes} \sum S^{-1}(r_2)dr_1| r\in R\}$. One has 
\[B_i \hat{\otimes} \Gamma^1(H)=(B_i \hat{\otimes}\Omega^1(H))/(B_i 
\hat{\otimes} J^1(H))= B_i \hat{\otimes} \Omega^1(H)/J^1(H).\]
Using $R \subset ker A_{l_i}$ it is easy to verify that the linear
maps $hor^{\Omega}_{l_i}$ sends $B_i \hat{\otimes} J^1(H)$
to zero, i.e. there exist corresponding the linear 
maps $hor_{l_i}$ on $(\Gamma(B_i) \hat{\otimes} \Gamma(H))^1$. As a  
consequence of there definition these linear maps are connections. One easily 
verifies the identity 
\begin{equation} \label{hord} hor_{l_i} \circ d = D_{l_i}|_{B_i \otimes H},
\end{equation} where the $D_{l_i}$ the local left covariant derivatives defined 
by (\ref{grob}).
\\Now we define a linear map $hor_l: \Gamma^1_c({\cal P}) \longrightarrow 
\oplus_{i \in I} \Gamma(B_i) \hat{\otimes} H$ by \[hor_l((\gamma_i)_{i\in 
I}):=(hor_{l_i}(\gamma_i))_{i \in I},~~(\gamma_i)_{i\in I} \in \Gamma^1_c({\cal 
P}.\] It remains to prove that the image of $hor_l$ lies in 
$\Gamma^1_c({\cal P})$. Then it follows immediately from the properties of the 
local connections $hor_{l_i}$ that $hor_l$ is a connection. 
\\To prove $hor_l(\Gamma^1_c({\cal P}) \subset \Gamma^1_c({\cal P})$ 
we need a lemma.
\begin{lem} \label{vker} $ hor_{l_i}((\chi_{i_{\Gamma}}(ker
\chi_{j_{\Gamma}}))^1) \subset
(\chi_{i_{\Gamma}}(ker \chi_{j_{\Gamma}}))^1 $ \end{lem}
Proof of the Lemma: Using the form of the generators of 
$J^i(B_{ij} \otimes H)$ (\ref{gener1})
-(\ref{gener5}) one finds easily that the differential calculus
$\Gamma^i(B_{ij} \otimes H)$ has the form $\Gamma^i(B_{ij} \otimes H)
=(\Gamma(B_{ij}) \hat{\otimes} \Gamma(H))/J$ where the differential ideal $J$
is generated by
\begin{eqnarray} \label{imchiij1}
&& \{\sum \tau_{ij}(r_2)d\tau_{ji}(r_1) \hat{\otimes} 1
+\sum \tau_{ij}(r_4)\tau_{ji}(r_1) \hat{\otimes} S^{-1}(r_3)dr_2| r \in R \},
\\ \label{imchiij2}
&& \{ \sum (a d\tau_{ji}(h) - (d\tau_{ji}(h))a) \hat{\otimes} 1| h
\in H,\,a \in B_i \}. \end{eqnarray} The identity 
$J= (\pi^i_{j_\Gamma} \otimes 
id)( \chi_{i_\Gamma} (ker \chi_{j_\Gamma})) $ is evident.
\\The factorization map $id^i_{{ij}_{\Gamma}}: \Gamma(B_{ij}) \hat{\otimes} 
\Gamma(H)) \longrightarrow \Gamma^i(B_{ij} \otimes H)$ fulfills 
\begin{equation} \label{idkreuz} id^i_{{ij}_{\Gamma}} \circ (\pi^i_{j_\Gamma} 
\otimes id)=(\pi^i_j \otimes id)_{\Gamma}.\end{equation}
Since $hor_{l_i}$ is a left modul homomorphism and $ker( \pi^i_{j_\Gamma} 
\otimes id)=ker \pi^i_{j_\Gamma} \hat{\otimes} \Gamma(H)$ one has 
\begin{equation} 
\label{nasoblub} hor_{l_i} ((ker 
(\pi^i_{j_\Gamma} \otimes id))^1) \subset (ker (\pi^i_{j_\Gamma} \otimes 
id))^1,\end{equation} thus 
$hor_{l_i}$ defines a connection $hor^{ij}_{l_i} : (\Gamma(B_{ij}) \hat{\otimes} \Gamma(H))^1 
\longrightarrow \Gamma^1(B_{ij}) \hat{\otimes} H$ by 
\begin{equation} \label{dasdach} hor^{ij}_{l_i} \circ (\pi^i_{j_\Gamma} \otimes id)=(\pi^i_{j_\Gamma} \otimes 
id) \circ hor_{l_i}. \end{equation}
Because of (\ref{krukru}), (\ref{nasoblub}) and \[(\pi^i_j \otimes id)_{\Gamma} 
\circ hor_{l_i} (\chi_{i_\Gamma} (ker
\chi_{j_\Gamma})) = id^i_{{ij}_{\Gamma}} \circ hor^{ij}_{l_i}(J)\] (which is 
immediate from (\ref{idkreuz}) and (\ref{dasdach})), to prove the assertion 
of the lemma we have to show that 
$id^i_{{ij}_{\Gamma}} \circ hor^{ij}_{l_i}(J)=0.$  
Note that the part of $J$ generated by (\ref{imchiij2}) lies in the horizontal
submodule $\Gamma^1(B_{ij}) \hat{\otimes} H$ and is therefore invariant under 
$hor^{ij}_{l_i}$. 
Now let us consider the part of $J$ generated by (\ref{imchiij1}).
Since $hor^{ij}_{l_i}$ is a left module homomorphism, it is sufficient to
consider the the product of the generators (\ref{imchiij1}) with 
a general element $(a \otimes h) \in B_{ij} \otimes H$ on the right.
Using $R \subset ker \varepsilon$, such an element can be written 
\begin{eqnarray*}&& \sum \tau_{ij}(r_2)d\tau_{ji}(r_1)a 
\hat{\otimes} h + \sum \tau_{ij}(r_4)\tau_{ji}(r_1)
\hat{\otimes} S^{-1}(r_3)dr_2)(a \otimes h)\\
&=& \sum \tau_{ij}(r_2)d\tau_{ji}(r_1)a
\hat{\otimes} h + \sum \tau_{ij}(r_4)\tau_{ji}(r_1)a
\hat{\otimes} S^{-1}(r_3)d(r_2h),~~r \in R, h \in H , a \in 
B_{ij}. \end{eqnarray*} 
\\Using (\ref{nasoeinhor1}), $R \subset ker \varepsilon$, (\ref{abl2}), 
$R \subset ker A_{l_j}$ and (\ref{dadumm})
one calculates \begin{eqnarray*} 
&& id^i_{{ij}_{\Gamma}}  \circ	hor^{ij}_{l_i}
 (\sum \tau_{ij}(r_2)
(d\tau_{ji}(r_1))a \hat{\otimes} h \\ &&+
\sum \tau_{ij}(r_4)\tau_{ji}(r_1)a \hat{\otimes}
S^{-1}(r_3)d(r_2h) ) \\
&=& \sum \iota^i_{{ij}_{\Gamma_m}}(a \tau_{ij}(r_2)d\tau_{ji}(r_1))
(1 \otimes h) \\
&&- \sum \iota^i_{{ij}_{\Gamma_m}} (a \tau_{ij}(r_3)\tau_{ji}(r_1) 
\pi^i_{j_{\Gamma_m}}(A_{l_i}(r_2h_1)))(1 \otimes h_2) \\
&=& \sum \iota^i_{{ij}_{\Gamma_m}} (a\tau_{ij}(r_2)d\tau_{ji}(r_1))
(1 \otimes h) \\
&&- \sum \iota^i_{{ij}_{\Gamma_m}}
(a \tau_{ij}(r_5)\tau_{ji}(r_1) \tau_{ij}(r_2h_1) 
\tau_{ji}(r_4h_3)\pi^j_{i_{\Gamma_m}}(A_{l_j}(r_3h_2)))(1 \otimes h_4) \\
&&- \sum \iota^i_{{ij}_{\Gamma_m}}
(a \tau_{ij}(r_4)\tau_{ji}(r_1) 
\tau_{ij}(r_2h_1)d\tau_{ji}(r_3h_2))(1 \otimes h_3)\\
&=& \sum \iota^i_{{ij}_{\Gamma_m}}
(a \tau_{ij}(r_2)d\tau_{ji}(r_1))(1 \otimes h) \\
&&- \sum \iota^i_{{ij}_{\Gamma_m}}
(a \tau_{ij}(h_1)\tau_{ji}(h_3)
\pi^j_{i_{\Gamma_m}}(A_{l_j}(rh_2)))(1 \otimes h_4) \\
&&- \sum \iota^i_{{ij}_{\Gamma_m}}
(a \tau_{ij}(r_2)d\tau_{ji}(r_1))(1 \otimes h)=0. \end{eqnarray*}
The last identity comes from the fact that $R$ is a right ideal. \hfill$\Boxneu$
\\\\Let $(\gamma_i)_{i \in I} \in \Gamma^1_c({\cal P})$.  
We have to prove that 
\[ (\pi^i_j \otimes id)_{\Gamma} \circ hor_{l_i} (\gamma_i) =
\phi_{{ij}_{\Gamma}} \circ (\pi^j_i \otimes id)_{\Gamma} \circ 
hor_{l_j} (\gamma_j).\] 
$\gamma_i$ has the general form \[\gamma_i=\sum_k
\chi_i(f_k^0)d\chi_i(f_k^1); 
\,f_k^0,f_k^1 \in {\cal P}.\] Using the compability condition of 
(\ref{glue}) and
(\ref{krukru}) one verifies that $\gamma_j$ has the form \[\gamma_j=\sum_k 
\chi_j(f^0_k)d\chi_j(f^1_k) + \rho, \;\; \rho \in 
\chi_{j_{\Gamma}}(ker\chi_{i_{\Gamma}}). \]
Now one obtains from  Lemma \ref{vker}, (\ref{smar}) and (\ref{hord})  
\begin{eqnarray*}
(\pi^i_j \otimes id)_{\Gamma}\circ hor_{l_i}(\gamma_i)
&=& (\pi^i_j \otimes id)_{\Gamma}\circ hor_{l_i}(\sum_k \chi_i(f^0_k)d
\chi_i(f^1_k))) \\ &=&
\sum_k (\pi^i_j \otimes id)_{\Gamma}(\chi_i(f^0_k) D_{l_i}(\chi_i(f^1_k))) \\
&=& \phi_{{ij}_{\Gamma}} \circ (\pi^j_i \otimes id)_{\Gamma} \circ
(\sum_k \chi_j(f^0_k)D_{l_j}(\chi_j(f^1_k)) \\ &=& \sum_k \phi_{{ij}_{\Gamma}}
\circ (\pi^j_i \otimes id)_{\Gamma} \circ hor_{l_j}(\chi_j(f^0_k)d
\chi_j(f^1_k)+\rho) \\ &=& \sum_k \phi_{{ij}_{\Gamma}} (\pi^j_i \otimes
id)_{\Gamma} \circ hor_{l_j} (\gamma_j),
\end{eqnarray*}
and the assertion is proved. \hfill$\Boxneu$

\begin{pr} 
There exists a bijection between left and right connections. 
\end{pr}
Proof: A left connection corresponds to a family of linear
maps $(A_{l_i})_{i \in I}$ satisfying ({\ref{abl1}), (\ref{abl2}) and 
(\ref{kerrl}). The linear maps $A_{r_i}:=-A_{l_i} \circ S$ fulfill (\ref{abl1}) 
and (\ref{kerrr}), thus the $A_{r_i}$ define right connections on the 
trivilizations. One has to prove that the
family $(A_{r_i})_{i \in I}$ satisfies (\ref{abl2}). Using
$\tau_{ij} \circ S=\tau_{ji}$ and $\sum d(\tau_{ij}(h_1)\tau_{ji}(h_2))=0$,
one calculates
\begin{eqnarray*} \pi^i_{j_{\Gamma_m}}(A_{r_i}(h)&=&-
\pi^i_{j_{\Gamma_m}}(A_{l_i}(S(h))) \\ &=& -\sum \tau_{ij}(S(h_3))
\pi^j_{i_{\Gamma_m}}(A_{l_j}(S(h_2)))\tau_{ji}(S(h_1))-\sum
\tau_{ij}(S(h_2))d\tau_{ji}(S(h_1)) \\
&=& -\sum \tau_{ji}(h_3)\pi^j_{i_{\Gamma_m}}(A_{l_j}(S(h_2)))\tau_{ij}(h_1)
- \sum \tau_{ji}(h_2)d\tau_{ij}(h_1) \\
&=& -\sum \tau_{ij}(h_1)
\pi^j_{i_{\Gamma_m}}(A_{l_j}(S(h_2)))\tau_{ji}(h_3) +
\sum \tau_{ij}(h_1)d \tau_{ji}(h_2) \\ &=& \sum
\tau_{ij}(h_1)\pi^j_{i_{\Gamma_m}}(A_{r_j}(h_2))\tau_{ji}(h_3) +
\sum \tau_{ij}(h_1)d \tau_{ji}(h_2).
\end{eqnarray*}
{}\hfill$\Boxneu$
\\\\Remark: A left (right) connection $hor_{l,r}$ and the corresponding 
left {right} covariant derivatives $D_{l,r}$ are connected by 
$hor_{l,r} \circ d = D_{l,r}|_{\cal P}$. Note that $hor_{l,r}$ can be extended 
to the submodule \[\Pi({\cal P}):=\{ \gamma \in \Gamma_c({\cal 
P})|~\chi_{i_{\Gamma_c}}(\gamma) \in (\Gamma(B_i) \hat{\otimes} H) \oplus 
(\Gamma(B_i) \hat{\otimes} \Gamma^1(H))\}.\] 
This means that the equation 
is valid on \[hor_{l,r} \circ d=D_{l,r}\] (equation on $hor \Gamma_c({\cal 
P})$).	      
\\\\To discuss curvatures of covariant derivatives and connections we introduce 
the notion of left (right) pre-connection forms.
\begin{de} A left (right) pre-connection form $\omega_{l,r}$ ist a linear map 
$\omega_{l,r} : H \longrightarrow \Gamma^1_c({\cal P})$ satisfying \begin{eqnarray}
\label{kov11} \omega_{l,r}(1) &=& 0, \\
\label{kov33} \Delta_{\cal P}^{\Gamma} (\omega_l(h)) &=& \sum \omega_l(h_2) 
\otimes S(h_1)h_3, \\ 
\Delta_{\cal P}^{\Gamma} (\omega_r(h)) &=& \sum \omega_r(h_2) \otimes h_3 
S^{-1}(h_1), \\ \label{kov44}
(1-hor_c) \circ \chi_{i_{\Gamma_c}}(\omega_l(h))&=& - \sum 1 \hat{\otimes} 
S(h_1)dh_2,~~\forall i \in I, \\
(1-hor_c) \circ \chi_{i_{\Gamma_c}}(\omega_r(h))&=& - \sum 1 \hat{\otimes}
(dh_2) S^{-1}(h_1) ~~\forall i \in I. \end{eqnarray} \end{de}
\begin{pr} \label{formbuform} Left (right) covariant derivatives are in 
bijective correspondence to 
left (right) pre-connection forms. \end{pr}
Proof: Let $\omega_l$ be a left pre-connection form.
$\omega_l$ determines a family of linear 
maps $A_{l_i}$ by \begin{equation} \label{immergleich}
A_{l_i}(h):=-(id \otimes \varepsilon) \circ hor_c \circ 
\chi_{i_{\Gamma_c}}(\omega_l(h)). \end{equation}
Because of (\ref{kov11}) the $A_{l_i}$ fulfill (\ref{abl1}).
\\Using 
\[ (1-hor_c) \circ \chi_{i_{\Gamma_c}}(\omega_l(h)) + hor_c \circ 
\chi_{i_{\Gamma_c}}(\omega_l(h))=\chi_{i_{\Gamma_c}}(\omega_l(h)),\] 
(\ref{kov44}), (\ref{kov33}), Lemma \ref{id} and (\ref{immergleich}) one 
verifies easily \begin{equation} \label{schwigleich} 
\chi_{i_{\Gamma_c}}(\omega_l(h))=-\sum 1 \hat{\otimes} S(h_1)dh_2 - 
\sum A_{l_i}(h_2) \hat{\otimes} S(h_1)h_3. \end{equation} 
Since \[ (\pi^i_j \otimes id) \circ \chi_{i_{\Gamma_c}} 
(\omega_l(h))=\phi_{{ij}_\Gamma} \circ	(\pi^j_i \otimes id) \circ 
\chi_{j_{\Gamma_c}}(\omega_l(h)), \] an easy calculation (using (\ref{phi}) and 
the projection $P_{inv}$ (\ref{pinv})) leads to (\ref{abl2}).  
\\We want to prove that the left covariant derivative $D_l$ determined by the 
$A_{l_i}$ is  
\begin{equation} \label{dassollstimmen} D_l(\gamma) = d \gamma + (-1)^{n}\sum \gamma_0 \omega_l(\gamma_1),~~\gamma \in 
hor \Gamma^n_c({\cal P}).\end{equation}
It is sufficiant to prove that for $\gamma \in hor \Gamma^n_c({\cal P})$ 
\[\chi_{i_{\Gamma_c}}(d \gamma +(-1)^n \sum \gamma_0 \omega_l(\gamma_1))= 
d\gamma \hat{\otimes} h +(-1)^{n+1} \gamma
A_{l_i}(h_1) \hat{\otimes} h_2=\chi_{i_{\Gamma_c}} \circ D_l (\gamma).\]
$\chi_{i_{\Gamma_c}}(\gamma)$ has 
the general form \[\chi_{i_{\Gamma_c}}(\gamma)= \sum_k \gamma^k_i \hat{\otimes}
h^k_i;~~\gamma^k_i \in \Gamma^n(B_i),~h^k_i \in H.\]
Using (\ref{schwigleich}) one obtains \begin{eqnarray*}
&&\chi_{i_{\Gamma_c}}(d \gamma +(-1)^n \sum \gamma_0 \omega_l(\gamma_1))\\
&=& (\sum_k d \gamma^k_i \hat{\otimes} h^k_i + (-1)^n \sum_k 
\gamma^k_i \hat{\otimes} dh^k_i + (-1)^n \sum_k \sum (\gamma^k_i \otimes 
h^k_{i1}) \chi_{i_{\Gamma_c}}(\omega_l(h^k_{i2})) \\
&=& \sum_k d\gamma^k_i \hat{\otimes} h^k_i + (-1)^n \sum_k \gamma^k_i 
\hat{\otimes} dh^k_i -(-1)^n 
\sum_k \sum (\gamma^k_i \otimes
h^k_{i1})(1 \hat{\otimes} S(h^k_{i2})dh^k_{i3}) \\
&&- (-1)^n \sum_k \sum (\gamma^k_i \otimes
h^k_{i1})(A_{l_i}(h^k_{i3}) \hat{\otimes} S(h^k_{i2})dh^k_{i4}) \\
&=& d\gamma \hat{\otimes} h +(-1)^{n+1} \gamma
A_{l_i}(h_1) \hat{\otimes} h_2. \end{eqnarray*} 
\\\\Note the identity \begin{equation} \label{dirksoso} D_{l_i}
(\gamma \hat{\otimes} h) = d\gamma \hat{\otimes} h +(-1)^{n+1} \gamma
A_{l_i}(h_1) \hat{\otimes} h_2 = d(\gamma \hat{\otimes} h) + (-1)^n \sum
(\gamma \hat{\otimes} h_1) \chi_{i_{\Gamma_c}}(\omega_l(h_2)). \end{equation}
Assume now there is given a left covariant derivative $D_l$. In terms of the 
corresponding linear maps $A_{l_i}$ one obtains a family of left pre-connection 
forms $\omega_{l_i}: H \longrightarrow (\Gamma(B_i) \hat{\otimes} \Gamma(H))^1$
by \[ \omega_{l_i}(h)=-\sum 1 \hat{\otimes} S(h_1)dh_2 - \sum A_{l_i}(h_2) 
\hat{\otimes} S(h_1)h_3.\] 
Using (\ref{abl2}) one obtains 
\[ (\pi^i_j \otimes id)_{\Gamma} (\omega_{l_i}(h)) = \phi_{{ij}_{\Gamma}} \circ 
(\pi^j_i \otimes id)_{\Gamma} (\omega_{l_j}(h), \] thus one has by \[
\omega_l(h)=(\omega_{l_i}(h))_{i \in I} \]  a left pre-connection 
form $\omega_l: H \longrightarrow \Gamma^1_c({\cal P})$. 
\\One easily verifies for $\gamma \hat{\otimes} h \in
\Gamma^n(B_i) \hat{\otimes} H$
\begin{equation} \label{dirk} D_{l_i}(\gamma \hat{\otimes} h) =
d\gamma \hat{\otimes} h +(-1)^{n+1} \gamma
A_{l_i}(h_1) \hat{\otimes} h_2 = d(\gamma \hat{\otimes} h) + (-1)^n \sum
(\gamma \hat{\otimes} h_1) \omega_{l_i}(h_2). \end{equation}
Using this formula it follows that 
\begin{equation} D_l(\gamma)=d\gamma + 
(-1)^n \sum \gamma_0 \omega_l(\gamma_1) \end{equation} for $\gamma \in hor 
\Gamma^n_c({\cal P})$. It is immediate from the formulas (\ref{dirksoso}) and 
(\ref{dirk}) and Proposition \ref{kovbij} that the correspondence is bijective.
\\For right covariant derivatives the proof is analogous.
\hfill$\Boxneu$
\\\\Remark: Note that the foregoing proof also shows the bijectiv 
correspondence between left (right) pre-connection forms and families of linear 
maps $A_{{l,r}_i} : H \longrightarrow \Gamma(B_i)$ fulfilling (\ref{abl1}) and 
(\ref{abl2}).
\begin{de} A left (right) pre-connection form $\omega_{l,r}$ is called 
left (right) connection form, if \begin{eqnarray} \label{rdrin}
R &\subset& ker((id \otimes \varepsilon) \circ hor_c \circ 
\chi_{i_{\Gamma_c}} \circ \omega_l),~~\forall i \in I, \\ \label{rmdrin}
S^{-1}(R) &\subset& ker((id \otimes \varepsilon) \circ hor_c \circ
\chi_{i_{\Gamma_c}} \circ \omega_r),~~\forall i \in I 
\end{eqnarray} is satisfied.\end{de}
\begin{pr} Left (right) connections are in bijective correspondence to 
left (right) connection forms. 
\end{pr}
Proof: The claim follows immedately from Proposition \ref{vi}, Proposition
\ref{formbuform}  and (\ref{immergleich}). \hfill$\Boxneu$
\\\\Remark: Note that classical connection forms are related to the connection 
forms considered above as follows: Let a classical principal bundle be given, 
let $X$ be a vector field on the total space $Q$, and let $h \in C^{\infty}(G)$ 
where $G$ is the structure group. A classical connection form is a Lie algebra 
valued 1-form $\tilde{\omega}$ of type Ad on $Q$. Then the formula 
\[\omega_l(h)(X)=-\tilde{\omega}(X)(h) \] defines a left connection form 
$\omega_l$ in the above sense. Condition (\ref{rdrin}) with $R=(ker 
\varepsilon)^2$ means that $\tilde{\omega})$ can be interpreted as a Lie algebra 
valued form. In this case (\ref{kov33}) and (\ref{kov44}) replace the usual 
conditions (type $Ad$, condition for fundamental vectors) for connection forms.
\begin{de} 
The left (right) curvature of a given left (right) covariant derivative is
the linear map $D^2_{l,r}: hor \Gamma_c({\cal P}) \longrightarrow 
hor \Gamma_c({\cal P})$.
\end{de}
\begin{de} Let $\omega_{l,r}$ be a left(right) pre-connection form of a left (right)
covariant derivative $D_{l,r}$. The linear maps $\Omega_{l,r} :H 
\longrightarrow \Gamma^2_c({\cal P})$ 
defined by \begin{eqnarray} \label{curv} \Omega_l(h)&:=& d\omega_l(h)-\sum
\omega_l(h_1) \omega_l(h_2), \\ \Omega_r(h) &:=& d\omega_r(h) + \sum 
\omega_r(h_2) \omega_r(h_1) \end{eqnarray} are called the left (right) 
curvature form of a given left (right) covariant derivative. \end{de}
Remark: In other words we take an analogue of the structure equation as 
definition of the curvature form.
\begin{pr} The left (right) curvature of a given left (right) covariant 
derivative is related to the left (right) curvature form by the identity 
\begin{eqnarray} \label{curv1}D_l^2 (\gamma) &=& \sum \gamma_0 
\Omega_l(\gamma_1),\;\; \gamma \in hor \Gamma_c({\cal P}),\\  \label{curv2}
 D_r^2( \gamma) &=& \sum \Omega_r(\gamma_1) \gamma_0,\;\;
\gamma \in hor \Gamma_c({\cal P}).\end{eqnarray}\end{pr} 
Proof: Because of the one-to-one correspondence between covariant 
derivatives on 
$\cal P$ and certain families of covariant derivatives on the trivilizations 
$B_i \otimes H$ it is sufficient to prove this assertion on a trivial bundle 
$B \otimes H$.
In this case the linear map $\omega_l$ belonging to a left covariant derivative 
has the form \[ \omega_l(h)=-(1 \otimes \sum S(h_1)dh_2)) -\sum (A_l(h_2)
\hat{\otimes} S(h_1)h_3). \] Therefore, one obtains for $\Omega_l$ 
\begin{eqnarray*}
d \omega_l(h) - \sum \omega_l(h_1)\omega_l(h_2) &=& -1 \hat{\otimes} \sum
dS(h_1)dh_2 - \sum dA_l(h_2) \hat{\otimes} S(h_1)h_3 \\ &&+ \sum A_l(h_2)
\hat{\otimes} 
(dS(h_1))h_3 + \sum A_l(h_2) \hat{\otimes} S(h_1)dh_3 \\ &&- 1
\hat{\otimes} 
\sum S(h_1)(dh_2)S(h_3)dh_4 -\sum A_l(h_2) \hat{\otimes} S(h_1)h_3 S(h_4)
dh_5 \\ &&+ \sum A_l(h_4) \hat{\otimes} S(h_1)(dh_2)S(h_3)h_5 -\sum
A_l(h_2)A_l(h_5) \hat{\otimes} S(h_1)h_3S(h_4)h_6 \\
&=& -\sum dA_l(h_2) \hat{\otimes} S(h_1)h_3 -\sum A_l(h_2)A_l(h_3)
\hat{\otimes} S(h_1)h_4, \end{eqnarray*} 
which leads for $\gamma \in \Gamma^n(B)$ to \[
\sum (\gamma \hat{\otimes} h_1) \Omega_l(h_2) = -\sum \gamma dA_l(h_1)  
\hat{\otimes} h_2- \sum \gamma A_l(h_1) A_l(h_2) \hat{\otimes} h_3. \]
On the other hand the left hand side of (\ref{curv1}) is 
\begin{eqnarray*} (D_l)^2(\gamma \otimes h) &=& 
D_l(d \gamma \hat{\otimes} h - \sum (-1)^n\gamma A_l(h_1)\hat{\otimes} h_2)\\
&=& -(-1)^{n+1} \sum (d \gamma) A_l(h_1)\hat{\otimes} h_2 -(-1)^n \sum	
(d \gamma) A_l(h_1) \hat{\otimes} h_2 \\ &&- \sum \gamma d A_l(h_1)
\hat{\otimes} h_2 
-\sum \gamma A_l(h_1)A_l(h_2) \hat{\otimes} h_3 \\
&=& -(\sum \gamma dA_l(h_1) \hat{\otimes} h_2+ \sum \gamma A_l(h_1) A_l(h_2) 
\hat{\otimes} h_3) \\
&=& \sum (\gamma \hat{\otimes} h_1) \Omega_l(h_2). \end{eqnarray*}
For right covariant derivatives the proof is analogous. \hfill$\Boxneu$
\\\\Remark: The proof shows that there is a linear map $F_{l,r}: H
\longrightarrow \Gamma^2(B)$ defined by \begin{eqnarray} 
\label{fl} F_l(h)&:=& dA_l(h) + \sum A_l(h_1) A_l(h_2), \\
\label{fr} F_r(h)&:=& dA_r(h) - \sum A_r(h_2) A_r(h_1) \end{eqnarray}
such that the left (right) curvature form of a given left (right) 
covariant derivative on a trivial QPFB has the form 
\begin{eqnarray} \Omega_l(h) &=& -\sum F_l(h_2) \hat{\otimes} S(h_1)h_3,\\
\Omega_r(h) &=& -\sum F_r(h_2) \hat{\otimes} h_3S^{-1}(h_1).\end{eqnarray} 
Using formula (\ref{abl2}) and the Leibniz rule (taking into account 
$\sum\tau_{ij}(h_1) \tau_{ij}(h_2) = \varepsilon(h)1$) it is easy to verify 
that  
the family of linear maps $F_{{l,r}_i}: H \longrightarrow \Gamma^2(B_i)$ 
corresponding to a left (right) curvature form on a locally trivial QPFB  
satisfies
\begin{equation} \pi^i_{j_{\Gamma_m}}(F_{{l,r}_i}(h))=\sum \tau_{ij}(h_1) 
\pi^j_{i_{\Gamma_m}}(F_{{l,r}_j}(h_2)) \tau_{ji}(h_3). \end{equation}
In general, an analogue of the Bianchi identity does not exist.
\\Now we make some remarks about the general form of the linear 
maps $A_{l_i}: H \longrightarrow \Gamma^1(B_i)$ corresponding to
connections on a locally trivial QPFB. For this we use the functionals
${\cal X}_i$ corresponding to the right ideal $R$, which  determines the right 
covariant 
differential calculus $\Gamma(H)$ (see \cite{wor1}, \cite{wor2} and the
appendix). 
Let the $h^i +R \in ker \varepsilon/R $ be a linear basis in $ker 
\varepsilon/R$. Then every element $h- \varepsilon(h)1 +R \in  ker
\varepsilon/R$ 
has the form ${\cal X}_i(h)h^i +R$. Since $1 \in ker A_{l_i}$ and $R \subset
ker A_{l_i}$ it follows that $A_{l_i}$ is determined by its values on the $h^k$,
\[A_{l_i}(h)=\sum_k {\cal X}_k(h) A_{l_i}(h^k). \]
In other words, to get a connection on the trivial pieces $B_i \otimes H$,
one chooses $A^k_i \in \Gamma^1(B_i)$ and defines the linear map $A_{l_i}$ by
\[ A_{l_i}(h)=\sum_k {\cal X}_k(h) A^k_i.\] 
The  connections so defined on the trivial pieces $B_i \otimes H$ do 
in general not give a connection on 
the locally trivial QPFB ${\cal P}$, because they do in general 
not fulfill the
condition (\ref{abl2}). If the right ideal $R$ fulfills (\ref{jbij}),
one can rewrite the condition  
(\ref{abl2}) as a condition for the one forms $A^k_i \in\Gamma^1(B_i)$. 
Recall that in this case 
$\sum \tau_{ij}(r_1)d\tau_{ji}(r_2)=0,~\forall r \in R$ (cf. (\ref{fish2})), 
thus \[ \sum 
\tau_{ij}(h_1)d\tau_{ij}(h_2)= \sum \sum_k {\cal X}_k(h)
\tau_{ij}(h^k_1)d\tau_{ji}(h^k_2).\] 
Furthermore, the 
condition (\ref{jbij}) leads to the identity 
\[ \sum \tau_{ij}(h_1) {\cal X}_l(h_2) \tau_{ji}(h_3)=
\sum \tau_{ji}(S(h_1)h_3) {\cal X}_l(h_2)=\sum_k {\cal X}_k(h) \sum  
\tau_{ji}(S(h^k_1)h^k_3){\cal X}_l(h^k_2). \]
Putting now $A_{l_i}(h)=\sum_k {\cal X}_k(h)A^k_i$ in (\ref{abl2}) leads to 
the following condition for the forms $A^k_i$:
\[ \pi^i_{j_{\Gamma_m}}(A^k_i)=\sum_l \tau_{ji}(S(h^k_1)h^k_3){\cal 
X}_l(h^k_2) \pi^j_{i_{\Gamma_m}}(A^l_i) + \tau_{ij}(h^k_1)d 
\tau_{ji}(h^k_2). \]
Note that, in the case $I=(1,2)$, it follows from the last formula that there 
exist connections. One can choose, e.g., one forms $A^l_2$ on the right, and 
solve the remaining equation for $A^k_1$ due to the surjectivity of 
$\pi^1_{2_{\Gamma_g}}$.
One can regard the set of all left (right) covariant derivatives ${\cal 
D}_{l,r}$ as a set with affine structure, where the corresponding vector space
is characterized by \begin{pr} \label{shmap} A linear map $C_{l,r}: 
hor\Gamma_c({\cal P}) \longrightarrow hor \Gamma_c({\cal P})$
is a difference of two left (right) covariant derivatives if and only if:
\begin{eqnarray}
\label{sh1} && C_{l,r}(1)=0, \\
&&C_{l,r}(hor \Gamma^n_c({\cal P})) \subset hor \Gamma^{n+1}_c({\cal P}),\\
\label{sh2} && C_{l}(\gamma \alpha) = (-1)^n \gamma C_l(\alpha);~~\gamma \in
\Gamma^n_c(B);~\alpha \in hor \Gamma_c({\cal P}),\\
\label{sh3} && C_r(\alpha \gamma) = (-1)^n C_r(\alpha) \gamma;~~\gamma \in
\Gamma_c(B);~\alpha \in hor \Gamma^n_c({\cal P}),\\
\label{ksh4} &&(C_{l,r} \otimes id) \circ \Delta_{{\cal P}_{\Gamma_c}} =
\Delta_{{\cal P}_{\Gamma_c}} \circ C_{l,r}, \\
\label{sh5} &&C_{l,r}(ker \chi_{i_{\Gamma_c}}|_{hor \Gamma_c({\cal P})})
\subset ker \chi_{i_{\Gamma_c}}|_{hor\Gamma_c({\cal P})};~~\forall i \in I .
\end{eqnarray}	\end{pr}
This is immediate from Definition \ref{kovohnz}.
\\Because of (\ref{sh5}) such a map $C_{l,r}$ defines a 
family of local maps $C_{{l,r}_i}$ by 
\[ C_{{l,r}_i} \circ \chi_{i_{\Gamma_c}} = \chi_{i_{\Gamma_c}} \circ C_{l,r}.\]
It is immediately that the set of left (right) connections is an affine 
subspace of ${\cal D}_{l,r}$. The elements of the corresponding vector space 
have the following additional property:
\begin{eqnarray*} 
(id \otimes \varepsilon) \circ C_{l_i}(1 \otimes 
r)&=&0,~~\forall i \in I,~~\forall r \in R, \\
(id \otimes \varepsilon) \circ C_{r_i}(1 \otimes
r)&=&0,~~\forall i \in I,~~\forall r \in S^{-1}(R). \end{eqnarray*}

\section{Example \label{ex}}

Here we present an example of a $U(1)$-bundle over the quantum space 
$S^2_{p q \phi}$. The quantum space $S^2_{p q \phi}$ is treated in detail
in \cite{cama} and we restrict ourselves here to a brief summary.
\\\\
The algebra $P(S^2_{p q \phi})$ of all polynomials over the quantum space 
$S^2_{p q \phi}$ is constructed by gluing together two copies of a quantum
disc along its classical subspace. 
\begin{de} The algebra $P(D_p)$ of all polynomials over the quantum disc
$D_p$ is defined as the algebra generated by the elements $x$ and $x^*$
fulfilling the relation 
\begin{equation} x^*x-pxx^*=(1-p)1, \end{equation} where $0<p<1$. \end{de}
Let $P(S^1)$ be the algebra generated by the elements $\alpha,\alpha^*$
fulfilling the 
relation \[ \alpha \alpha^*=\alpha^* \alpha=1.\] $P(S^1)$ can be considered as 
the algebra of all trigonometrical polynomials over the circle $S^1$.
\\There exists a surjective homomorphism $\phi_p: P(D_p) \longrightarrow
P(S^1)$ defined by \begin{eqnarray*}
\phi_p(x) &=& \alpha, \\ \phi_p(x^*) &=& \alpha^*, \end{eqnarray*} and one
can consider this homomorphism as the ``pull back'' of the embedding of 
the circle into the quantum disc. 
\\The algebra $P(S^2_{p q \phi})$ of all polynomials over the quantum
space $S^2_{p q \phi}$ is defined as \begin{equation}
P(S^2_{p q \phi}):=\{ (f,g) \in P(D_p) \bigoplus
P(D_q)|\phi_p(f)=\phi_q(g) \}. 
\end{equation}
 It was shown in \cite{cama} that one can also regard this algebra as the 
 algebra generated by the 
 elements $f_1$, $f_{-1}$ and $f_0$ fulfilling the relations
 \begin{eqnarray} f_{-1}f_1- q f_1 f_{-1} &=& (p-q)f_0 + (1-p)1, \\
 f_0 f_1- p f_1 f_0 &=& (1-p) f_1, \\
 f_{-1} f_0 -p f_0 f_{-1} &=& (1-p) f_{-1}, \\
 (1-f_0)(f_1f_{-1} - f_0) &=& 0, 
 \end{eqnarray}
where the isomorphism is given by $f_1 \rightarrow (x,y)$, $f_{-1} \rightarrow
(x^*,y^*)$ and $f_0 \rightarrow (xx^*,1)$. (Here, the generators of $P(D_q)$ are
denoted by $y$ and $y^*$.) 
It is proved in \cite{cama} that the $C^*$-closure $C(S^2_{p q \phi})$ of 
$P(S^2_{p q \phi})$ is isomorphic to the $C^*$-algebra $C(S^2_{\mu c})$
over the Podles sphere $S^2_{\mu c}$ for $c>0$.
\\\\Now, let us construct a class of QPFB's with structure group $U(1)$ and
base space $S^2_{p q \phi}$. The algebra of polynomials $P(U(1))$ over 
$U(1)$ is the same algebra as $P(S^1)$.
With $\Delta(a)=\alpha \otimes \alpha$, $\varepsilon(\alpha)=1$ and 
$S(\alpha)=\alpha^*$, $P(U(1))$ is a Hopf algebra. According to Proposition 
\ref{transf} we need just one transition 
function $\tau_{12}: P(U(1)) \longrightarrow P(S^1)$ to obtain a
locally trivial QPFB. We define a class of transition functions 
$\tau^{(n)}_{12}$ as follows: \begin{eqnarray*}
\tau^{(n)}_{12}(\alpha) &:=& \alpha^n, \\
\tau^{(n)}_{12}({\alpha^*}) &:=& {\alpha^*}^n. \end{eqnarray*}
It follows that \begin{eqnarray*} \tau^{(n)}_{21}(\alpha)&:=&{\alpha^*}^n, \\
\tau^{(n)}_{21}({\alpha^*})&:=&\alpha^n. \end{eqnarray*}
We obtain a class of locally trivial QPFB's 
$({\cal P}^{(n)}, \Delta_{{\cal P}^{(n)}},P(U(1)), P(S^2_{p q \phi}), \iota,
((\chi_p,ker \pi_p), (\chi_q, ker \pi_q)) )$ 
corresponding to these transition functions (see formulas (\ref{phi}) and 
(\ref{P})), where $\iota$ is the canonical embedding $P(S^2_{p q \phi})
\subset {\cal P}^{(n)}$ and $\pi_{p,q}: P(S^2_{p q \phi}) \longrightarrow 
P(D_{p,q})$ and 
$\chi_{p,q}: {\cal P}^{(n)} \longrightarrow P(D_{p,q}) \otimes P(U(1))$ are the
restrictions of the canonical projections on $P(S^2_{p q \phi})$ and ${\cal
P}^{(n)}$ respectively. 
\\In the following, we restrict ourselves 
to the case $n=1$. 
\begin{pr} Let $J \subset P(D_p) \otimes P(D_q)$ be the ideal
generated by the element \[(xx^*-1) \otimes (yy^*-1).\] Then
${\cal P}^{(1)}$ is algebra isomorphic to $(P(D_p) \otimes P(D_q))/J$. \end{pr}
Proof: $(P(D_p) \otimes P(D_q))/J$  by \begin{eqnarray*} 
a &=& 1 \otimes y + J, \\
a^* &=& 1 \otimes y^* +J, \\ b &=& x \otimes 1 +J \\
b^* &=& x^* \otimes 1 +J. \end{eqnarray*}
It is easy to see that $(P(D_p) \otimes P(D_q))/J$ can be considered as the 
algebra $\C<a,a^*,b, b^*>/\tilde{J}$, where the ideal $\tilde{J}$ is generated 
by the relations
\begin{eqnarray}
{a}^*{a} - q{a} {a}^* &=& (1-q) 1, \nonumber\\
{b}^* {b}-p{b} {b}^* &=& (1-p) 1, \nonumber \\
{b} {a}={a} {b},~~ {b} {a}^*={a}^*{b}^*,~~{b}^*{a}={a}
{b}^*~~,{b}^*{a}^*&=& {a}^*{b}^* \label{prel}, \\
(1-{a}{a}^*)(1-{b} {b}^*)&=& 0. \nonumber \end{eqnarray}
Further consider the following elements in ${\cal P}^{(1)}$:
\begin{eqnarray*} \tilde{a} &=& (1 \otimes \alpha , y \otimes \alpha), \\
\tilde{a}^* &=& (1 \otimes \alpha^*, y^* \otimes \alpha^*), \\
\tilde{b} &=& (x \otimes \alpha^*, 1 \otimes \alpha^*), \\
\tilde{b}^* &=& (x^* \otimes \alpha, 1 \otimes \alpha). \end{eqnarray*}
A short calculation shows that these elements fulfill the same relations 
(\ref{prel}) as the $a$, $a^*$, $b$ and $b^*$.
Thus, there exists a homomorphism $F: (P(D_p) \otimes P(D_q))/I
\longrightarrow {\cal P}^{(1)}$ defined by \[ F(a):=\tilde{a},~~F(b):= \tilde{b},~~
F(a^*):=\tilde{a}^*,~~F(b^*):=\tilde{b}^*.\]
We will show that $F$ is an isomorphism. For surjectivity it is sufficient to 
show that the elements $\tilde{a}$, $\tilde{a}^*$, $\tilde{b}$ and 
$\tilde{b}^*$ generate the algebra ${\cal P}^{(1)}$.
It is shown in \cite{cama} Lemma 2 that the elements $x^k{x^*}^l,~~k,l \geq 0$ form
a vector space basis of $P(D_p)$. Analogous the elements 
$\alpha^i, i \in \Z$ ($a^{-1}=a^*$), form a vector space basis in $P(U(1))$. 
Thus a general element $f \in P(D_p) \otimes P(U(1)) \bigoplus P(D_q)
\otimes  P(U(1))$ has the form 
\[ f=(\sum_{k,l \geq 0,i} c^p_{k,l,i} x^k {x^*}^l \otimes \alpha^i,
\sum_{m,n \geq 0,j \in \Z} c^q_{m,n,j} y^m{y^*}^n \otimes \alpha^j). \]
$f \in {\cal P}^{(1)}$ means that there is the restriction 
\[ \sum_{k,l \geq 0,i \in \Z} c^p_{k,l,i} \alpha^{k-l} \otimes \alpha^i =
\sum_{m,n \geq 0,j} c^q_{m,n,j} \alpha^{m-n-j} \otimes \alpha^j, \]
which leads to the following condition for the coefficients $c^p_{k,l,i}$ and
$c^q_{m,n,j}$.
\begin{equation} \label{soeinc} \sum_{l \geq 0, s+l \geq 0} c^p_{s+l,l,t}= \sum_{n \geq 0, n+s+t \geq 0}
c^q_{s+t+n,n,t},~~\forall s,t \in \Z. \end{equation}
$f \in {\cal P}^{(1)}$ has the form $f= \sum_{s,t} f_{s,t}$, where	
\[ f_{s,t}= (\sum_{l \geq 0, l+s \geq 0}  c^p_{s+l,l,t}  x^{l+s} {x^*}^l 
\otimes \alpha^t,\sum_{n \geq 0, n+s+t \geq 0} c^q_{s+t+n,n,t}
y^{n+s+t}{y^*}^n \otimes \alpha^t) \in {\cal P}^{(1)} \] due to (\ref{soeinc}).
Because of (\ref{soeinc}) one can write $f_{s,t}$ as 
\begin{eqnarray*} f_{s,t} &=& \sum_{l \geq 0, l+s \geq 0}  c^p_{s+l,l,t}
(x^{l+s} {x^*}^l 
\otimes \alpha^t, y^{m+s+t}{y^*}^m \otimes \alpha^t) 
\\ &&+ \sum_{n \geq 0,
n+s+t \geq 0} c^q_{s+t+n,n,t}
 (x^{k+s}{x^*}^k \otimes \alpha^t, y^{n+l+t}{y^*}^n
\otimes \alpha^t) \\
&&- \sum_{l \geq 0, l+s \geq 0}  c^p_{s+l,l,t} (x^{k+s}{x^*}^{k} 
\otimes \alpha^t, y^{m+s+t}{y^*}^m \otimes \alpha^t). \end{eqnarray*}
The identity  
\[(x^s x^l {x^*}^l \otimes \alpha^t, y^{s+t} y^n {y^*}^n \otimes
\alpha^t)= \tilde{a}^{s+t+n} \tilde{a}^*{}^n  \tilde{b}^{s+l} 
\tilde{b}^*{}^l, \] which is a direct consequence of the definition of 
$\tilde{a}$, $\tilde{a}^*$, $\tilde{b}$ and $\tilde{b}^*$, shows that $F$ is surjective.
\\To show the injectivity of $F$ we define the homomorphisms 
$F_{p,q} : (P(D_p) \otimes P(D_q))/I \longrightarrow 
P(D_{p,q}) \otimes P(U(1))$ by $F_{p,q}:= \chi_{p,q} \circ F$. 
Because of $ker \chi_p \cap ker \chi_q =\{0\}$, $ker F=\{0\}$ if and only if 
$ker F_p \cap ker F_q=\{0\}$. First let us describe the ideals $ker F_{p,q}$. 
Let $I_p$ and $I_q$ be the ideals generated by 
$1-aa^*$ and $1-bb^*$ respectively. From (\ref{prel}) 
it is immediate that the algebras 
$((P(D_p)\otimes P(D_q))/I)/I_{p,q}$ are isomorphic to 
$P(D_{p,q}) \otimes P(U(1))$, where the isomorphism 
$((P(D_p)\otimes P(D_q))/I)/I_p \longrightarrow P(D_p) \otimes P(U(1))$ is 
defined by $a \rightarrow 1 \otimes \alpha$, $b \rightarrow x \otimes 1$, and 
the isomorphism   
$((P(D_p)\otimes P(D_q))/I)/I_q \longrightarrow P(D_q) \otimes P(U(1))$ is 
defined by $a \rightarrow y \otimes 1$, $b \rightarrow 1 \otimes \alpha$. 
\\Moreover, there are automorphisms $\tilde{F}_{p,q}: P(D_{p,q}) \otimes P(U(1)) 
\longrightarrow P(D_{p,q}) \otimes P(U(1))$ defined by 
\begin{eqnarray*}
\tilde{F}_p(1 \otimes \alpha)&:=& 1 \otimes \alpha,~~\tilde{F}_q(1 \otimes 
\alpha):= 1 \otimes \alpha, \\
\tilde{F}_p(1 \otimes \alpha^*)&:=& 1 \otimes \alpha^*,~~\tilde{F}_q(1
\otimes \alpha^*):= 1 \otimes \alpha^*, \\
\tilde{F}_p(x \otimes 1)&:=& x \otimes \alpha^*,~~\tilde{F}_q(y \otimes 
1):= y \otimes \alpha, \\
\tilde{F}_p(x^* \otimes 1)&:=& x^* \otimes \alpha,~~\tilde{F}_q(y^* \otimes
1):= y^* \otimes \alpha^*. \end{eqnarray*}
Let $\eta_{p,q}$ be the quotient
maps with respect to the ideals $I_{p,q}$. A short calculation shows that 
\[ F_{p,q}=\tilde{F}_{p,q} \circ \eta_{p,q}, \] thus we have found 
$ker F_{p,q}=I_{p,q}$. It remains to show $I_p \cap I_q=\{0\}$.
There are the following identities in $P(D_q) \otimes P(D_p)/I$:
\begin{eqnarray*}
(1-aa^*)a &=& qa(1-aa^*),~~(1-aa^*)a^*=q^{-1}a^*(1-aa^*), \\
(1-bb^*)b &=& p b(1-bb^*),~~(1-bb^*)b^*=p^{-1} b^*(1-bb^*). 
\end{eqnarray*} From these relations and the definition of $I_p=ker F_p$
follows that for  $f \in ker F_p$ there exists an element 
$\tilde{f}$ such that $f=(1-aa^*)\tilde{f}$. 
$ker F_q$ has an anlogous property with $1-bb^*$, instead of $1-aa^*$.
Using that $1-xx^*$ is not a zero 
divisor in $P(D_p)$, see \cite{cama} Lemma 3, it is now easy to see that $f \in 
ker F_p \cap ker F_q$ is of the form $f=(1-aa^*)(1-bb^*)\tilde{f}$.
Thus $f=0$, i.e. $ker F_p \cap ker F_q =0$.\hfill$\Boxneu$.
\\\\The proof has shown that ${\cal P}^{(1)}\simeq \C < a,a^*,b,b^* > /J$, where
$J$ is the ideal generated by the relations (\ref{prel}). Under this 
identification, the mappings belonging to the bundle can be given explicitly. 
\begin{eqnarray*}
&&\Delta_{{\cal P}^{(1)}}(a) = a \otimes \alpha,~~\Delta_{{\cal P}^{(1)}}(a^*)
= a^*\otimes \alpha^*, \\
&&\Delta_{{\cal P}^{(1)}}(b) = b \otimes \alpha^*,~~\Delta_{{\cal
P}^{(1)}}(b^*) = b^*\otimes \alpha,\\[.1cm]
&&\chi_p(a)= 1 \otimes \alpha,~~\chi_q(a) = y \otimes \alpha, \\
&&\chi_p(a^*)= 1 \otimes \alpha^*,~~\chi_q(a^*) = y^* \otimes \alpha^*, \\
&&\chi_p(b)= x \otimes \alpha^*,~~\chi_q(b) = 1 \otimes \alpha^*, \\
&&\chi_p(b^*)= x^* \otimes \alpha,~~\chi_q(b^*) = 1 \otimes \alpha, \\[.1cm]
&&\iota(f_1) = ba,~\iota(f_{-1})= a^*b^*,~\iota(f_0)=bb^*.\end{eqnarray*}
In the classical limit $p,q\lra 1$ the algebra becomes commutative and only 
the relation $(1-aa^*)(1-bb^*)=0$ remains. It is easy to see that this  
relation,
together with the natural requirement $|a|\leq 1$, $|b|\leq 1$, describes a
subspace of ${\R}^4$ homeomorphic to $S^3$. The right $U(1)$-action
is a simultaneous rotation in $a$ and $b$, and the orbit through $b=0$ is
the fibre over the top $(0,0,0)$ of the base space 
(see the discussion in \cite{cama}).
\\\\To build a connection on this locally trivial QPFB, first we have to  
construct an adapted covariant differential structure on ${\cal P}^{(1)}$. By 
Definition \ref{strukstrik}, the adapted covariant differential structure is
defined by giving differential calculi $\Gamma(P(D_p))$ and
$\Gamma(P(D_q))$ and a right covariant 
differential calculus $\Gamma(P(U(1))$ on the Hopf algebra $P(U(1))$. 
\\As the differential calculi $\Gamma((P(D_{p,q}))$ on the quantum discs 
$D_{p,q}$ we choose the calculi already used in \cite{cama} and described 
in detail in \cite{vaks}. The 
differential ideal $J(P(D_p) \subset \Omega(P(D_p))$ determining 
$\Gamma(P(D_p))$ is generated by the elements 
\begin{eqnarray*}
x(dx) &-& p^{-1} (dx)x,~~x^*(dx^*) - p(dx^*)x^*, \\
x(dx^*) &-& p^{-1} (dx^*) x,~~x^*(dx) - p(dx) x^*. \end{eqnarray*}
Exchanging $x$ with $y$ and $p$ with $q$ one obtains the differential ideal
$J(P(D_q)) \subset \Omega(P(D_q))$ determining $\Gamma(P(D_q))$. The 
corresponding calculus $\Gamma(P(S^2_{pq\phi}))$ on the basis was explicitely 
described in \cite{cama}.
\\Furthermore we use the right covariant differential calculus 
$\Gamma(P(U(1)))$ determined by the right ideal $R$ generated by the 
element \[ \alpha + \nu \alpha^*-(1+ \nu)1, \] where $0< \nu \leq 1$.
One easily verifies that $R$ fulfills (\ref{jbij}). Thus the differential
ideal $J(P(S^1))$ is generated by the sets (\ref{fish1}), 
(\ref{fish2}) and (\ref{fish3}). Using these generators in the present case
one obtains the following relations in $\Gamma_m(P(S^1))$: \begin{eqnarray*}
(d \alpha^*)\alpha &=& \alpha (d \alpha^*), \\
(d \alpha^*) \alpha &=& \nu \alpha (d \alpha^*), \\
(d \alpha^*) \alpha &=& p\alpha (d \alpha^*), \\
(d \alpha^*) \alpha &=& q \alpha (d \alpha^*).	 
\end{eqnarray*}
Therefore $d \alpha^*=d\alpha=0$ for $\nu,p,q \neq 1$, and the LC
differential algebra $\Gamma_m(P(S^2_{p q \phi}))$ has the following form:
\begin{eqnarray*} \Gamma^0_m( P(S^2_{p q \phi})) &=& P(S^2_{p q \phi}), \\
\Gamma^n_m(P(S^2_{p q \phi})) &=& \Gamma^n(P(D_p)) \bigoplus
\Gamma^n(P(D_q));~~n >0. \end{eqnarray*}
$\Gamma(P(S^2_{p q \phi}))$ coincides with $\Gamma_m(P(S^2_{p q \phi}))$ for $p 
\not= q$, and is embedded as a subspace defined by the gluing for $p=q$ (cf. 
\cite{cama}).
\\Now we want to construct a connection on the bundle ${\cal P}^{(1)}$ which
can be regarded as the connection corresponding to the quantum magnetic 
monopole with strength $g=-\frac{1}{2}$.
\\The  functionals $\cal X$ and $f$ on $P(U(1))$ corresponding to the basis 
element $(\alpha-1) +R \in ker \varepsilon/R$ are given by 
\begin{eqnarray*}
{\cal X}(\alpha) &=& 1;~~{\cal X}(\alpha^*)= -\nu^{-1}, \\
f(\alpha)&=& \nu;~~f(\alpha^*)=\nu^{-1}, \\
f(hk)&=& f(h)f(k);~~h,k \in {\cal A}(U(1)), \\
{\cal X}(hk)&=&{\cal X}(h)f(k)+\varepsilon(h){\cal X}(k). \end{eqnarray*}
$\cal X$ is a linear basis in the space of functionals annihilating $1$
and the right ideal $R$ (see also the appendix and \cite{wor2}), 
i.e. ${\cal X}$ is a basis of the $\nu$-deformed Lie algebra
corresponding to the differential calculus on $U(1)$.
We define the linear maps $A_{l_1}: P(U(1)) \longrightarrow
\Gamma(P(D_p))$ and $A_{l_2} : P(U(1)) \longrightarrow \Gamma(P(D_q))$ 
corresponding to a left connection on ${\cal P}^{(1)}$ by 
\begin{eqnarray} \label{conn}
A_{l_1}(h)&=&{\cal X}(h)\frac{1}{4}(xdx^*-x^*dx)~~h \in P(U(1)), \\
\label{cons} A_{l_2}(h)&=&{\cal X}(h) \frac{1}{4}(y^*dy-ydy^*)~~h \in P(U(1)).
\end{eqnarray} 
Because of ${\cal X}(R)=0$ and ${\cal X}(1)=0$, $A_{l_1}$ and $A_{l_2}$ fulfill
the conditions (\ref{abl1}) and (\ref{kerrl}). Since there is no gluing 
$\Gamma^1_m(B)$ the condition (\ref{abl2}) is also fulfilled.
\\Moreover any choice of one forms to the right of ${\cal X}$ gives a 
connection.
\\\\
A short calculation shows (see formula (\ref{fl})) that the linear maps 
$F_1: P(U(1)) \longrightarrow \Gamma^2(P(D_p))$ and $F_2: P(U(1)) 
\longrightarrow \Gamma^2(P(D_q))$ corresponding to the curvature have the 
following form: \begin{eqnarray*} 
F_1(h)&=&{\cal X}(h)\frac{1}{4}(1+p) dxdx^*+ 
\sum {\cal X}(h_1) {\cal X}(h_2) \frac{1}{16}(xx^* - px^*x)dxdx^*,
\\
F_2(h)&=&-{\cal X}(h)\frac{1}{4}(1+q) dydy^*+ \sum {\cal X}(h_1) 
{\cal X}(h_2) \frac{1}{16}(yy^* - qy^*y)dydy^* \end{eqnarray*}
 \\[.3cm]
In the classical case, the local connection forms ${A_l}_1$ and ${A_l}_2$
can be transformed, using suitable local coordinates, from the classical
unit discs to the upper and lower hemispheres of the classical $S^2$.
The resulting local connection forms on $S^2$ just coincide with the
well-known magnetic potentials of the Dirac monopole of charge $-1/2$.
\\To explain this we will briefly describe the classical Dirac monopole 
(see \cite{naka}).
\\\\The classical Dirac monopole is defined on ${\R}^3-\{0\}$, which is of 
the same homotopy type as $S^2$. The corresponding gauge theory is 
a $U(1)$ theory, and the Dirac monopole is described as a connection on a 
$U(1)$ principal fibre bundle over $S^2$.
\\\\Let $\{U_N,U_S\}$ be a covering of $S^2$, where $U_N$ respectively $U_S$ is 
the closed northern respectively southern hemisphere, $U_N \cap 
U_S =S^1$. One can write $U_N$ and $U_S$ in polar coordinates (up to the poles):
\begin{eqnarray*}
U_N &=& \{(\theta, \phi),~0 < \theta \leq \pi/2,~0 \leq \phi <2 \pi \} \cup 
\{N\},\\
U_S &=& \{(\theta, \phi),~~\pi/2 \leq \theta < \pi,~~0 \leq \phi <2 \pi\} \cup 
\{S\}. \end{eqnarray*} 
By 
\[ g^{(n)}_{12}(\phi)=
\exp(i\, n \phi),~~0\leq \phi < 2\pi,~~n \in \Z  
\]
 a family of 
transition functions $g^{(n)}_{12}: S^1 \longrightarrow U(1),~~n \in \Z $ is given.
A standard procedure defines a corresponding family of $U(1)$ principal fibre 
bundles $Q^{(n)}$.
\\Let $\xi_i : S^1 \longrightarrow U_i,~i=N,S$ be the embedding defined by 
$\xi_i(\phi)=(\pi/2,\phi)$. 
A connection on $Q^{(n)}$ is defined by two Lie algebra valued one forms 
$\tilde{A}_N$ und $\tilde{A}_S$ fulfilling 
\[ \xi^*_N(\tilde{A}_N)=\xi^*_S(\tilde{A}_S)+ i\,n d\phi. \]  
The Wu-Yang forms defined by
\begin{eqnarray*}
\tilde{A}^{(n)}_N &=& i\frac{n}{2}(1-\cos \theta)d\phi, \\
\tilde{A}^{(n)}_S &=& -i\frac{n}{2}(1+ \cos \theta)d \phi \end{eqnarray*} fulfill 
these condition.
$\tilde{A}^{(n)}_N$ and $\tilde{A}^{(n)}_S$ are vector potentials generating 
the magnetic field $B=\frac{n}{2} \frac{\vec{r}}{|\vec{r}|^3}$.
The strength of the Dirac monopole is $n/2$.
\\\\The classical analogue $\tilde{P}^{(n)}$ to the above constructed locally 
trivial QPFB ${\cal P}^{(n)}$ 
is $U(1)$ principal fibre bundles over a space constructed by 
gluing together two discs over their boundaries.
A disc $D$ can be regarded as a subspace of $\C$: 
\[ D:=\{ x \in \C,~xx^* \leq 1 \}.\] The space resulting from the gluing 
together two copies of $D$ over $S^1=\{x \in D,~xx^*=1 \}$ is 
topologically isomorphic to the sphere $S^2$. Every $x \in 
S^1$ has the form $x=\exp(i \phi),~0 \leq \phi < 2 \pi$. The classical 
$U(1)$ bundles $\tilde{P}^{(n)}$ are given by transition functions 
$\tilde{g}^{(n)}_{12} : S^1 \longrightarrow U(1)$, which are obtained by 
$\tau^{(n)}_{12}=(\tilde{g}^{(n)}_{21})^*$ ($*$ means here the pull-back) 
from the above transition functions of QPFB. The exchange of the indices comes 
from formula (\ref{phi}.
One has $\tilde{g}^{(n)}_{12}(\exp(i \phi))=\exp (-i\,n \phi),~n \in \N$. 
Obviously, the $\tilde{P}^{(n)}$ are topologically isomorphic to $Q^{(-n)}$.
\\The classical analogue to the above defined connection on ${\cal P}^{(1)}$
is given by the following one forms on $D$ (see (\ref{conn}) and (\ref{cons})):
\begin{eqnarray*} \tilde{A}_1 &=& \frac{1}{4}(xdx^*-x^*dx), \\
\tilde{A}_2 &=& \frac{1}{4}(x^*dx - xdx^*). \end{eqnarray*}
Let $\xi : S^1 \longrightarrow D$ the embedding. A short calculation shows
that $\tilde{A}_1$ and $\tilde{A}_2$ fulfill 
\[ \xi^*(\tilde{A}_1)= \xi^*(\tilde{A}_2)-i d\phi. \] 
Now one defines the following maps  
$\eta_N: U_N \setminus \{N\}  \longrightarrow D$ and $\eta_S: U_S \setminus 
\{S\} \longrightarrow D$ by
\begin{eqnarray*} \eta_N(\theta, \phi) &:=& \sqrt{1-\cos \theta} \exp{(i \phi)}, 
\\ \eta_S(\theta, \phi) &:=& \sqrt{1+ \cos \theta} \exp{(i \phi)}, 
\end{eqnarray*} 
and one easily verifies \begin{eqnarray*}
\tilde{A}^{-1}_N &=& \eta^*_N(\tilde{A}_1), \\
\tilde{A}^{-1}_S &=& \eta^*_S(\tilde{A}_2). \end{eqnarray*}
 
\section{Final remarks}
We have developed the general scheme of a theory of connections on locally 
trivial QPFB, including a reconstruction theorem for bundles and a nice
characterization of connections in terms of local connection forms.
Here we make some remarks about questions and problems arising in our context,
and about possible future developments.\\[.2cm]
1. It is very important to look for more examples. Our example of a $U(1)$
bundle over a glued quantum sphere is essentially the same as the example of
\cite{Kon} of an $SU_q(2)$ bundle over an anlogous glued quantum sphere.
(Indeed, in \cite{Kon} another quantum disc is used, which, however,
is isomorphic to the disc used in our paper -- both are isomorphic to the
shift algebra. The bundles are equivalent in the sense of the main Theorem
of \cite{Kon}, which says that a QPFB with structure group $H$ is determined
by a bundle with the classical subgroup of $H$ as structure group.)
For other examples, one has to look for algebras with a covering (or being
a gluing) such that the $B_{ij}$ are ``big enough'' to allow for nontrivial
transition functions $\tau_{ij}:H\lra B_{ij}$: $B_{ij}$ must contain in its
center subalgebras being the homomorphic image of the algebra $H$.
This seems to be possible only if $H$ has nontrivial classical subgroups
and $B_{ij}$ contains suitable classical subspaces, as in our example.
The following (almost trivial) example of a gluing along two noncommutative
parts indicates that one may fall back to a gluing along classical subspaces
in many cases:
Let $A_1=C^*({\frak S})\oplus_\sigma C^*({\frak S})=A_2$ be two copies of a
quantum sphere being glued together from shift algebras via the symbol map
$\sigma$, as described in \cite{cama}. Then the gluing
$A_1\oplus_{pr_{1,2}}A_2:=\{((a_1,a_2),(a_1',a_2'))\in A_1\oplus A_2~|~ 
a_2=a_1'\}$ (gluing of two quantum spheres along hemispheres) is obviously
isomorphic to $\{(a_1,a_2,a_3)\in C^*({\frak S})\oplus C^*({\frak S})
\oplus C^*({\frak S})~|~ \sigma(a_1)=\sigma(a_2)=\sigma(a_3)\}=
C^*({\frak S})\oplus_\sigma C^*({\frak S})\oplus_\sigma C^*({\frak S})$. 
This is a glued quantum sphere with a (quantum disc) membrane inside,
glued along the classical subspaces. (This corresponds perfectly to the
classical picture of gluing two spheres along hemispheres.)\\[.2cm]
2. The permanent need to work with covering completions is an unpleasant
feature of the theory. It would therefore be very important to find some
analogue of algebras of smooth functions in the noncommutative situation which 
have a suitable class of ideals forming a distributive lattice with respect to
$+$ and $\cap$ (cf. \cite[Proposition 2]{cama}). It is not clear if such a 
class exists even in classical algebras of differentiable functions.\\[.5cm]
3. Principal bundles are in the classical case of utmost importance
in topology and geometry. In the above approach, one could e. g. ask for
characteristic classes (trying to generalize the Chern-Weil construction),
and for a notion of parallel transport defined by a connection 
(a naive idea would be to call a horizontal form parallel, if its covariant
derivative vanishes).\\[.5cm]
4. For locally trivial QPFB, a suitable notion of locally trivial associated 
quantum vector bundle (QVB) exists (\cite{cama2}). Its definition (via
cotensor products) is designed to have the usual correspondence between
vector valued horizontal forms (of a certain ``type'') and sections of the
associated bundle. To a connection on a QPFB one can also associate 
connections on the corresponding QVB (assuming a certain differential structure
there).\\[.5cm]
5. The notion of gauge transformation in our context is considered in
\cite{cama3}. Gauge transformations are defined as isomorphisms of the left 
(right) $B$-module ${\cal P}$, with natural compatibility conditions. It turns out
that the set of covariant derivatives is invariant under gauge transformations,
whereas connections are not always transformed into connections.
 
\section{Appendix}
The purpose of this appendix is to collect some results about covariant
differential calculi on quantum groups (\cite{dipl}, \cite{wor2},
\cite{klsch})
and about coverings and gluings of algebras and differential algebras 
\cite{cama}.
\subsection{Covariant calculi on Hopf algebras}
We freely use standard facts about Hopf algebras, including the Sweedler
notation (e. g. $\Delta(h)=h_1\otimes h_2$). We assume that the antipode
is invertible.

A differential algebra over an algebra $B$ is a $ \N $-graded algebra
$\Gamma(B)=\oplus_{i\in \N }\Gamma^i(B)$, $\Gamma^0((B)=B$, equipped
with a differential $d$, i. e. a graded derivative of degree $1$ with 
$d^2=0$. It is called differential calculus if it is generated as an algebra
by the $db$, $b\in B$.
A differential ideal of a differential algebra is a $d$-invariant graded  
ideal. 
There is always the universal differential calculus
$\Omega(B)$ determined by the property that every differential calculus
$\Gamma(B)$ is of the form $\Gamma(B)\simeq\Omega(B)/J(B)$ for some 
differential ideal $J(B)$.

If two algebras $A$, $B$ and differential algebras $\Gamma(A)$, $\Gamma(B)$
are given, an algebra homomorphism $\psi:A\lra B$ is said to be 
differentiable
with respect to $\Gamma(A)$, $\Gamma(B)$, if there exists a homomorphism 
$\psi_\Gamma:\Gamma(A)\lra\Gamma(B)$ of
differential algebras extending $\psi$. For $\Gamma(A)=\Omega(A)$ this
extension, denoted in this case by $\psi_{\Omega \rightarrow \Gamma}$, always 
exists. If, in addition, $\Gamma(B)=\Omega(B)$, the notation $\psi_{\Omega}$
is used. $J(B)=ker id_{\Omega \rightarrow \Gamma}$ is a differential ideal $J(B) 
\subset \Omega(B)$ such that $\Gamma(B)=Omega(B)/J(B)$. $J(B)$ is called the 
differential ideal corresponding to $\Gamma(B)$.

Now we list some facts about covariant differential calculi.   
\bde A differential calculus $\Gamma(H)$ over a 
Hopf algebra $H$ is called right covariant, if $\Gamma(H)$ is a right $H$ 
comodule algebra with right coaction $\Delta^{\Gamma}$ such that 
\begin{equation} 
\Delta^{\Gamma}(h_0dh_1 ...dh_n)=\Delta(h_0)(d 
\otimes id)\circ \Delta(h_1)...(d \otimes id) \circ \Delta(h_n). 
\end{equation} 
$\Gamma(H)$ is called left covariant, if $\Gamma(H)$ is a left $H$-comodule
algebra with left coaction ${}^{\Gamma}\Delta$ such that 
\begin{equation} 
{}^{\Gamma}\Delta(h_0dh_1 ...dh_n)=\Delta(h_0)(id 
\otimes d)\circ \Delta(h_1)...(id \otimes d) \circ \Delta(h_n).
\end{equation}
$\Gamma(H)$ is called bicovariant if it is left and right 
covariant. \ede
 Because of the universality property the universal differential calculus
over any 
 Hopf algebra is bicovariant. In the sequel we list some properties of
right covariant differential calculi. The construction of left covariant 
 differential algebras is analogous.
 \\
 Let $\Delta^{\Omega}$ be the right coaction of the 
 universal differential calculus $\Omega(H)$ and let $\Gamma(H)$ be a 
 differential algebra over the Hopf algebra $H$.
 $\Gamma(H)$ is right covariant if and only if the corresponding
differential 
 ideal $J(H) \subset \Omega(H)$ has the property \[ \Delta^{\Omega}(J(H)) 
 \subset J(H) \otimes H. \]
 Let us consider a right-covariant differential calculus $\Gamma(H)$. Let 
 $\Gamma^1_{inv}(H):=\{\gamma \in \Gamma(H)| \Delta^{\Gamma}(\gamma)
=\gamma 
 \otimes 1\}$. There exists a projection $P: \Gamma^1(H) \longrightarrow 
 \Gamma^1_{inv}(H)$ defined by \[ P(\sum
h^0dh^1)=S^{-1}(h^0_2h^1_2)h^0_1dh^1_1. 
 \] Now one can define a linear map $\eta_{\Gamma} : H \longrightarrow
\Gamma^1(H)$ by 
 \[ \eta_{\Gamma}(h):=P(dh)=\sum S^{-1}(h_2)dh_1. \] By an easy calculation
one obtains 
 the identity $dh=\sum h_2\eta_{\Gamma}(h_1)$. 
 The linear map $\eta_{\Gamma}$ has the following properties: 
\begin{eqnarray*}\Delta^{\Gamma}(\eta_{\Gamma}(h))&=&\eta_{\Gamma}(h)
\otimes 1, \\
\eta_{\Gamma}(h)k&=&\sum k_2(\eta_{\Gamma}(hk_1)-\varepsilon
(h)\eta_{\Gamma}(k_1)), \\
d\eta_{\Gamma}(h) &=& -\sum \eta_{\Gamma}(h_2)\eta_{\Gamma}(h_1)
\end{eqnarray*} 
In the case $\Gamma(H)=\Omega(H)$ we use the symbol $\eta_\Omega$.
\\The first degrees of 
right-covariant differential algebras are in one-to-one correspondence to 
right ideals $R \subset ker \varepsilon \subset H$ in the following sense:
First, if a differential calculus is given, $R:=ker \eta_\Gamma\cap ker 
\varepsilon$ is a right ideal with the property $R \subset 
 ker\varepsilon$, and one can prove that the subbimodule $J^1(H)$
corresponding 
 to  $\Gamma^1(H) \cong  \Omega^1(H)/J^1(H)$ is generated by the space 
 $\eta_\Omega(R)=\{\sum S^{-1}(r_2)dr_1|r \in R \}$. On the other hand, 
every right ideal
$R \subset ker \varepsilon$ defines a right covariant differential algebra 
 $\Gamma(H)=\Omega(H)/J(H)$, where the 
 differential ideal $J(H) \subset \Omega(H)$ is generated by the set 
$\eta_\Omega(R)$. Analogously, right ideals $R \subset ker 
 \varepsilon$ also correspond to left covariant differential calculi.
In this case, the
 differential ideal $J(H)$ corresponding to $R$ is generated by 
$\{\sum S(r_1)dr_2|r \in R \}$. Bicovariant differential 
 calculi are given by right ideals $R$ with the property $\sum S(r_1)r_3  
\otimes r_2 \subset H \otimes 
 R;~~\forall r \in R$ (Ad-invariance).\\[.3cm]
Now one can choose a linear basis $h_i +R$ in $ker \varepsilon /R$. This
leads to 
a set of functionals ${\cal X}_i$ on $H$ annihilating $1$ and $R$ such
that $\eta_{\Gamma}(h)={\cal X}_i(h) \eta_{\Gamma}(h_i)$, $h\in H$. The 
set of elements $\eta_{\Gamma}(h_i)$ is a left and right $H$ module basis
in 
$\Gamma^1(H)$, and the set of the ${\cal X}_i$ is a linear basis in the 
space of all 
functionals annihilating $1$ and $R$. It is obvious that $dh=\sum h_2 
{\cal X}_i(h_1) \eta_{\Gamma}(h_i)$. 
Besides the functionals ${\cal X}_i$ the linear basis in $ker \varepsilon/R$ 
determines also functionals
$f_{ij}$ on $H$ satisfying 
\begin{eqnarray*} f_{ij}(1) &=& \delta_{ij}, \\
f_{ij}(hk)&=& \sum_l f_{il}(h) f_{lj}(k) 
\end{eqnarray*} 
\[ {\cal X}_i(hk)=\sum_l {\cal X}_l(h) f_{li}(k) + \varepsilon(h) {\cal
X}_i(k).\]

\begin{de} 
Let $A$ be a vector space and let $H$ be a Hopf algebra such
that there exists linear map $\Delta_A : A \longrightarrow A \otimes H$.
$\Delta_A$ is called right $H$-coaction and $A$ is called right $H$ comodule if
\begin{eqnarray} (\Delta_A \otimes id) \circ \Delta_A &=& (id \otimes \Delta)
\circ \Delta_A ,\\
(id \otimes \varepsilon) \circ \Delta_A &=& id. 
\end{eqnarray} 
\end{de}
If $A$ is an algebra and $\Delta_A$ is an homomorphism of algebras then $A$ is
called a right $H$ comodule algebra. The left coaction is defined analogously.
\\The definition of covariant differential calculi over Hopf algebras is easily 
generalized to $H$ comodule algebras:

\begin{de} \label{mocov}
 A differential calculus $\Gamma(A)$ over a right $H$ comodule 
algebra $A$ is called right covariant if the right coaction $\Delta^{\Gamma}_A : 
\Gamma(A) \longrightarrow \Gamma(A) \otimes H$ defined by 
\begin{equation}
\Delta^{\Gamma}_A (a_0da_1 ...da_n)=\Delta_A(a_0)(d
\otimes id)\circ \Delta_A(a_1)...(d \otimes id) \circ \Delta_A(a_n)
\end{equation} exists. 
\end{de}

\subsection{Covering and gluing}
Let finite families $(B_i)_{i\in I}$, $(B_{ij})_{(i,j)\in I\times I\setminus
D}$, $D$ the diagonal in $I\times I$, $B_{ij}=B_{ji}$, and homomorphisms
$\pi^i_j:B_i\lra B_{ij}$ be given. Then the algebra
\[
\textstyle
B=\{(b_i)_{i\in I}\in \bigoplus_i
B_i|\pi^i_j(b_i)=\pi^j_i(b_j)~\forall
 i\neq j\}=:\bigoplus_{\pi^i_j}B_i
\]
is called gluing of the $B_i$ along the $B_{ij}$ by means of the $\pi^i_j$.
Special cases of gluings arise from coverings:
A finite covering of an algebra $B$ is a finite family $(J_i)_{i\in I}$
of ideals in $B$ with $\cap_iJ_i=0$. Taking now $B_i=B/J_i$, $B_{ij}=
B/(J_i+J_j)$, $\pi^i_j:B_i\lra B_{ij}$ the canonical projections
$b+J_i\mapsto b+J_i+J_j$, one can form the gluing
\[
\textstyle
B_c=\bigoplus_{\pi^i_j}B_i,
\] 
which is called the covering completion
of $B$ with respect to the covering $(J_i)_{i \in I}$.
$B$ is always embedded in $B_c$ via the map $K:b\mapsto (b+J_i)_{i\in I}$. 
The covering $(J_i)_{i\in I}$ is called complete if $K$ is also surjective, 
i.e. $B$ is isomorphic to $B_c$. Every two-element covering is complete, as 
well as every covering of a C*-algebra. On the other hand, if 
$B=\bigoplus_{\pi^i_j}
B_i$ is a general gluing, and $p_i:B\lra B_i$ are the restrictions of the 
canonical projections, then $(ker~p_i)_{i\in I}$ is a complete covering of 
$B$.
					  
If $\Gamma(B)$ is a differential algebra, a covering $(J_i)_{i\in I}$ of
$\Gamma(B)$ is said to be differentiable if the $J_i$ are differential 
ideals.
A differential algebra $\Gamma(B)$ with differentiable covering 
$(J_i)_{i\in I}$ is called LC differential algebra (LC = locally calculus),
if the factor differential algebras $\Gamma(B)/J_i$ are differential calculi
 over $B/J_i^0$ ($J_i^0$ the degree zero component of $J_i$) and $J_i^0\neq
0,~\forall i$.

\begin{de} \label{auer} Let $(B,(J_i)_{i\in I})$ be an algebra with complete covering,
let $B_i=B/J_i$, let $\pi_i:B\lra B_i$ be the natural surjections, and
let $\Gamma(B)$
and $\Gamma(B_i)$ be differential calculi such that $\pi_i$ are 
differentiable and $(ker{\pi_i}_\Gamma)_{i\in I}$ is a covering
of $\Gamma(B)$. Then $(\Gamma(B),(\Gamma(B_i))_{i\in I})$ is called adapted
to $(B,(J_i)_{i\in I})$. \end{de}
The following proposition is essential for Definition \ref{strukstrik}:

\begin{pr} \label{diff}
Let $(B,(J_i)_{i \in I})$ be an algebra with
complete covering, and let $\Gamma(B_i)$ be differential calculi over the
algebras $B_i$. Up to isomorphy there exists a unique differential
calculus $\Gamma(B)$ such that $(\Gamma(B),(\Gamma(B_i))_{i \in I}))$ is  
adapted to $(B,(J_i)_{i \in I})$. 
\end{pr}
As shown in \cite{cama}, the differential ideal corresponding to
$\Gamma(B)=\Omega(B)/J(B)$ is just $J(B)=\cap_{i\in 
I}ker\,{\pi_i}_{\Omega
\rightarrow\Gamma}$.\\
Finally, there is a proposition concerning the covering completion of 
adapted 
differential calculi:

\begin{pr} \label{diff1}
Let $( \Gamma(B),(\Gamma(B_i))_{i \in I}) $ be
adapted to $(B,(J_i)_{i \in I})$.
Then the covering completion of $( \Gamma(B), (ker \pi_{i_\Gamma})_{i \in
I})$ is an LC differential algebra over $B_c$. 
\end{pr}


\begin{thebibliography}{100}
\bibitem{br} Brzezinski, T.: Translation map in quantum principal bundles,
{\it J. Geom. and Phys.} {\bf 20} (1996), 349-370
\bibitem{brma} Brzezi\'nski, T., and S. Majid: Quantum group gauge theory
on quantum spaces, {\it Commun. Math. Phys.} {\bf 157} (1993), 591--638,
{\sf hep-th/9208007}, Preprint DAMTP/92-27
\bibitem{Kon} Budzy\'nski, R. J. and W. Kondracki: Quantum principal fiber
bundles: Topological aspects, {\it Rep. Math. Phys.} {\bf 37} (1996), 365--385,
preprint 517 PAN Warsaw 1993, {\sf hep-th/9401019}
\bibitem{dipl} Calow, D., Differentialkalk\"ule auf Quantengruppen,
Diplomarbeit,
Leipzig 1995
\bibitem{cama} Calow, D. and R. Matthes: Covering and gluing of algebras
and differential algebras, {it J. Geom. and Phys.} {\bf 32} (2000),
364--396,
{\sf math.QA/9910031}, Preprint NTZ 25/1998
\bibitem{cama2} Calow, D. and R. Matthes: Locally trivial quantum vector bundles,
{\sf math.QA/00}
\bibitem{cama3} Calow, D. and R. Matthes: Gauge transformations on locally trivial
quantum principal bundles, {\sf math.QA/00}
\bibitem{di1} Dixmier, J.: {\it Les $C^*$-algebres et leurs representations},
Gauthier-Villars, Paris 1964
\bibitem{dop} Doplicher, S.: Quantum spacetime, {\it Ann. Inst. Henri
Poincare, Physique theorique} {\bf 64} (1996), 543--553
\bibitem{dfr} Doplicher, S., Fredenhagen, K. and J. E. Roberts: The quantum
structure of spacetime at the Planck scale and quantum fields,
{\it Commun. Math. Phys.} {\bf 172} (1995), 187--220
\bibitem{dur} Durdevic, M.: Geometry of quantum principal bundles I,
{\it Commun. Math. Phys.} {\bf 175} (1996), 457-521, {\sf q-alg/9507019}
\bibitem{dur1} Durdevic, M.: Geometry of quantum principal bundles II,
{\it Rev. Math. Phys.} {\bf 9} (5) (1997), 531--607, {\sf q-alg/9412005}
\bibitem{fgr} Fr\"ohlich, J., Grandjean, O. und A. Recknagel: Supersymmetric 
quantum theory, non-commutative geometry, and gravitation, {\it Sym\'etries
quantiques (Les Houches, 1995)}, 221--385, North Holland, Amsterdam, 1998,
ETH-TH/97-19, 
{\sf hep-th/9706132}
\bibitem{haj} Hajac, P. M.: Strong connections on quantum principal
bundles,
{\it Commun. Math. Phys.} {\bf 182} (1996), 579-617
\bibitem{kem4} Kempf, A.: String/quantum gravity motivated uncertainty
relations and regularisation in field theory, {\sf hep-th/9612082},
DAMTP/96-101
\bibitem{klile2} Klimek, S. and A. Lesniewski: A two-parameter quantum 
deformation of the unit disc, {\it J. Funct. Anal.} {\bf 115} (1993), 1--23
\bibitem{klsch} Klimyk, A. U. and K. Schm\"udgen: {\it Quantum
groups and their representations}, Texts and Monographs in Physics, Springer
1997
\bibitem{korn} Koornwinder, T. H.: General Compact Quantum Groups, a
Tutorial, Preprint
University of Amsterdam, Faculty of Mathematics and Computer Science
\bibitem{mue} M\"uller, A.: Classifying spaces for quantum principal
bundles, {\it Commun. Math. Phys.} {\bf 149} (1992), 495--512
\bibitem{naka} Nakahara M.: Geometry, Topology and Physics, Graduate Student
Series in Physics, Institute of Physics Publishing, Bristol and Philadelphia,
1990
\bibitem{pf1} Pflaum, M. J.: Quantum groups on fibre bundles,
{\it Comm. Math. Phys.} {\bf 166} (1994), 279-315, {\sf hep-th/9401085}
\bibitem{pfsch} Pflaum, M. J. and P. Schauenburg: Differential calculi on
noncommutative bundles, {\it Z. Phys. C} {\bf 6} (1997), 733--744,
{\sf q-alg/9612030}, Preprint gk-mp-9407/7 M\"unchen 1994
\bibitem{scha}
Schauenburg, P.: Zur nichtkommutativen Differentialgeometrie von
Hauptfaserb\"undeln - Hopf-Galois-Erweiterungen von De Rham-Komplexen,
Dissertation M\"unchen 1993
\bibitem{schn} Schneider, H. J.: Principal
homogeneous spaces for arbitrary Hopf algebras, {\it Isr. J. Math.} {\bf 72}
(1990), 167--195
\bibitem{vaks} Sinel'shchikov, S. and L. Vaksman: On
q-analogues of bounded symmetric
domains and Dolbeault Complexes, {\it Mathematical Physics, Analysis and
Geometry,}
{\bf 1} (1)(1998), 75--100, {\sf q-alg/9703005}
\bibitem{wor1} Woronowicz, S. L.: Twisted SU(2) Group. An Example of a
Non-Commutative Differential Calculus, {\it Publ. RIMS, Kyoto University}
{\bf 23} (1987), 117--181
\bibitem{wor3} Woronowicz, S. L.: Compact matrix pseudogroups,
{\it Commun.
Math. Phys.} {\bf 111} (1987), 613--665
\bibitem{wor2} Woronowicz, S. L.: Differential calculus on compact
matrix pseudogroups (quantum groups), {\it Commun. Math. Phys.} {\bf 122}
(1989), 125--170
\end{thebibliography}
\end{document}